\documentclass[10pt]{article}

\usepackage{amsmath,amssymb,epsfig}
\usepackage{theorem,array,latexsym}
\usepackage[latin1]{inputenc}
\usepackage{color}
\usepackage[active]{srcltx}

\newcommand{\D}{\protect\displaystyle}
\newcommand{\T}{\protect\textstyle}
\newcommand{\ipl}{\langle} 
\newcommand{\ipr}{\rangle} 

\newtheorem{theor}{Theorem}

\newtheorem{propo}[theor]{Proposition}

\begin{document}

\title{On inverse problems for semiconductor equations}

\setcounter{footnote}{3}
\author{\normalsize\renewcommand{\thefootnote}{\arabic{footnote}}
M.~Burger$^1$      \quad
H.W.~Engl$^{1,2}$  \quad
A.~Leit\~ao$^2$    \quad
P.A.~Markowich$^{2,3}$}

\date{\small\today}
\maketitle

\vskip-0.8cm \noindent
\begin{center} \begin{minipage}{12cm} \begin{small}
$^1$ Institut f\"ur Industriemathematik, Johannes Kepler Universit\"at,
A--4040 Linz, Austria \ {\tt (engl@indmath.uni-linz.ac.at)} \\[1ex]
$^2$ Johann Radon Institute for Computational and Applied Mathematics, 
Austrian Academy of Sciences c/o Johannes Kepler Universit\"at, 
A-4040 Linz, Austria \ {\tt (antonio.leitao@oeaw.ac.at)} \\[1ex]
$^3$ Institut f\"ur Mathematik, Universit\"at Wien, Boltzmanngasse 9,
A-1090 Vienna, Austria {\tt (peter.markowich@univie.ac.at)}
\end{small} \end{minipage} \end{center}

\medskip
\begin{abstract}
This paper is devoted to the investigation of inverse problems related
to stationary drift-diffusion equations modeling semiconductor devices.
In this context we analyze several identification problems corresponding
to different types of measurements, where the parameter to be reconstructed
is an inhomogeneity in the PDE model (doping profile). \\
For a particular type of measurement (related to the voltage-current map)
we consider special cases of drift-diffusion equations, where the
inverse problems reduces to a classical {\em inverse conductivity problem}.
A numerical experiment is presented for one of these special situations
(linearized unipolar case).
\end{abstract}

%---------------------------------------------------------------------------
\section{Introduction}

The drift diffusion equations are the most widely used model to describe
semiconductor devices. From the point of view of applications, there
is great interest in replacing laboratory testing by numerical simulation
in order to minimize development costs. For the current state of technology,
the drift diffusion equations represent a realistic compromise between
computational efficiency (to solve this nonlinear system of partial
differential equations) and an accurate description of the underlying
device physics.

The name {\em drift diffusion equations} of semiconductors originates
from the type of dependence of the current densities on the carrier
densities and the electric field. The current densities are the sums of
drift terms and diffusion terms.
It is worth mentioning that, with the increased miniaturization of
semiconductor devices, one comes closer and closer to the limits of
validity of the drift diffusion equation. This is due to 
the fact that in ever smaller devices the assumption that the
free carriers can be modeled as a continuum becomes invalid. On the
other hand, the drift diffusion equations are derived by a scaling limit
process, where the mean free path of a particle tends to zero.

The mathematical modeling of semiconductor equations has developed
significantly, together with their manufacturing. The {\em basic
semiconductor device equations} where first presented, in the level
of completeness described in this paper, by W. Van Roosbroeck (see
\cite{VanR}) in 1950. Since then they have been subject of intensive
mathematical  and numerical investigation (cf. \cite{MRS} for an overview).

This paper is devoted to the investigation of inverse problems related
to stationary drift-diffusion equations modeling semiconductor devices.
In this context we analyze several inverse problems related to the
identification of doping profiles. In all these inverse problems the
parameter to be identified corresponds to the so called {\em doping profile}
(a parameter function in a system of PDE's). However, the reconstruction
problems are related to data generated by different types of measurement
techniques.

The paper is organized as follows. In section~\ref{sec:sceq} we describe
the stationary and transient drift diffusion equations. Some existence
and uniqueness results (needed further in the text) are presented and
some particular models (derived from different simplification assumptions)
are investigated. In Section~\ref{sec:idp} the inverse doping problems are
presented. We address the inverse problems modeled by the voltage-current
map, by capacitance measurements, and by Laser-beam-induced measurements.
We also address the identification issue for some of the inverse problems
mentioned above. In Section~\ref{sec:num} we present some new numerical
results for an identification problem related to the voltage-current map
(linearized unipolar case). The results are obtained using the
Landweber-Kaczmarz method.

There are other relevant inverse problems for semiconductor equations that
are not covered in this paper:

\begin{itemize}

\item The inverse problem of identifying transistor contact resistivity of planar
electronic devices, such as MOSFETs (metal oxide semiconductor field-effect
transistors) is treated in \cite{FC}. It is shown that a one-point boundary
measurement of the potential is sufficient to identify the resistivity from
a one-parameter monotone family, and such identification is both stable and
continuously dependent on the parameter. Because of the device
miniaturization, it is impossible to measure the contact resistivity in a direct way
to satisfactory accuracy. There are extensive experimental and simulation
studies for the determination of contact resistivity by certain accessible
boundary measurements.

\item A similar problem of determining the contact resistivity of a semiconductor device
from a single voltage measurement is investigated in \cite{BF}. It can be modeled
as an inverse problem for the elliptic differential equation $\Delta V -
p \chi(S)u = 0$ in $\Omega \subset \mathcal R^2$, $\partial V / \partial n
= g \geq 0$ but $g \not\equiv 0$ on $\partial\Omega$, where $V(x)$ is the
measured voltage, $S \subset \Omega$ and $p>0$ are unknown. In this paper,
the authors consider the identification of $p$ when the contact location $S$
is also known.

\item  The problem of optimal design of devices, where the aim is to find a
doping profile that can reach certain design goals, e.g., maximum drive
current while keeping the leakage current below a certain threshold. From a
computational point of view, this problem exhibits many similarities to the
inverse doping problems considered in this paper, with the difference that in
optimal design one usually has to solve the drift-diffusion equations for only
one or two different applied voltages. We refer to \cite{BuPi,HiPi01,St}
\end{itemize}

%---------------------------------------------------------------------------
\section{Semiconductor equations} \label{sec:sceq}

The basic semiconductor
device equations consist of the Poisson equation (\ref{eq:VanR1}),
the continuity equations for electrons (\ref{eq:VanR2}) and holes
(\ref{eq:VanR3}), and the current relations for electrons
(\ref{eq:VanR4}) and holes (\ref{eq:VanR5}). For some applications,
in order to account for thermal effects in semiconductor devices, its
also necessary to add to this system the heat flow equation (\ref{eq:VanR6}).
\begin{eqnarray}
&& {\rm div} (\epsilon \nabla V)  =  q(n - p - C)        \label{eq:VanR1} \\
&& {\rm div}\, J_n  =  q ( \partial_t n + R)             \label{eq:VanR2} \\
&& {\rm div}\, J_p  =  q (-\partial_t p - R)             \label{eq:VanR3} \\
&& J_n  =  q ( D_n(E,T) \nabla n - \mu_n(E,T) n \nabla V)          \label{eq:VanR4} \\
&& J_p  =  q (-D_p(E,T) \nabla p - \mu_p(E,T) p \nabla V)          \label{eq:VanR5} \\
&& \rho\ c(T)\ \partial_t T - H = {\rm div}\,k(T) \nabla T, \label{eq:VanR6}
\end{eqnarray}
defined in a domain $\Omega \subset \mathbb{R}^d$ ($d=1,2,3$) representing
the semiconductor device. Here $V$ denotes the electrostatic potential
($- \nabla V$ is the electric field, $E=\vert\nabla V\vert$), $n$ and $p$ are the concentration
of free carriers of negative charge (electrons) and positive charge (holes)
respectively, and $J_n$ and $J_p$ are the densities of the electron and the
hole current respectively. $D_n$ and $D_p$ are the diffusion coefficients
for electrons and holes respectively. $\mu_n$ and $\mu_p$ represent the
mobilities of electrons and holes respectively. The positive constants
$\epsilon$ and $q$ denote the permittivity coefficient (for silicon)
and the elementary charge respectively. The function $R = R(n,p,x)$ denotes
the recombination-generation rate. The function $C = C(x)$ represent the
doping concentration, which is produced by diffusion of different materials
into the silicon crystal and by implantation with an ion beam.
The constants $\rho$ and $c$ represent the specific mass density and
specific heat of the material. $k$ and $H$ denote the thermal
conductivity and the locally generated heat.

This set of equations is considered in a domain $\Omega \subset
\mathbb{R}^d$ ($d=1,2,3$) representing the semiconductor device. We
assume the boundary $\partial\Omega$ of $\Omega$ to be divided into two
nonempty disjoint parts: $\partial\Omega = \overline{\partial\Omega_N} \cup
\overline{\partial\Omega_D}$. The Dirichlet part of the boundary $\partial\Omega_D$
models the Ohmic contacts, where the potential $V$ as well as the
concentrations $n$ and $p$ are prescribed. The Neumann part
$\partial\Omega_N$ of the boundary corresponds to insulating surfaces,
thus a zero current flow and a zero electric field in the normal
direction are prescribed.

In the next subsection, when we turn our attention to the stationary
drift-diffusion equations, we shall discuss in more detail the system
(\ref{eq:VanR1})--(\ref{eq:VanR6}) as well as corresponding boundary
conditions. Detailed expositions of the subject of modeling,
analysis and simulation of semiconductor equations can be found in the
books of S. Selberherr \cite{Se}, P. Markowich \cite{Ma} and P. Markowich
et al \cite{MRS}.

\subsection{Stationary drift diffusion equation}

We shall consider system (\ref{eq:VanR1}) -- (\ref{eq:VanR6}) under the
following assumptions: Thermal effects will not be taken into account, i.e. we
shall work under the assumption of constant particle temperature. Further, we
consider the carrier concentrations $n$ and $p$ and the potential $V$
as time-independent functions.

Under the above assumptions, if we substitute the current relations for
electrons and holes (\ref{eq:VanR4}) and (\ref{eq:VanR5}) into the
corresponding continuity equations (\ref{eq:VanR2}) and (\ref{eq:VanR3}),
we obtain a coupled system of partial differential equations, the so
called {\em stationary drift diffusion equation}:
\begin{eqnarray}
 {\rm div} (\epsilon_s \nabla V) & = & q(n - p - C),\ {\rm in}\ \Omega
 \label{eq:stat-dd-1} \\
 {\rm div} (D_n \nabla n - \mu_n n \nabla V) & = & R,\ {\rm in}\ \Omega
 \label{eq:stat-dd-2} \\
 {\rm div} (D_p \nabla p - \mu_p p \nabla V) & = & R,\ {\rm in}\ \Omega .
 \label{eq:stat-dd-3}
\end{eqnarray}

Next we briefly discuss the modeling of the recombination-generation rate.
The {\em bandgap} is relatively large for semiconductors (gap between
valence and conduction band), and a significant amount of energy is
necessary to transfer electrons from the valence and to the conduction
band. This process is called generation of electron-hole pairs. On the
other hand, the reverse process corresponds to the transfer of a
conduction electron into the lower energetic valence band. This process
is called recombination of electron-hole pairs. In our model these
phenomena are described by the {\em recombination-generation rate} $R$. Various
models can be found in the literature (see, e.g., \cite{Se}). For the sake
of simplicity, we shall consider either the {\em Shockley Read Hall rate}
$$  R_{SRH} \ = \ \D\frac{n p - n_i^2}{\tau_p(n+n_i) + \tau_p(p+n_i)} $$
or the {\em Auger recombination-generation rate}
$$    R_{AU} \ = \ (C_n n + C_p p) \, (n p - n_i^2) \, , $$
where $n_i$ denotes the intrinsic density, $\tau_n$ and $\tau_p$ are the
lifetimes of electrons and holes respectively (see Table~\ref{tab:typ-val}
for some typical values for recombination generation parameters).
In both cases we can write
$$ R \ = \ {\cal R}(n,p,x) \, (n p - n_i^2) \, . $$

\begin{table}
\begin{center} \begin{tabular}{cl}
\hline \\[-2.3ex]
{ Parameter} & { \hfil Typical value \hfil} \\[0.3ex]
\hline \\[-2.0ex]
$\epsilon_s$ & $11.9 \ \epsilon_0$ \\
$\mu_n$      & $\approx 1500 \ {\rm cm}^2 \ {\rm V}^{-1} \ {\rm s}^{-1}$ \\
$\mu_p$      & $\approx  450 \ {\rm cm}^2 \ {\rm V}^{-1} \ {\rm s}^{-1}$ \\
$C_n$    & $2.8 \times 10^{-31} \ {\rm cm}^6 {\rm / s}$ \\
$C_p$    & $9.9 \times 10^{-32} \ {\rm cm}^6 {\rm / s}$ \\
$\tau_n$ & $10^{-6}\, {\rm s}$ \\
$\tau_p$ & $10^{-5}\, {\rm s}$ \\
\hline
\end{tabular} \end{center}
\caption{Properties of silicon at room temperature (physical constants:
Permittivity in vacuum $\epsilon_0 = 8.85 \times 10^{-14} {\rm As \, V}^{-1}
\, {\rm cm}^{-1}$; Elementary charge $q = 1.6 \times 10^{-19} {\rm As}$).}
\label{tab:typ-val}
\end{table}

Now we shall introduce the boundary conditions. As already mentioned,
the boundary $\partial\Omega$ of $\Omega$ is divided in two nonempty
parts: $\partial\Omega = \partial\Omega_N \cup \partial\Omega_D$. Due to
the thermal equilibrium assumption it follows $np = n_i^2$, and the assumption of vanishing space charge
density gives $n-p-C = 0$, for $x \in \partial\Omega_D$. On the Dirichlet part of the boundary this implies
the following type of boundary conditions:

\begin{eqnarray}
V & \hskip-2.6cm
    = \ V_D(x) \ := \ \ U(x) + V_{\rm bi}(x) & \rm on \ \partial\Omega_D
      \label{eq:dd-dbc1} \\[1ex]
n & \hskip-0.4cm
    = \ n_D(x) \ := \ 
      \T\frac{1}{2} \left( C(x) + \sqrt{C(x)^2 + 4 n_i^2} \right) &
      \rm on \ \partial\Omega_D \label{eq:dd-dbc2} \\
p & = \ p_D(x) \ := \ 
      \ \T\frac{1}{2} \left(-C(x) + \sqrt{C(x)^2 + 4 n_i^2} \right) &
      \rm on \ \partial\Omega_D, \label{eq:dd-dbc3}
\end{eqnarray}
where $U(x)$ is the applied potential, (differences in $U(x)$ between
different segments of $\partial\Omega_D$ correspond to the applied bias
between these two contacts), $V_{\rm bi}(x) := U_T \, \ln
\big( \frac{n_D(x)}{n_i} \big)$ and $U_T$ is the thermal voltage.

Since the Neumann part of the boundary $\partial\Omega_N = \partial\Omega 
- \partial\Omega_D$ models insulating or artificial surfaces, a zero
current flow and a zero electric field in the normal direction are
prescribed. Thus, the following homogeneous boundary conditions are
supplied (in terms of $J_n$ and $J_p$):
\begin{eqnarray}
 \nabla V \cdot \nu & = & 0 \ \ \rm on \ \partial\Omega_N
 \label{eq:dd-nbc1} \\
 J_n \cdot \nu      & = & 0 \ \ \rm on \ \partial\Omega_N
 \label{eq:dd-nbc2} \\
 J_p \cdot \nu      & = & 0 \ \ \rm on \ \partial\Omega_N \, .
 \label{eq:dd-nbc3}
\end{eqnarray}

\begin{figure}[t]
\centerline{ \epsfxsize10cm \epsfbox{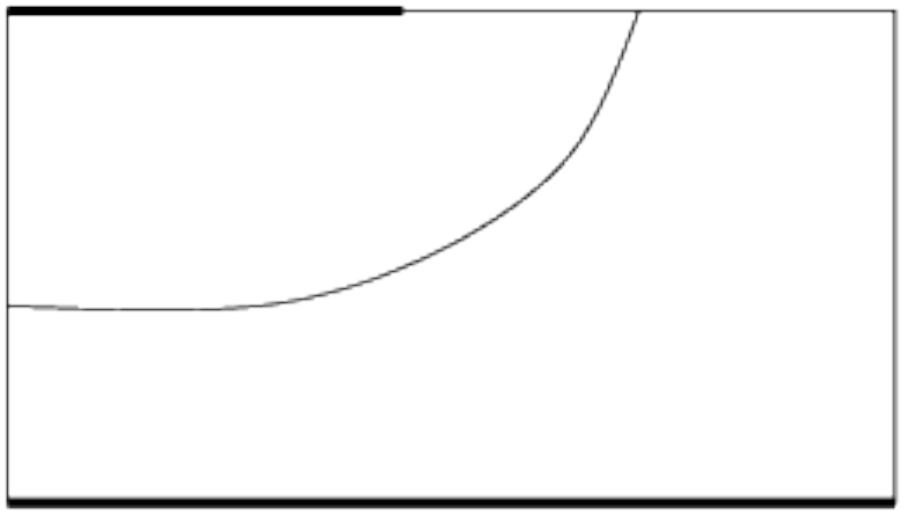} }
\caption{P-N diode. Example of P-N junction.} \label{fig:diode}
\vskip-5cm \unitlength1cm
\centerline{
\begin{picture}(13,5)
% \put(0,0){\dashbox{0.1}(13,7){}}
\put(0.7,4){$\partial\Omega_N$}
\put(11.5,4){$\partial\Omega_N$}
\put(9,6.1){$\partial\Omega_N$}
\put(3,6){$\Gamma_1 \subset \partial\Omega_D$}
\put(6,1.3){$\partial\Omega_D$}
\put(3.5,4.5){\bf P-region}
\put(7.5,3.0){\bf N-region}
\end{picture} }
\end{figure}
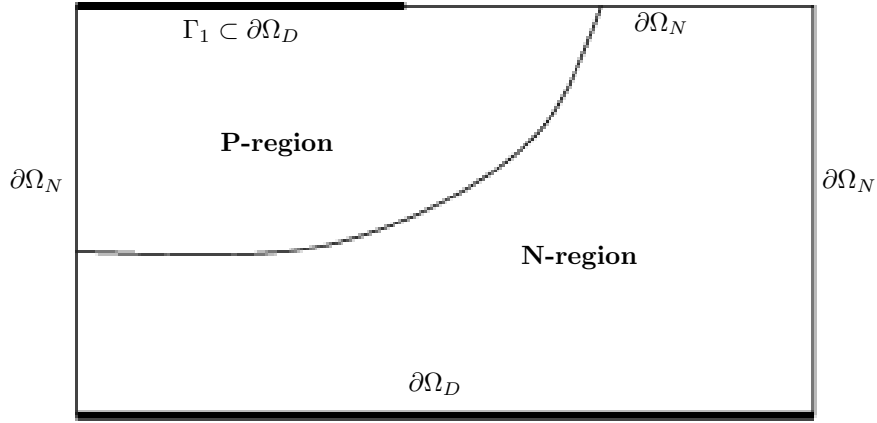

Next we briefly address the modeling of the doping profile. The function
$C(x)$ models a preconcentration of ions in the crystal, so $C(x) = C_{+}(x)
- C_{-}(x)$ holds, where $C_{+}$ and $C_{-}$ are concentrations of negative
and positive ions respectively. In those subregions of $\Omega$ in which the
preconcentration of negative ions predominate (P-regions), we have
$C(x) < 0$. Analogously, we define the N-regions, where $C(x) > 0$ holds.
The boundaries between the P- and N-regions (where $C$ change sign) are
called P-N junctions. An example of a device with a very simple P-N
junction is shown in Figure~\ref{fig:diode}, where a two-dimensional
P-N diode is represented.

Now we introduce an important special change of variables, in order to rewrite
system (\ref{eq:stat-dd-1})--(\ref{eq:stat-dd-3}) as well as boundary
conditions (\ref{eq:dd-dbc1})--(\ref{eq:dd-dbc3}) and
(\ref{eq:dd-nbc1})--(\ref{eq:dd-nbc3}) in a more convenient way. This
variable transformation is (partially) motivated by the Einstein relations
$$ D_n \ = \ U_T \, \mu_n, \ \ \ \ \ D_p \ = \ U_T \, \mu_p, $$
which are a standard assumption about the mobilities and diffusion
coefficients. The so called {\em Slotboom variables} $(u,v)$ are
defined by the relations
\begin{equation} \label{eq:slotboom}
 n(x) \ = \ n_i\, \exp\left(\frac{ V(x)}{U_T}\right)\, u(x),\ \ \
 p(x) \ = \ n_i\, \exp\left(\frac{-V(x)}{U_T}\right)\, v(x) \, .
\end{equation}
For convenience, we rescale the potential and the mobilities:
$$ V(x) \ \leftarrow \ V(x) / U_T,\ \ \
\tilde\mu_n \ := \ q U_T \mu_n,\ \ \ \tilde\mu_p \ := \ q U_T \mu_p \, . $$
Note that the current relations now read:
$$  J_n \ = \   \tilde\mu_n n_i \, e^{ V} \nabla u,\ \ \ \ \ 
    J_p \ = \ - \tilde\mu_p n_i \, e^{-V} \nabla v \, . $$
Therefore, we can write the stationary drift diffusion equation in the
form

\begin{eqnarray}
\lambda^2 \, \Delta V & \hskip-0.3cm \label{eq:dd-sys1}
  = \ \delta^2 \big(e^V u - e^{-V} v \big) - C(x), & {\rm in}\ \Omega \\[1ex]
{\rm div}\, J_n & \hskip-0.4cm \label{eq:dd-sys2}
  = \ \delta^4 \, Q(V,u,v,x) \, (u v - 1), & {\rm in}\ \Omega \\[1ex]
{\rm div}\, J_p & \label{eq:dd-sys3}
  = \ - \delta^4 \, Q(V,u,v,x) \, (u v - 1), & {\rm in}\ \Omega \\[1ex]
V & \hskip-1.9cm = \ V_D \ = \ U + V_{\rm bi} , & \rm on \ \partial\Omega_D
\label{eq:dd-sys4} \\
u & \hskip-2.55cm = \ u_D \ = \ e^{-U} ,         & \rm on \ \partial\Omega_D
\label{eq:dd-sys5} \\
v & \hskip-2.85cm = \ v_D \ = \ e^{U} ,          & \rm on \ \partial\Omega_D
\label{eq:dd-sys6} \\
\nabla V \cdot \nu & \hskip-4.5cm = \ 0 & \rm on \ \partial\Omega_N
\label{eq:dd-sys7} \\
J_n \cdot \nu      & \hskip-4.5cm = \ 0 & \rm on \ \partial\Omega_N
\label{eq:dd-sys8} \\
J_p \cdot \nu      & \hskip-4.5cm = \ 0 & \rm on \ \partial\Omega_N \, ,
                     \label{eq:dd-sys9}
\end{eqnarray}
where $\lambda^2 := \epsilon/(q U_T)$, $\delta^2 := n_i$ and the function
$Q$ is defined implicitly by the relation $Q(V,u,v,x) = {\cal R}(n,p,x)$.
Notice the applied potential has also to be rescaled: $U(x) \leftarrow U(x) / U_T$.

\subsection{Some existence and uniqueness results for stationary drift
diffusion equation} \label{ssec:ex-uniq}

In this subsection we discuss the solution theory for the system of
drift diffusion equations (\ref{eq:dd-sys1})--(\ref{eq:dd-sys9}). First
we present an existence result, which can be found in the monography
\cite[Theorem 3.3.16]{MRS}.

\begin{propo} \label{prop:MRS-3316}
Let $\kappa > 1$ be a constant satisfying
$$ \kappa^{-1} \ \le \ u_D(x)\, ,\ v_D(x) \ \le \ \kappa, \
   \forall x \, \in \, \partial\Omega_D \, , $$
and let $-\infty < C_m \leq C_M < + \infty$. Then for any \ $C \in
\{ L^\infty(\Omega)\, ; \ C_m \le C(x) \le C_M,\ x \in \Omega \}$
\ the boundary value problem (\ref{eq:dd-sys1})--(\ref{eq:dd-sys9})
admits a weak solution $(V,u,v) \in ( H^1(\Omega) \cap L^\infty(\Omega) )^3$
satisfying
$$ \kappa^{-1} \ \le \ u(x)\, ,\ v(x) \ \le \ \kappa, \
   \forall x \, \in \, \Omega \, , $$
furthermore
\begin{eqnarray*}
V(x) & \ge & \min \left( \inf_{\partial\Omega_D} V_D \, ,\ U_T
     \ln \left[ \frac{1}{2\kappa n_i} (C_m + (C_m^2 + 4 n_i^2)^{1/2} ) \right]
     \right) ,\ {\rm in}\ \Omega \\
V(x) & \le & \max \left( \sup_{\partial\Omega_D} V_D \, ,\ U_T \ln
     \left[ \frac{\kappa}{2 n_i} (C_M + (C_M^2 + 4 n_i^2)^{1/2} ) \right]
     \right) ,\ {\rm in}\ \Omega \, .
\end{eqnarray*}
\end{propo}
{\it Sketch of the proof:} \\
Solving the Poisson equation and the continuity equations (three elliptic
mixed boundary value problems), one at a time, it is possible to define
an operator in an appropriate $L^2$--space, whose fixed point is a weak
solution of (\ref{eq:dd-sys1})--(\ref{eq:dd-sys9}). The existence of a
fixed point is established by the Schauder Fixed point operator, once one
proves that the fixed point operator is completely continuous, which is
accomplished by the use of standard elliptic theory. \hfill $\Box$
\bigskip

As far as uniqueness of solutions of system 
(\ref{eq:dd-sys1})--(\ref{eq:dd-sys9}) is concerned, a corresponding result
can be obtained if the applied voltage is small (in the norm of 
$L^\infty(\partial\Omega_D) \cap H^{3/2}(\partial\Omega_D)$). The
following uniqueness result corresponds to \cite[Theorem 2.4]{BEMP}.

\begin{propo} \label{prop:BEMP-24}
Let the voltage $U$ be such that $\|U\|_{L^\infty(\partial\Omega_D)}
+ \|U\|_{H^{3/2}(\partial\Omega_D)}$ is sufficiently small. Then system
(\ref{eq:dd-sys1})--(\ref{eq:dd-sys9}) has a unique solution $(V,u,v) \in
( H^1(\Omega) \cap L^\infty(\Omega) )^3$.
\end{propo}

Since existence and uniqueness of solutions for system
(\ref{eq:dd-sys1})--(\ref{eq:dd-sys9}) can be guaranteed for
small applied voltages only, it is reasonable to consider instead of
this system its linearized version around the equilibrium point
$U \equiv 0$ instead. We shall follow this approach through a large part of this paper.

Under stronger assumptions on the boundary parts $\partial\Omega_D$, 
$\partial\Omega_N$ as well as on the boundary conditions $V_D$, $u_D$,
$v_D$, it is even possible to show $H^2$-regularity for a solution $(V,u,v)$
of (\ref{eq:dd-sys1})--(\ref{eq:dd-sys9}). Next we shall discuss a
corresponding result; for the proof details we refer to the monography
by Markowich \cite{Ma}. First we have to consider the following assumptions:

\begin{enumerate}
\item[A1)] $\Omega$ is a bounded domain of class $C^{0,1}$ in $\mathcal R^d$
and the $(d-1)$-dimensional Lebesgue measure of $\partial\Omega_D$ is
positive;
\item[A2)] The Dirichlet boundary data $(V_D,u_D,v_D)$ in
(\ref{eq:dd-sys4})--(\ref{eq:dd-sys6}) satisfy
$$ (V_D,u_D,v_D) \in (H^2(\Omega))^3, \  (V_D,u_D,v_D)|_{\partial\Omega_D}
\in (L^\infty(\partial\Omega_D))^3\, . $$
Furthermore, $(V_D,u_D,v_D)|_{\partial\Omega_N} = (0,0,0)$ and there is
$U_+ \ge 0$ such that
$$ e^{-U_+} \le \D\inf_{\partial\Omega_D} u_D, \ \inf_{\partial\Omega_D} v_D;
\ \ \D\sup_{\partial\Omega_D} u_D, \ \sup_{\partial\Omega_D} v_D \le
e^{U_+}\, ; $$
\item[A3)] The doping profile satisfies  $C \in L^\infty(\Omega)$;
\item[A4)] The function $Q$ in (\ref{eq:dd-sys2}), (\ref{eq:dd-sys3}) is such
that $Q(\cdot,\cdot,\cdot,x) \in C^1(\mathcal R \times (0,\infty)^2)$ for all
$x \in \Omega$; $Q(V,u,v,\cdot)$, $\nabla_{(V,u,v)}Q(V,u,v,\cdot) \in L^\infty
(\Omega)$ uniformly for $(V,u,v)$ in bounded sets of $\mathcal R \times
(0,\infty)^2$; $Q(V,u,v,x) \ge 0$ in $\mathcal R \times (0,\infty)^2 \times
\Omega$;
\item[A5)] The mobilities $\mu_n$, $\mu_p$ satisfy: 
$\mu_n = \mu_n(x)$, $\mu_p = \mu_p(x)$, $\mu_n, \mu_p \in W^{1,\infty}
(\Omega)$; functions $\mu_n$, $\mu_p$ are both positive and uniformly
bounded (away from zero) in $\Omega$;
\item[A6)] The solution of
$$ \Delta w = f \ {\rm in} \ \Omega\, , \ \ \ \ w|_{\partial\Omega_D} =
w_\nu|_{\partial\Omega_N} = 0\, , $$
satisfies \ $\| w \|_{2,q,\Omega} \le K_1 \| f \|_{q,\Omega}$ for every
$f \in L^q(\Omega)$ with $q=2$ and $q=3/2$.
\end{enumerate}

\begin{propo} \mbox{\bf \cite[Theorem 3.3.1]{Ma}}
\label{prop:Ma-331}
Let's assume that assumptions A1) -- A6) hold. Then every weak solution
$(V,u,v) \in (H^1(\Omega) \cap L^\infty(\Omega))^3$ satisfies
$$ (V,u,v) \in (H^2(\Omega))^3 \, . $$
\end{propo}

Notice that, in the 2-dimensional case, assumption $A6)$ can only be
satisfied (for $q=2$) if the angle between the Neumann and Dirichlet
parts of $\partial\Omega$ is smaller than $\pi/2$. Otherwise, the
solutions of the elliptic mixed boundary value problems in $A6)$
will not belong to any space $H^{1+\epsilon}(\Omega)$, $\epsilon>0$
(see \cite{Gr} for details).

If assumption $A6)$ holds only for $q=2$, it is still possible to prove
$H^2$-regularity for the solution $V$ of the Poisson equation. However,
it is not possible to prove square-integrability of all second derivatives
of $u$ and $v$, but this is usually not needed for the formulation of the
inverse problem.

\subsection{The equilibrium case}

In this subsection we analyze the equilibrium case for the stationary
drift diffusion equations, which corresponds to the assumption
$U(x) \equiv 0$. In this particular case several simplifications are
possible. It is immediate to see that the solution of
(\ref{eq:dd-sys1})--(\ref{eq:dd-sys9}) is of the form
$(V=V^0, u \equiv 1, v \equiv 1)$, where
\begin{equation} \label{eq:poiss-equil}
\left\{ \begin{array}{rcll}
   \lambda^2 \, \Delta V^0 & = & e^{V^0} - e^{-V^0} - C(x)
   & {\rm in}\ \Omega \\
   V^0 & = & V_{\rm bi}(x) & {\rm on}\ \partial\Omega_D \\
   \nabla V^0 \cdot \nu & = & 0 & {\rm on}\ \partial\Omega_N \, .
\end{array} \right.
\end{equation}

As already mentioned in subsection (\ref{ssec:ex-uniq}), we shall be
interested in the linearized drift diffusion system at the equilibrium.
Keeping this in mind, we compute the derivative of the solution of
(\ref{eq:dd-sys1})--(\ref{eq:dd-sys9}) with respect to the voltage $U$ at $U \equiv 0$ in the direction
$h$. This directional derivative is given by the solution $(\hat V,
\hat u, \hat v)$ of
\begin{equation} \label{eq:poiss-equil-lin}
\left\{ \begin{array}{ll}
   \lambda^2 \, \Delta \hat V \ = \ e^{V^0} \hat u + e^{-V^0} \hat v
   + \big(e^{V^0} + e^{-V^0}\big) \hat V & {\rm in}\ \Omega \\[1ex]
   {\rm div}\, (\mu_n e^{V^0} \nabla \hat u) \ = \ Q_0(V^0,x)
   (\hat u + \hat v) & {\rm in}\ \Omega \\[1ex]
   {\rm div}\, (\mu_p e^{-V^0} \nabla \hat v) \ = \ Q_0(V^0,x)
   (\hat u + \hat v) & {\rm in}\ \Omega \\[1ex]
   \hat V \ = \  h  & {\rm on}\ \partial\Omega_D \\[1ex]
   \hat u \ = \ -h & {\rm on}\ \partial\Omega_D \\[1ex]
   \hat v \ = \  h  & {\rm on}\ \partial\Omega_D \\
   \D\frac{\partial V^0}{\partial\nu}   \ = \
   \D\frac{\partial \hat u}{\partial\nu}\ = \
   \D\frac{\partial \hat v}{\partial\nu}\ = \ 0 & {\rm on}\
   \partial\Omega_N \, ,
\end{array} \right.
\end{equation}
where $Q_0(V^0,x) \ = \ Q(V^0,1,1,x)$. Notice that, in the linearized
case close to equilibrium, the solutions $(\hat u, \hat v)$ of the
continuity equations do not depend on the electrostatic potential
$\hat V$.

\subsection{Unipolar and bipolar cases} \label{ssec:ubc}

In this subsection we introduce two special cases, which are going to
play a key rule in the modeling of some of the inverse problems analyzed in
this paper. We start by introducing the operator called {\em voltage-current}
(V--C) map:
$$ \begin{array}{rcl}
   \Sigma_C: H^{3/2}(\partial\Omega_D) & \to & H^{1/2}(\Gamma_1) \\
   U & \mapsto & J\cdot\nu|_{\Gamma_1}\ =\ (J_n+J_p)\cdot\nu |_{\Gamma_1}\, ,
   \end{array} $$
where $\Gamma_1 \subset \partial\Omega_D$ (see Figure~\ref{fig:diode}
for an example). The map $\Sigma_C$ takes the applied voltage $U$ into
the outflow current density on $\Gamma_1$. The linearized unipolar case (close
to equilibrium) corresponds to the model obtained from the unipolar drift diffusion
equations by linearizing the V--C map at $U \equiv 0$. This simplification
is motivated by the fact that the V--C map can only be defined as a single-valued function in a
neighborhood of $U=0$. Furthermore, the following assumptions are also taken
into account
\begin{itemize}
\item[{\it i)}] The concentration of holes satisfy $p = 0$ (or, equivalently,
$v = 0$ in $\Omega$);
\item[{\it ii)}] No recombination-generation rate is present, i.e.
${\cal R} = 0$ (or $Q = 0$).
\end{itemize}
Under this assumptions, system (\ref{eq:dd-sys1})--(\ref{eq:dd-sys9})
reduces to the decoupled system:
\begin{equation} \label{eq:unipolar}
\left\{ \begin{array}{ll}
   \lambda^2 \, \Delta V^0 \ = \ e^{V^0} - C(x) & {\rm in}\ \Omega \\
   {\rm div}\, (e^{V^0} \nabla \hat{u}) \ = \ 0 &  {\rm in}\ \Omega \\[1ex]
   V^0 \ = \ V_{\rm bi}(x) & {\rm on}\ \Omega_D \\
   \hat{u} \ = \ U(x) & {\rm on}\ \Omega_D \\[1ex]
   \nabla V^0 \cdot \nu \ = \ 0 & {\rm on}\ \Omega_N \\
   \hat{J}_n \cdot \nu \ = \ 0 & {\rm on}\ \Omega_N
\end{array} \right.
\end{equation}

The inverse problem of identifying the doping profile in the linearized
unipolar model (\ref{eq:unipolar}) corresponds to identification of $C(x)$
from the map
$$ \Sigma_C'(0): U \mapsto (\hat{J}_n \cdot \nu) |_{\Gamma_1} \, , \qquad \hat{J}_n := \mu_n e^{V_0} \nabla \hat{u}. $$
Notice that, since $V = V_{\rm bi}$ is known at $\partial\Omega_D$,
the current data (output) $J_n \cdot \nu = \mu_n e^{V^0} u_\nu$ can be 
directly replaced by the Neumann data $u_\nu$.

We shall return to this identification problem in Section~\ref{sec:idp},
where the inverse problem described above is considered as a generalization
of the well known {\em electrical impedance tomography} or {\em inverse
conductivity problem} (see, e.g., \cite{Bo,Is} for a survey on these inverse
problems).

Next we concentrate on deriving the so called {\em bipolar case}. As in
the unipolar case, we will be interested in reconstructing the doping
profile $C$ in (\ref{eq:dd-sys1})--(\ref{eq:dd-sys9}) from the linearized 
V--C map at $U \equiv 0$. This is an interesting case, due to the fact that
the Poisson Equation and the continuity equations decouple.

From (\ref{eq:poiss-equil-lin}) we see that the Gateaux derivative of
the V--C map $\Sigma_C$ at the point $U=0$ in the direction $\Phi$ is
given by the expression
$$ \Sigma'_C(0) \Phi \ := \ \left( \mu_n \, e^{V_{\rm bi}} \hat{u}_\nu -
\mu_p \, e^{-V_{\rm bi}} \hat{v}_\nu \right) |_{\Gamma_1} \, , $$
where $(u,v)$ solve
\begin{equation} \label{eq:bip-lin}
   \left\{ \begin{array}{ll}
     {\rm div}\, (\mu_n e^{V^0} \nabla \hat{u}) \ = \ Q_0(V^0,x) (\hat{u} + \hat{v})
     & {\rm in}\ \Omega \\
     {\rm div}\, (\mu_p e^{-V^0} \nabla \hat{v}) \ = \ Q_0(V^0,x) (\hat{u} + \hat{v})
     & {\rm in}\ \Omega \\[1ex]
     \hat{u} \ = \ -\Phi              & {\rm on}\ \partial\Omega_D \\
     \hat{v} \ = \ \Phi               & {\rm on}\ \partial\Omega_D \\
     \D\frac{\partial \hat{u}}{\partial \nu} \ = \
     \frac{\partial \hat{v}}{\partial \nu} \ = \ 0  & {\rm on}\ \partial\Omega_N
   \end{array} \right.
\end{equation}
and $V^0$ is the solution of the equilibrium problem (\ref{eq:poiss-equil}).

Notice that the solution of the Poisson equation can be computed a priori,
since it does not depend on $\Phi$. The linear operator $\Sigma'_C(0)$
is continuous. Actually we can prove more: since $(u,v)$ depend
continuously (in $H^2(\Omega)^2$) on the boundary data $\Phi$ (in
$H^{3/2} (\partial\Omega_D)$), it follows from the boundedness and
compactness of the trace operator $\gamma: H^2(\Omega) \to
H^{1/2}(\Gamma_1)$ that $\Sigma'_C(0)$ is a bounded and compact
operator. The application $\Sigma_C'(0)$ maps the Dirichlet data for
$(u,v)$ to a weighted sum of their Neumann data and can be compared with the
identification problem in the electrical impedance tomography.

%---------------------------------------------------------------------------
\subsection{Flipped bipolar case} \label{ssec:fbc}

In this subsection we introduce another special case, which will be
relevant for the formulation of the inverse problem related to the {\em
laser-beam-induced current} (LBIC) measurements. We start by introducing
the LBIC functional defined by the boundary integral
$$ \begin{array}{rcl}
   \mathcal{I}: L^2(\Omega) & \to & \mathcal{R} \\
   g & \mapsto & \D\int_{\Gamma_1}
   \left\{ \mu_n e^{V^0} \hat{u}_\nu - \mu_p e^{-V^0} \hat{v}_\nu \right\} ds
   \end{array} $$
where $\Gamma_1 \subset \partial\Omega_D$ is defined as in Subsection%
~\ref{ssec:ubc}, $V^0$ is the solution of the equilibrium problem
(\ref{eq:poiss-equil}) and $(u,v)$ solve
\begin{equation} \label{eq:fbip-lin1}
   \left\{ \begin{array}{ll}
     {\rm div}\, (\mu_n e^{V^0} \nabla \hat{u}) \ = \ Q_0(V^0,x) (\hat{u} + \hat{v}) + g
     & {\rm in}\ \Omega \\
     {\rm div}\, (\mu_p e^{-V^0} \nabla \hat{v}) \ = \ Q_0(V^0,x) (\hat{u} + \hat{v}) + g
     & {\rm in}\ \Omega \\[1ex]
     \hat{u} \ = \ \hat{v} \ = \ 0          & {\rm on}\ \partial\Omega_D \\
     \D\frac{\partial \hat{u}}{\partial \nu} \ = \
     \frac{\partial \hat{v}}{\partial \nu} \ = \ 0  & {\rm on}\ \partial\Omega_N
   \end{array} \right.
\end{equation}

Notice that the only differences between systems (\ref{eq:fbip-lin1}) and 
(\ref{eq:bip-lin}) (from the bipolar case) are: 1) the $L^2(\Omega)$ function
$g$, appearing on the right hand side of the linearized continuity
equations and representing the applied laser beam (see Section~\ref{sec:idp}
for details on the problem formulation); 2) the Dirichlet boundary condition
at $\partial\Omega_D$.

The inverse problem of reconstructing the doping profile from measurements
of the LBIC functional was considered in \cite{FI1,FI2}. An alternative
representation for the functional $\mathcal I$ was derived in \cite{FI1}.
Analyzing the variational formulation of the system constituted by
(\ref{eq:poiss-equil}) and (\ref{eq:fbip-lin1}) and using standard
functional analytical arguments as well as basic elliptic theory (see
\cite{GT}), the authors proved that the LBIC functional can be written as
\begin{equation} \label{eq:lbic-repr}
\mathcal I(g) \ = \ \ipl \hat \tilde{v} - \tilde{u}, \, g \ipr_{L^2(\Omega)} \, ,
\end{equation}
where $V^0$ is defined as before and $(\hat u, \hat v)$ solve the system
\begin{equation} \label{eq:fbip-lin2}
   \left\{ \begin{array}{ll}
     {\rm div}\, (\mu_n e^{V^0} \nabla \tilde{u}) \ = \ Q_0(V^0,x)
     (\tilde{u} - \tilde{v}) & {\rm in}\ \Omega \\
     {\rm div}\, (\mu_p e^{-V^0} \nabla \tilde{v}) \ = \ Q_0(V^0,x)
     (\tilde{v} - \hat \tilde{u}) & {\rm in}\ \Omega \\[1ex]
     \tilde{u} \ = \ \tilde{v} \ = \ 1  & {\rm on}\ \Gamma_1 \\
     \tilde{u} \ = \ \tilde{v} \ = \ 0  & {\rm on}\ \partial\Omega_D / \Gamma_1 \\
     \D\frac{\partial \tilde{u}}{\partial \nu} \ = \
     \frac{\partial \tilde{u}}{\partial \nu} \ = \ 0  & {\rm on}\ \partial\Omega_N
   \end{array} \right.
\end{equation}

We shall refer to system (\ref{eq:poiss-equil}), (\ref{eq:fbip-lin2})
as {\em flipped bipolar case}. As in the bipolar case, the solution of the
Poisson equation can be computed a priori, since $V^0$ does not depend
on $g$. Therefore, to evaluate $\mathcal I$ in (\ref{eq:lbic-repr}) one
needs only to solve the coupled system (\ref{eq:fbip-lin2}). Moreover,
from the representation formula (\ref{eq:lbic-repr}), it follows that
$\mathcal I$ is a linear continuous functional on $L^2(\Omega)$.

%---------------------------------------------------------------------------
\section{Inverse doping problem} \label{sec:idp}

The so called {\em inverse doping profile} corresponds to the problem
of identifying a doping profile $C(x)$ in system
(\ref{eq:dd-sys1})--(\ref{eq:dd-sys9}) from indirect measurements. In
practical applications, the following types of measurements are available (cf. \cite{Khetal95}):

\begin{enumerate}
\item  Current flow through a contact $\Gamma_1 \subset \partial\Omega_D$:
$$ I(U) \ = \ \int_{\Gamma_1} (J_n + J_p) .\nu \, ds \, , $$
where $U \in H^{3/2}(\partial\Omega_D)$ \ with \ $\| U \|$ \ small. \\
Under the (idealized, but technologically realizable) assumption that we not only know the averaged
flow through $\Gamma_1$, but the actual flow $J \cdot \nu$ on $\Gamma_1$,
this type of measurement corresponds to the voltage-current map introduced
in Subsection~\ref{ssec:ubc}:
$$  \Sigma_C(U) \ := \ (J_n + J_p) \cdot \nu \big|_{\Gamma_1} \ \in \
    H^{1/2}(\Gamma_1)\, . $$
\item  Mean capacitance of a contact $\Gamma_1 \subset \partial\Omega_D$:
$$ Cap(U) \ = \ \frac{\partial}{\partial U}
   \left( \int_{\Gamma_1} \nabla V . \nu \, ds \right) \, . $$
We shall consider the idealized (but again technologically realizable) data corresponding to measurements of the
variation of the electric flux (in the normal outward direction) with respect
to an applied voltage $U$ at $\partial\Omega_D$. This data corresponds to
the so called {\em capacitance measurements}
$$ {\cal T}_C(U) \ := \ \frac{\partial }{\partial U} \ \frac{\partial V}{\partial\nu} \Big|_{\Gamma_1}
   = \frac{\partial \hat{V}}{\partial\nu} \Big|_{\Gamma_1}
   \ \in \ H^{1/2}(\Gamma_1) \, , $$
here $V$ is the solution of the Poisson equation for an applied voltage
$U \in H^{3/2}(\partial\Omega_D)$.
\item Measurements of the total current $i(x)$ flowing out through one
contact induced by a laser beam applied at different locations $x \in \Omega$:
$$ i(x) \ := \ \, \mathcal I( \delta(\cdot - x)) \ = \
   \hat v(x) - \hat u(x)\, , $$
where $\mathcal I$ is the LBIC functional defined in Subsection~%
\ref{ssec:fbc}, and $(\hat u, \hat v)$ is the solution of system
(\ref{eq:fbip-lin2}). These data correspond to the so called {\em
laser-beam-inducted current measurements}.
\end{enumerate}
In all cases we assume that $\Gamma_1 \subset \partial\Omega_D$ is
sufficiently regular with non zero measure. The first step in the
investigation of the inverse problems modeled by operators $\Sigma_C$
and ${\cal T}_C$ consists in analyzing whether these operators are
well defined in appropriate spaces. The next three subsections are devoted
to the analysis of each of these operators. In the last subsection we
discuss in details the inverse problem related to the V--C map for the
linearized unipolar case close to equilibrium and its relation with the
electrical impedance tomography.

\subsection{The voltage-current map} \label{ssec:vcmap}

in this subsection we analyze the V--C map introduced above.
The map $\Sigma_C$ takes (for a fixed doping profile $C$) the
applied voltage  $U$ into the corresponding current density.
The non-linear operator $\Sigma_C$ is well-defined, when
considered as a map between suitable Sobolev spaces. This
assertion is a consequence of the following result:

\begin{propo} \mbox{\bf \cite[Proposition 3.1]{BEMP}}
\label{prop:bemp31}
For each applied voltage $U \in B_r(0) \subset H^{3/2}(\partial\Omega_D)$
with $r>0$ sufficiently small, the current $J \cdot \nu \in H^{1/2}
(\Gamma_1)$ is uniquely defined. Furthermore, $\Sigma_C:
H^{3/2}(\partial\Omega_D) \to H^{1/2}(\Gamma_1)$ is continuous
and is continuously differentiable in $B_r(0)$.
\end{propo}
{\it Sketch of the proof:} \\
The first part of the proof follows basically from the uniqueness of
solutions for system (\ref{eq:dd-sys1})--(\ref{eq:dd-sys9}) in
$H^2(\Omega)^3$ together with regularity properties of the {\em
Neumann trace operator} $\gamma: H^2(\Omega) \to H^{1/2}(\Gamma_1)$.
The Fr\'echet-differentiability follows from standard estimates of the
residual in the Taylor expansion of the operator $\Sigma_C$. \hfill $\Box$

\bigskip

By iterating the argumentation in Proposition~\ref{prop:bemp31}, one 
can even prove that $\Sigma_C$ is of class $C^\infty$ in $B_r(0)
\subset H^{3/2}(\partial\Omega_D)$ for $r$ sufficiently small.

Proposition~\ref{prop:bemp31} establishes a basic property to consider
the inverse problem of reconstructing the doping profile $C$ from the
V--C map. In the sequel we shall consider two possible inverse problems
for the V--C map.

In the first inverse problem we assume that, for each $C$, the output
corresponds to the map $\Sigma_C$. A realistic experiment corresponds
to measure, for given $\{ U_j \}_{j=1}^N$, with $\|U_j\|=1$, the outputs
$$ \big\{ \Sigma_C(t U_j)\ |\ \ j=1,\cdots,N; \ \ t \in [0,r] \big\} \, . $$
This data corresponds to the assumption that the amplitude of an applied
voltage $U_j$ can be varied continuously in a practical experiment. In 
practice, the functions $U_j$ are chosen to be piecewise constant.

Notice that, for fixed $U_j$, the continuity of $\Sigma_C$ implies the
continuity of the function $t \mapsto \Sigma_C (t U_j)$. Therefore,
we have
$$ \Sigma_C(\cdot U_j) \ \in \ C(0,r; H^{1/2}(\Gamma_1)) \ \subset \
   L^2(0,r; L^2(\Gamma_1)) \, . $$
This allow the following abstract formulation of the inverse problem for
the V--C map:
\begin{equation} \label{eq:ip-abstract}
 F(C) \ = \ Y \, ,
\end{equation}
where
\begin{enumerate}
\item[1)] Parameter: \ $C = C(x) \ \in \ L^2(\Omega) =: \mathcal X$;
\vskip-5ex
\item[2)] Output: \ $Y = \big\{ \Sigma_C(\cdot U_j) \big\}_{j=1}^N \in
[L^2((0,r) \times \Gamma_1)]^N =: \mathcal Y$;
\vskip-5ex
\item[3)] Parameter-to-output map: \ $F: \mathcal X \to \mathcal Y$.
\end{enumerate}
The domain of definition of the operator $F$ is
$$ D(F) := \{ C \in L^2(\Omega) ; \, C_m \le C(x) \le C_M,
   \mbox{ a.e. in } \Omega \} \, , $$
where $\underline{C}$ and $\overline{C}$ are appropriate positive constants.

This choice of spaces is motivated by Propositions~\ref{prop:BEMP-24}
and~\ref{prop:Ma-331}, which guarantee, for each $t \|U_j\| < r$ ($r$
small) and \ $C \in D(F)$, the existence and uniqueness of a $H^2$-%
solution $(V,n,p)$ for system (\ref{eq:dd-sys1})--(\ref{eq:dd-sys9}).
Therefore, the map
$$ \begin{array}[t]{rcl}
     F: D(F) \subset {\mathcal X} & \to & {\mathcal Y} \\
     C & \mapsto & \big\{ \Sigma_C(\cdot U_j) \big\}_{j=1}^N \end{array} $$
is well defined. Furthermore, $F$ is also Fr\'echet-differentiable in
$D(F)$. Indeed, we already know that the map $(V,u,v) \mapsto J \cdot \nu
|_{\Gamma_1}$ is continuously differentiable (this is included in the
proof of Proposition~\ref{prop:bemp31}). Thus, it is enough to verify
the differentiability of the map $\mathcal U_j: D(F) \ni C \mapsto (V,u,v)
\in H^2(\Omega)^3$, for fixed $U_j$. The variation of the solution $(V,u,v)$
of system (\ref{eq:dd-sys1})--(\ref{eq:dd-sys9}) with respect to a
variation of the doping profile $C$ can be deduced similarly as in
(\ref{eq:poiss-equil-lin}). To prove Fr\'echet-differentiability of
$\mathcal U_J$, we only have to estimate the residual in the Taylor
expansion of this map, as in the proof of Proposition~\ref{prop:bemp31}.

This inverse problem is addressed in the literature as identification
of doping profiles from {\em full voltage-current data}. Next we shall
investigate a different formulation of the same inverse problem related
to the V--C map.

In practical applications, the V--C map can only be defined in a
neighborhood of $U=0$ (due to hysteresis defects for large  applied
voltages). This motivates the analysis of the problem of identifying
the doping profile $C$ from the linearized V--C map at $U=0$. (see
unipolar and bipolar cases in Subsection~\ref{ssec:ubc}).

As described in Subsection~\ref{ssec:ubc}, the Gateaux derivative of
the V--C map at $U=0$ in direction $\Phi$ is given by
$$ \Sigma'_C(0) \Phi \ = \ \left( \mu_n \, e^{V_{\rm bi}} u_\nu -
   \mu_p \, e^{-V_{\rm bi}} v_\nu \right) |_{\Gamma_1} \, , $$
where $(u,v)$ solve the system in (\ref{eq:bip-lin}) and $V^0$ is
the solution of the equilibrium case (\ref{eq:poiss-equil}). In
Subsection~\ref{ssec:ubc} we have already verified the boundedness and
compactness of $\Sigma'_C(0)$. Contrary to the case of full data, the
solution of the Poisson equation can be computed {\em a priori}, since
it is independent of $\Phi$.

The next step to complete the setup of this second inverse problem for
the V--C map is to define the problem data. The data for the problem
can be obtained from the full V--C data:
$$ Y \ := \ \big\{ \Sigma'_C(0) U_j \big\}_{j=1}^N \ \in \
   \big[ L^2(\Gamma_1) \big]^N \, . $$
In the literature, this inverse problem in called identification of
doping profiles from {\em reduced voltage-current data}. Notice that
the functions $U_j$ are defined as before. Therefore, we obtain for
the inverse problem with reduced data the same abstract formulation
as in (\ref{eq:ip-abstract}) with
\begin{enumerate}
\item[1)] Parameter: \ $C = C(x) \ \in \ L^2(\Omega) =: \mathcal X$;
\vskip-5ex
\item[2)] Output: \ $Y = \big\{ \Sigma_C'(\cdot U_j) \big\}_{j=1}^N \in
\big[ L^2(\Gamma_1) \big]^N =: \mathcal Y$;
\vskip-5ex
\item[3)] Parameter-to-output map: \ $F: \mathcal X \to \mathcal Y$.
\end{enumerate}
The domain of definition of the operator $F$ is the same as in the case
of full data. Notice that the parameter-to-output operator for reduced
data is given by:
$$ F : \begin{array}[t]{rcl}
D(F) \subset \mathcal X & \to     & \mathcal Y \\
C                        & \mapsto & \big\{ \Sigma'_C(0) U_j \big\}_{j=1}^N
\end{array} $$
Analogously as in the full V--C data case, one can prove that the
non-linear parameter-to-output operator is well defined and Fr\'echet
differentiable in its domain of definition $D(F)$.

As already observed, the solution of the Poisson equation can be
computed {\em a priori}. The remaining problem (coupled system for
$(u,v)$) is quite similar to the problem of {\em electrical impedance
tomography}. In this inverse problem the aim is to identify the
conductivity $q = q(x)$ in the equation:
$$ -{\rm div}\,(q \nabla u) \ = \ f \ \  {\rm in}\ \Omega \, , $$
from measurements of the {\em Dirichlet-to-Neumann map}, which maps
the applied voltage $u|_{\partial\Omega}$ to the electrical flux
$q u_\nu|_{\partial\Omega}$. The application $\Sigma_C'(0)$ maps the
Dirichlet data for $\hat{u}$ and $\hat{v}$ to the weighted sum of their
Neumann data. It can be seen as the counterpart of electrical impedance
tomography for common conducting materials.

We close this subsection discussing yet another inverse problem for
the capacitance measurements. This problem again concerns the reduced
V--C map and arrises in a limiting case of the drift diffusion equations,
called limit of {\em zero space charge}, which is mathematically represented
by the scaling limit $\lambda \to 0$. In this case the Poisson equation reduces to an algebraic
relation between $V$ and $C$ and existence of solutions of the zero-space-charge problem in $L^\infty(\Omega)$
has been proven in \cite{MRS}.

Notice that, in the limiting case $\lambda = 0$, without further regularity
assumptions on the doping profile $C$ we can only guarantee $H^1$ regularity
for a solution $(u,v)$ of (\ref{eq:bip-lin}). Therefore, $J \cdot \nu \in
H^{-1/2}(\Gamma_1)$ follows.
However, as already observed in \cite{BEMP}, if $\nabla C \in L^p(\Omega)$
for $p$ sufficiently large ($p \ge 6$), one can show that the reduced V--C
map exists and maps continuously to $L^2(\Gamma_1)$.

From the Poisson equation in equilibrium we obtain $\sinh V = 2C$ and the
linearized continuity equations can be written in the form
\begin{equation} \label{eq:zsc-lin}
   \left\{ \begin{array}{ll}
     {\rm div}\, (\mu_n a \nabla \hat u) \ = \ q(a,x) (\hat u + \hat v)
     & {\rm in}\ \Omega \\
     {\rm div}\, (\mu_p a^{-1} \nabla \hat v) \ = \ q(a,x) (\hat u + \hat v)
     & {\rm in}\ \Omega \\[1ex]
     u \ = \ -\Phi              & {\rm on}\ \partial\Omega_D \\
     v \ = \ \Phi               & {\rm on}\ \partial\Omega_D \\
     \D\frac{\partial u}{\partial \nu} \ = \
     \frac{\partial v}{\partial \nu} \ = \ 0  & {\rm on}\ \partial\Omega_N
   \end{array} \right.
\end{equation}
where
$$  a \, = \, a(C) \, = \, e^{{\rm arcsinh }(2C)} \, ,\ \ \ \
    q(a,x) \, = \, Q(\ln(a), 1, 1, x) \, . $$

Thus, in this limiting case, the inverse doping profile problem reduces
to the identification of the conductivity $a$ in the coupled system
(\ref{eq:zsc-lin}) from the reduced V--C map. Once we have reconstructed
the coefficient $a$, the doping profile can be obtained from the relation
$C = \frac{1}{2} \sinh( \ln a)$.

\subsection{Capacitance measurements}

In this subsection we address the inverse problem modeled by the
operator ${\cal T}_C$, introduced at the beginning of Section~\ref{sec:idp}.

The operator ${\cal T}_C$ maps an applied voltage $U$ at $\partial\Omega_D$
to the idealized data corresponding to the Neumann trace of the electric
potential $\hat V$ at $\Gamma_1 \subset \partial\Omega_D$, i.e.
$$ \begin{array}{rcl}
     {\cal T}_C: H^{3/2}(\partial\Omega_D) & \to & H^{1/2}(\Gamma_1) \\
     U &  \mapsto & \D\frac{\partial \hat V}{\partial\nu} \Big|_{\Gamma_1}
   \end{array}  $$
where $\hat V$ solves:
$$ \left\{ \begin{array}{rcll}
   \lambda^2 \, \Delta \hat V & = & \big( e^{V^0} + e^{-V^0} \big)\, \hat V
                 + e^{V^0} \hat u + e^{-V^0} \hat v & \mbox{ in } \Omega \\
   \hat V & = & U & \mbox{ on } \partial\Omega_D \\
   \nabla \hat V \cdot \nu & = & 0 & \mbox{ on } \partial\Omega_N
   \end{array} \right. $$
here $V^0$ is the solution the equilibrium case (\ref{eq:poiss-equil}) and
$(\hat u, \hat v)$ is the solution of the system
$$ \left\{ \begin{array}{ll}
     {\rm div}\, (\mu_n e^{V^0} \nabla \hat u) \ = \ Q_0(V^0,x)
     (\hat u + \hat v)   & \mbox{ in } \Omega \\
     {\rm div}\, (\mu_p e^{-V^0} \nabla \hat v) \ = \ Q_0(V^0,x)
     (\hat u + \hat v)   & \mbox{ in } \Omega \\[1ex]
     \hat u \ = \ -U     & \mbox{ on } \partial\Omega_D \\
     \hat v \ = \ U      & \mbox{ on } \partial\Omega_D \\
     \D\frac{\partial \hat u}{\partial \nu} \ = \
     \frac{\partial \hat v}{\partial \nu} \ = \ 0
                         & \mbox{ on } \partial\Omega_N
   \end{array} \right.  $$

Using {\em a priori estimates} of the solution of the Poisson equation,
we conclude that $\hat V$ depends continuously on the boundary data as
well as on the functions $\hat u$ and $\hat v$, appearing on the right
hand side of the PDE. Further, we know that the map $U \mapsto (\hat u,
\hat v)$ is well-defined. Therefore, we can deduce the well-definedness
of the application ${\cal T}_C$, for each doping profile $C$ in
$$ \{ C \in L^2(\Omega) ; \, C_m \le C(x) \le C_M,
   \mbox{ a.e. in } \Omega \} \, . $$
The continuity of ${\cal T}_C$ can be proved in an analogous way.
Furthermore, repeating the argumentation used for the operator
$\Sigma'(0)$, one can prove boundedness and compactness of the linear
operator ${\cal T}_C$.

To obtain the abstract formulation of the inverse problem related to
the operator ${\cal T}_C$, we take into account the more realistic case
of a finite number of measurements:

\begin{enumerate}
\item[1)] Parameter: \ $C = C(x) \ \in \ L^2(\Omega) =: \mathcal X$;
\vskip-5ex
\item[2)] Output: \ $Y = \big\{ {\cal T}_C(U_j) \big\}_{j=1}^N \in
\big[ L^2(\Gamma_1) \big]^N =: \mathcal Y$;
\vskip-5ex
\item[3)] Parameter-to-output map: \ $F: \mathcal X \to \mathcal Y$;
\end{enumerate}
for fixed $U_j \in H^{3/2}(\partial\Omega_D)$ with $\| U_j \|$ small.
The domain of definition of the operator $F$ is the same as in the case
of the V--C map. The parameter-to-output operator is defined by
$$  \begin{array}[t]{rcl}
    F: \ D(F) \subset \mathcal X  & \to &  \mathcal Y \\
    C  & \mapsto &  \big\{ {\cal T}_C U_j \big\}_{j=1}^N
    \end{array} $$

The well-definedness of the operator $F$ follows from the one of
${\cal T}_C$. The Fr\'echet-differentiability of the parameter-to-output
operator can be proved analogously as in the case of full voltage-current
data.

\subsection{Laser-beam-inducted current measurements} \label{ssec:up-lbic}

In this subsection we analyze the inverse problem related to the
laser-beam-inducted current (LBIC) image. This is a newly developed
non-destructive optical technique for the detection of semiconductor
properties. In this technique a laser beam is applied to the semiconductor
body to induce currents to flow through the ohmic contacts on the boundary.
The LBIC image consists of measurements of the local current $i(x)$ flowing
out through one contact ($\Gamma_1 \subset \partial\Omega_D$) induced by a
laser beam applied at location $x$ for all $x \in \Omega$. This image,
considered as a mapping $\Omega \ni x \mapsto i(x) \in \mathcal R$, is
considered to contain information about the doping profile. Therefore,
the relation between the LBIC image and the doping profile can be modeled
as an inverse problem for the system of drift-diffusion equations.

Let $(V^0, \hat u, \hat v)$ be the solution of system (\ref{eq:poiss-equil}),
(\ref{eq:fbip-lin2}). According to the representation formula
(\ref{eq:lbic-repr}), the LBIC image can be rewritten in the form
$$  i(x) \ = \ \tilde{v}(x) - \tilde{u}(x) \, ,\ x \in \Omega\, . $$

The equilibrium potential $V^0$ satisfying (\ref{eq:poiss-equil}) is
determined uniquely by the doping profile $C(x)$ and vice versa. Therefore,
reconstructing the doping profile $C(x)$ from the LBIC image is equivalent
to reconstructing the exponential of equilibrium potential $e^{V^0}$ from
the representation $i(x)$ of the LBIC image.

In \cite{FI1,FI2} the uniqueness of the inverse problem is analyzed. In
\cite{FI2} a one dimensional problem is considered and the non-uniqueness
of the inverse problem is proven. We shall next address this result.

A measurement $i(x)$ for $x \in \Omega$ is said to be attainable, if
$i(x) = \hat v(x) - \hat u(x)$ with $(\hat u, \hat v)$ being the solution
of (\ref{eq:fbip-lin2}) for some potential $V^0$. Let us for a moment
consider the one dimensional version of system (\ref{eq:fbip-lin2}) for
$\Omega = (0,1)$
\begin{equation} \label{eq:fbip-lin-1D}
   \left\{ \begin{array}{ll}
     (\mu_n e^{V^0} \tilde{u}')' \ = \ Q_0(V^0,x)
     (\tilde{u} - \tilde{v}) & x \in (0,1)  \\
     (\mu_p e^{-V^0} \hat \tilde{v}')' \ = \ Q_0(V^0,x)
     ( \tilde{v} - \tilde{u}) & x \in (0,1) \\[1ex]
     \tilde{u}(0) \ = \ \tilde{v}(0) \ = \ 1  &  \\
     \tilde{u}(1) \ = \ \tilde{v}(1) \ = \ 0  &
   \end{array} \right.
\end{equation}
In this case, the objective is to reconstruct $V^0$ (or alternatively
$e^{V^0}$) from given $i(x) = \tilde{u}(x) - \tilde{v}(x)$. The next result
establishes a necessary and sufficient condition for attainability of
a measurement $i(x)$.

\begin{propo} \mbox{\bf \cite[Theorem 2.1]{FI2}}
\label{prop:if2-21}
A measurement $i(x)$ is attainable if and only if there exists constants
$c_1$ and $c_2$ so that the equation
\begin{equation} \label{eq:fi-c1c2}
\left( c_1 - \frac{Q_0}{\mu_n} I(x) \right) Y(x) + i'(x) -
\left( c_2 + \frac{Q_0}{\mu_p} I(x) \right) Y(x)^{-1} \ = \ 0\, ,
\end{equation}
has a positive solution $Y(x)$ for each $x \in (0,1)$, and $Y(x)$ satisfies
the integral equation
$$  1 + \int_0^1 \left( c_1 - \frac{Q_0}{\mu_n} I(x) \right) Y(x) dx
    \ = \ 0\, . $$
Here $I(x) = \int_0^x i(\xi) d\xi$. Furthermore, if $i(x)$ is attainable
then the constants $c_1$ and $c_2$ are nonpositive.
\end{propo}
{\it Sketch of the proof:} \\
To prove the necessity, one integrates the differential equations in
(\ref{eq:fbip-lin-1D}) and obtain a representation for $i'(x)$. The
attainability of $i(x)$  follows from the fact that $Y(x) = e^{-V^0(x)}$
satisfies both the quadratic equation and the integral equation of the
proposition.

To prove the sufficiency, one sets
\begin{eqnarray*}
  \hat u(x) & = & 1 + \int_0^x \left(c_1 - \frac{Q_0}{\mu_n} I(\xi)\right)
  Y(\xi) d\xi \\
  \hat v(x) & = & 1 + \int_0^x \left(c_2 - \frac{Q_0}{\mu_p} I(\xi)\right)
  Y(\xi)^{-1} d\xi
\end{eqnarray*}
and obtain in a straightforward way that $(\hat u, \hat v)$ solve
(\ref{eq:fbip-lin-1D}) for $V^0(x) = - \ln Y(x)$. From an obvious
substitution follows $\hat v - \hat u = \int_0^x i'(\xi) d\xi$.
\hfill $\Box$

\bigskip

According to this result, the attainability of a measurement $i(x)$ is
equivalent to the solvability of a quadratic equation for $Y(x)$. Notice
that, if the constants $c_1$ and $c_2$ are known, then the potential is
obtained simply by $V^0(x) = - \ln Y(x)$.

A first identifiability result is given in \cite{FI1}. In this paper, the
authors prove that $i(x) \equiv 0$ if and only if $V^0(x) \equiv c$, for
some real constant $c$ (see Theorem~3.2 in the reference above). Therefore,
in general there is no uniqueness for the inverse problem. Using 
Proposition~\ref{prop:if2-21}, the same authors manage to extend this
first non-uniqueness result for the one-dimensional case presented above,
as follows

\begin{propo} \mbox{\bf \cite[Theorem 2.3]{FI2}}
\label{prop:if2-23}
Let $i(x) \in C^1_0[0,1]$ be an attainable measurement and $V^0$ be the
corresponding potential. Moreover, assume that the constants $c_1$ and
$c_2$ found by Proposition~\ref{prop:if2-21} with respect to $Y(x) =
e^{-V^0(x)}$ satisfy
$$  c_1 \ < \ \frac{Q_0}{\mu_n} I_{\rm min} \ \ \ and \ \ \
    c_2 \ < \ -\frac{Q_0}{\mu_p} I_{\rm max} \, , $$
where $I_{\rm min}$ and $I_{\rm max}$ are respectively the maximum and the
minimum of $I(x)$ in $[0,1]$. Then there is a one-parameter family
$\{ V(x)\}$, containing $V^0$ and strictly monotone in the parameter, that
produces the same measurement $i(x) = \hat v(x) - \hat u(x)$ from system
(\ref{eq:fbip-lin-1D}).
\end{propo}
%
%{\it Sketch of the proof:} \\
%The proof is based on the necessary and sufficient condition given in
%Proposition~\ref{prop:if2-21}. Solving Equation (\ref{eq:fi-c1c2}) for
%$Y$, one obtains
%
%$$ Y(x;c_1,c_2) \ = \ \frac{-i'(x)-\sqrt{y(x;c_1,c_2)}}
%                      {2(c_1 - Q_0 \mu_n^{-1} I(x))} \, , $$
%
%where
%$$ y(x;c_1,c_2) \ := \ i'(x)^2 + 4
%   \left(c_1 - \frac{Q_0}{\mu_n} I(x)\right)
%   \left(c_2 - \frac{Q_0}{\mu_p} I(x)\right) \, . $$
%
%Using this representation for $Y$, the integral equation in Proposition%
%~\ref{prop:if2-21} can be written as
%
%$$  F(c_1,c_2) \ := \ \frac{1}{2} \int_0^1 \sqrt{y(x;c_1,c_2)}\ dx - 1
%    \ \equiv \ 0 \, . $$
%
%Since $\partial F/\partial c_1$ and $\partial F/\partial c_2$ are both
%negative, the implicit function theorem can be used to obtain the
%representation $c_1 = f(c_2)$, where $f$ is given by $F(c_1,c_2) = 0$.
%The one parameter family of functions can be obtained by
%
%$$ V(x;c_2) \ := \ - \ln Y(x; f(c_2), c_2) \, . $$
%
%The monotony of $V$ (with respect to $c_2$) follows from the monotony of $f$,
%which is a consequence of the implicit function theorem. \hfill $\Box$

\bigskip

Proposition~\ref{prop:if2-23} characterizes the nonuniqueness of the
one-dimensional inverse problem for the LBIC operator. Therefore, more
information is needed to possibly recover $V^0(x)$ from $i(x)$ uniquely.
In the LBIC technique it is reasonable to assume that doping profile is
known on the boundary where the ohmic contacts are made ($\partial 
\Omega_D \subset \partial\Omega$). 
Thus, if we assume that $V^0$ is given at $x=0$, we gain another constraint
for $c_1$ and $c_2$, namely
\begin{equation} \label{eq:fi-constr}
  c_1 e^{-V^0(0)} - i'(0) - c_2 e^{V^0(0)} \ = \ 0 \, .
\end{equation}
This additional constraint ensures the unique recovery of $V^0(x)$ from
the LBIC image $i(x)$ among the monotone one-parameter family $\{ V(x)\}$
described in Proposition~\ref{prop:if2-23}. Notice that this class of
potentials does not include all possible solutions to the inverse
problem. Therefore, the above constraint does not lead to uniqueness
of solutions of the one-dimensional inverse problem in general.

In \cite{FI2} the authors also propose an algorithm for the reconstruction
of $V^0(x)$ from $i(x)$ based on Proposition~\ref{prop:if2-23} and the
additional constraint (\ref{eq:fi-constr}). As discussed above, it is enough
to reconstruct the constants $c_1$ and $c_2$. The proposed algorithm consists
in a Gauss-Newton method for the minimization of a least square functional
$J$ associated to the residual of the pairs $(c_1, c_2)$ in both
(\ref{eq:fi-c1c2}) and (\ref{eq:fi-constr}), namely
$$ J(c_1,c_2) \ := \ \frac{1}{2} \big( J_1(c_1,c_2)^2 +
                   J_2(c_1,c_2)^2 \big)\, , $$
where
\begin{eqnarray*}
 J_1(c_1,c_2) & := & \frac{1}{2} \int_0^1 \sqrt{ i'(x)^2 + 4
\left(c_1 - \frac{Q_0}{\mu_n} I(\xi)\right)
\left(c_2 - \frac{Q_0}{\mu_p} I(\xi)\right) }\ dx - 1 \\[1ex]
J_2(c_1,c_2) & := & \Big( c_1 e^{-V^0(0)} - i'(0) - c_2 e^{V^0(0)} \Big)
   e^{-|V^0(0)|} \, .
\end{eqnarray*}

In \cite{FIR}, a similar model based on the drift diffusion equations is
used in order to analyze several parameter identification problems
for semiconductor diodes by LBIC imaging. Numerical methods are developed
for the simulation of the LBIC images of a diode as well as for the
identification of parameters (junction depth, diffusion length equilibrium
potential) from the LBIC image by least-squares formulation.

\subsection{Inverse doping profile: Identification}
\label{ssec:up-id}

In this subsection we consider the identification question related to the
inverse doping profile problem, i.e. we shall focus on the following
fundamental issue concerning the parameter identification problems:

\begin{quote}
{\em Is the available data enough to determine uniquely the doping profile,
or (alternatively) which set of data is sufficient to determine uniquely
the doping profile?}
\end{quote}

In the one-dimensional case (i.e. $\Omega = (0,L)$) the identification
problem was considered in \cite{BEM}. One can assume that the voltage is
applied at $x=0$ and the measurements of both current and capacitance
are taken at $x=L$. Therefore, a single measurement (reduced data) consists
of two real numbers and full data, in this case, correspond to measure
the current and/or the capacitance as a function of the applied voltage
$U \in (-r,r)$, with appropriate $r \in \mathcal R$.

Arguing with the dimensionality of the parameter and data spaces as well as
with structural properties of the operators $\Sigma_C$ and $\mathcal T_C$,
the authors are able to fully analyze the one-dimensional inverse doping
profile problem. The corresponding results are summarized in the following
proposition.

\begin{propo}
Let us consider the inverse doping profile for system (\ref{eq:dd-sys1})--%
(\ref{eq:dd-sys9}) at the one-dimensional dimensional domain $\Omega = (0,L)$.
The following assertion hold:
\begin{enumerate}
\item If one has access only to restricted data, even if it is possible to
measure both, current and capacitance, the data are not sufficient to identify
the doping profile;
\item If one has access to full data, it is not possible to uniquely identify
the doping profile neither from current measurements nor from capacitance
measurements.
\end{enumerate}
\end{propo}

One should notice that the doping profile $C = C(x)$ in this case is a
function of a one-dimensional space variable.

In the same paper, the authors also consider the transient case of the
one-dimensional inverse doping profile problem. They prove, under special
assumptions, that if both current and capacitance measurement are
available, then the doping profile can be uniquely reconstructed from
the data. Since we consider only the stationary drift-diffusion system
in this paper, we shall not investigate this result here. For details,
we refer to \cite[Theorem~3]{BEM}.

In the one-dimensional case, the special problem in which the doping profile
is a piecewise constant function of position is treated in \cite{BEM}. In
this very particular case, the domain $\Omega$ can be split as $\bar\Omega
= \bar\Omega_n \cup \bar\Omega_p$, such that $C(x) \equiv C_+$ in $\Omega_n$
and $C(x) \equiv C_-$ in $\Omega_p$. This problem is also known as {\em
identification of P-N junctions}. The authors prove that reduced current
data suffice to uniquely identify the exact location of the P-N junctions
(i.e., $\bar\Omega_n \cap \bar\Omega_p$) if the number of junctions is lower
or equal to two (see \cite[Theorem~4]{BEM}).

The two-dimensional case is considered in \cite{BEMP,BEM}. Particularly
interesting is the inverse problem related to the V--C map for the
linearized unipolar case close to equilibrium (see Subsection~%
\ref{ssec:ubc}), which can be directly related to the inverse problem in electrical
impedance tomography.

The inverse problem in the unipolar case corresponds to the determination
of the doping profile $C$ from the map
\begin{equation} \label{eq:sigma_C}
\Sigma_C'(0):
   \begin{array}[t]{rcl}
     H^{3/2}(\partial\Omega_D) & \to & H^{1/2}(\Gamma_1) \\
     U & \mapsto & (\hat{J}_n\cdot\nu) |_{\Gamma_1}
   \end{array}
\end{equation}
where $(u,V^0)$ is the solution of the system in (\ref{eq:unipolar}).

As already observed in Subsection~\ref{ssec:ubc}, it follows from the
fact that $V^0|_{\partial\Omega_D} = V_{\rm bi}$ is a known function,
that the current data $J_n\cdot\nu |_{\Gamma_1}$ can be directly
substituted by the Neumann data $u_\nu |_{\Gamma_1}$.
Therefore, the inverse problem can be divided in 2 steps:
\begin{enumerate}
\item[1)] Define $\gamma := e^{V^0}$ and identify $\gamma$ in
$$ \left\{ \begin{array}{ll}
     {\rm div} (\gamma \nabla u) \ = \ 0  & {\rm in}\ \Omega \\
     u \ = \ U                         & {\rm on}\ \partial\Omega_D \\
     u_\nu \ = \ 0              & {\rm on}\ \partial\Omega_N
  \end{array} \right. $$
from the Dirichlet-to-Neumann map: \ $u|_{\partial\Omega_D} \mapsto
u_\nu |_{\Gamma_1};$
\item[2)] Obtain the doping profile $C(x)$ from: \ $C = \gamma -
\lambda^2 \Delta \, (\ln \gamma)$.
\end{enumerate}

The identification problem in 1) corresponds to the {\em electrical
impedance tomography} (or {\em inverse conductivity problem}) in elliptic
equations with mixed boundary data. For the case of the full
Dirichlet-to-Neumann operator, i.e. $\Gamma_1 = \partial\Omega_D =
\partial\Omega$, this inverse problem has been intensively analyzed
in the literature over the last fifteen years. Using different
regularity assumptions on the conductivity $\gamma$, many authors
proved that the coefficient $\gamma(x)$ of the elliptic equation
$\nabla \cdot (\gamma\nabla u) = 0$ is uniquely determined by the
corresponding Dirichlet-to-Neumann map on the boundary (a historical
overview can be found in \cite{Bo}).

In the sequel we mention a result due to A. Nachman for two-dimensional
domains. The proof of this theorem gives a constructive procedure for
recovering $\gamma$ from the Dirichlet-to-Neumann map.

\begin{propo} {\bf \cite[Theorem~1]{Na}}
\label{pr:Na-1}
Let $\Omega$ be bounded and Lipschitz. Further, let \ $\gamma_i \in
L^\infty(\Omega) \cap W^{2,p}(\Omega)$, $i=1,2$, for some $p > 1$ with
positive lower bound. Then, the equality of the Dirichlet-to-Neumann maps
$$ \Lambda_i : \begin{array}[t]{rcl}
H^{1/2}(\partial\Omega) & \to     & H^{-1/2}(\partial\Omega) \\
u                       & \mapsto & u_\nu
\end{array} $$
for the solutions of \ ${\rm div} (\gamma_i \nabla u) = 0$, \ 
implies $\gamma_1 = \gamma_2$.
\end{propo}

According to Proposition~\ref{prop:Ma-331}, $H^2$-regularity of the
solution $(V,u,v)$ of system (\ref{eq:dd-sys1})--(\ref{eq:dd-sys9})
can be obtained under stronger regularity assumptions on both the mixed
boundary conditions, and the domain. Using this regularity result, it is
possible to adapt Proposition~\ref{pr:Na-1} for the identification problem
in the unipolar case for the operator $\Sigma_C'(0): H^{3/2}(\partial\Omega_D)
\to H^{1/2}(\Gamma_1)$ in the idealized case $\Gamma_1 = \partial\Omega_D
= \partial\Omega$, as follows

\begin{propo} {\bf \cite[Theorem~4.2]{BEMP}} \label{pr:ident-up}
Let $\Omega \subset \mathbb{R}^2$ be bounded and Lipschitz. Further,
let $\Gamma_1 = \partial\Omega_D = \partial\Omega$. Then, given two
doping profiles $C_1, C_2 \in D(F)$, the equality \ $\Sigma'_{C_1}(0) =
\Sigma'_{C_2}(0)$ implies \ $C_1 = C_2$.
\end{propo}

If we consider the solution of (\ref{eq:dd-sys1})--(\ref{eq:dd-sys9})
to be only in $H^1$ (see Proposition~\ref{prop:MRS-3316}), we can
alternatively consider the following identifiability result from Brown
and Uhlmann for $W^{1,p}(\Omega)$, $p > 2$, conductivities:

\begin{propo} {\bf \cite[Theorem~4.1]{BU}}
Let $\Omega$ be bounded and Lipschitz. Further, let $\gamma_1$ and
$\gamma_2$ be two conductivities with $\nabla\gamma_i$ in $L^p(\Omega)$,
$p>2$. Then, the equality of the Dirichlet-to-Neumann maps
$$ \Lambda_i :
   \begin{array}[t]{rcl}
     H^{1/2}(\partial\Omega) & \to     & H^{-1/2}(\partial\Omega) \\
     u                       & \mapsto & u_\nu
   \end{array} $$
for the solutions of \ ${\rm div} (\gamma_i \nabla u) = 0$, \
implies $\gamma_1 = \gamma_2$.
\end{propo}

Using this identifiability result, it is
possible to deduce, for the operator $\Sigma'_{C}(0):
H^{1/2}(\partial\Omega_D) \to H^{-1/2}(\Gamma_1)$ an analog
result to the one presented in Proposition~\ref{pr:ident-up}.
Notice that this result is particularly interesting for the case
of zero space charge (see Subsection~\ref{ssec:vcmap}) for the
V--C map, allowing to prove identifiability of doping profiles
$C \in L^\infty(\Omega) \cap W^{1,p}(\Omega)$.

%---------------------------------------------------------------------------
\section{Numerical experiments} \label{sec:num}

In this section we derive a numerical method to identify the doping profile
in the linearized unipolar case close to equilibrium (\ref{eq:unipolar}).
In this particular case, due to the assumptions $p \equiv 0$ and $Q
\equiv 0$, the Poisson equation and the continuity equation for the electron
density $n$
decouple, and we have to identify $C = C(x)$ in
$$
\left\{ \hskip-0.2cm \begin{array}{rcl@{\ }l}
   \lambda^2 \, \Delta V^0 & = & e^{V^0} - C(x) & {\rm in}\ \Omega \\
   V^0 & = & V_{\rm bi}(x) & {\rm on}\ \partial\Omega_D \\
   \nabla V^0 \cdot \nu & = & 0 & {\rm on}\ \partial\Omega_N \\
\end{array} \right.
\hskip0.4cm
\left\{ \hskip-0.2cm \begin{array}{rcl@{\ }l}
   {\rm div}\, (e^{V^0} \nabla \hat{u}) & = & 0 &  {\rm in}\ \Omega \\
   \hat{u} & = & U(x) & {\rm on}\ \partial\Omega_D \\
   \nabla \hat{u} \cdot \nu & = & 0 & {\rm on}\ \partial\Omega_N \, .
\end{array} \right.
$$
Notice that, due to the relation $\hat{J} \cdot \nu = \mu_n e^{V^0} \hat{u}$, the
Neumann boundary condition $\hat{J}_n \cdot \nu |_{\partial\Omega_N} = 0$
in (\ref{eq:unipolar}) can be substituted by $\nabla \hat{u} \cdot \nu = 0$ on
$\partial\Omega_N$. As already observed in Subsection~\ref{ssec:up-id},
we can write $\gamma(x) := e^{V^0(x)}$; solve the parameter identification
problem
\begin{equation} \label{eq:num-d2n}
\left\{ \begin{array}{rcll}
   {\rm div}\, (\gamma \nabla \hat{u}) & = & 0 &  {\rm in}\ \Omega \\
   \hat{u} & = & U(x) & {\rm on}\ \Omega_D \\
   \nabla \hat{u} \cdot \nu & = & 0 & {\rm on}\ \Omega_N \, ,
\end{array} \right.
\end{equation}
for the function $\gamma$; and finally evaluate
$$ C(x) \ = \ \gamma - \lambda^2 \, \Delta (\ln \gamma) \, . $$
Since the evaluation of $C$ from $\gamma$ can be explicitely performed
(a direct problem) and is a well posed procedure, we shall focus on
the problem of identifying the function parameter $\gamma$ in
(\ref{eq:num-d2n}). Therefore, the inverse problem of identifying the
doping profile $C(x)$ in the linearized unipolar model (\ref{eq:unipolar})
corresponds to the identification of $\gamma(x)$ in (\ref{eq:num-d2n})
from the Dirichlet to Neumann (DtN) map
$$  \Lambda_\gamma : \begin{array}[t]{rcl}
    H^{3/2}(\partial\Omega_D) & \to & H^{1/2}(\Gamma_1) \\
    U & \mapsto & \gamma\, \D\frac{\partial u}{\partial\nu} \Big|_{\Gamma_1}
    \end{array} $$
As we saw in Subsection~\ref{ssec:up-id}, the DtN operator is given by
$\Lambda_\gamma = \Sigma_C'(0)$.

Notice that, due to the nature of the physical problem related to the
drift-diffusion equations, we can consider as {\em inputs} for the DtN
map only functions of the type:
$$ U \ = \ \left\{ \begin{array}{rl}
           \tilde U, & \mbox{ at }\ \partial\Omega_D \backslash \Gamma_1 \\
           0, & \mbox{ at }\ \Gamma_1 \end{array} \right. . $$
Furthermore, the {\em outputs} or measurements are only available at
$\Gamma_1$. This is the basic difference between the parameter
identification problem in (\ref{eq:num-d2n}) and the inverse problem
in {\em electrical impedance tomography}, i.e. the fact that both
Dirichlet (input) and Neumann (output) are prescribed only at specific
parts of the boundary. For this special inverse problem (with mixed
boundary data) there are so far no analytical results concerning
identifiability and the few numerical results in the literature are
those discussed in \cite{BEMP, BEM, FIR}.

We shall work with a reduced set of data, as described in Subsection~%
\ref{ssec:vcmap}, i.e. within the following framework:
\begin{enumerate}
\item[1)] Parameter: \ $\gamma = \gamma(x) \ \in \ H^2(\Omega) =: \mathcal X$;
\vskip-5ex
\item[2)] Input (fixed): \ $U_j \in H^{1/2}(\partial\Omega_D)$, \
$U_j |_{\Gamma_1} = 0$, \ $1 \le j \le N$;
\item[3)] Output (data): \ $Y = \big\{ \gamma
\frac{\partial \hat{u}_j}{\partial \nu} |_{\Gamma_1} \big\}_{j=1}^N \in
[L^2(\Gamma_1)]^N =: \mathcal Y$; \\
(here $u_j$ is the solution of (\ref{eq:num-d2n}) for $U = U_j$)
\vskip-5ex
\item[4)] Parameter-to-output map: \ $F: \mathcal X \to \mathcal Y$.
\end{enumerate}
The domain of definition of the operator $F$ is
$$ D(F) := \{ \gamma \in H^2(\Omega) ; \,  \gamma(x) \ge \gamma_- > 0,
   \mbox{ in } \Omega \} \, , $$
where $\gamma_-$ is an appropriate positive constant.
We shall denote the noisy data by $Y^\delta$ and assume that the data error
is bounded by
$$ \| Y - Y^\delta \| \ \le \ \delta\, . $$
Thus, we are able to represent the inverse doping problem in the abstract
form
\begin{equation} \label{eq:ip-dd}
F(\gamma) \ = \ Y^\delta \, .
\end{equation}

A common technique to solve the inverse problem in (\ref{eq:ip-dd}) is
the output least-square family of methods. Basically, all output
least-square methods minimize iteratively the residual functional
related to (\ref{eq:ip-dd}) with some Newton-type method \cite{Ba,Ba1,Ba2,
BG,BG1,DES,EHN,Ha,Ho,Ka,Ka1,Sc}. In the literature, one can find several
applications of such methods for the {\em electrical impedance tomography
problem} (see, e.g., \cite{Bo1,DL,Do,YWT}).

A simple and robust iterative method to solve the problem in (\ref{eq:ip-dd})
is the so called {\em Landweber iteration} \cite{DES,EHN,ES,HNS}, in which
the k-step is described by
$$ \gamma^\delta_{k+1} \ = \ \gamma_k^\delta - F'(\gamma_k^\delta)^*
   \big( F(\gamma_k^\delta) - Y^\delta \big)\, . $$
This iteration is known to generate a regularization method for the inverse
problem, the stopping index playing the rule of the regularization parameter
(for regularization methods see, e.g., \cite{EHN,EKN,ES,Mo,TA}).

For our numerical experiments, we propose an iterative method of adjoint
type in order to solve the identification problem (\ref{eq:num-d2n}), the
so called {\em Landweber--Kaczmarz method}. This method derives from the
coupling of the strategies of the Landweber iteration and the Kaczmarz
method. The Kaczmarz method is a fixed point algorithm which has been
proven to be efficient for solving inverse problems in Tomography
\cite{Bau,BS,Gro,Ki,Nt}. For a detailed analysis of the method we refer to
\cite{Bau,Me} for the finite dimensional setting and to \cite{BR,Mc,Mc1}
for the infinite dimensional setting.

A detailed analysis of the Landweber--Kaczmarz method can be found in
\cite{KS}. It is worth mentioning that this method has already been
successfully applied to the {\em electrical impedance tomography} problem
\cite{Na}. To formulate the method, we need first define the
parameter-to-output maps
\begin{enumerate}
\item[4')] Operators \ $ \mathcal F_j$ for $j = 1, \dots, N$:
$$ \mathcal F_j: \begin{array}[t]{rcl}
   H^2(\Omega) & \to & L^2(\Gamma_1) \\ \gamma & \mapsto & \gamma \frac{\partial \hat{u}_j}{\partial \nu}|_{\Gamma_1}
   \end{array} $$
\end{enumerate}
Now, setting $Y_j := \mathcal F_j(\gamma)$ for $1 \le j \le N$, the
Landweber--Kaczmarz iteration can be written as:
\begin{equation} \label{eq:land-kacz}
  \gamma^\delta_{k+1} \ = \ \gamma_k^\delta
                          - \mathcal F_k'(\gamma_k^\delta)^*
  \big( \mathcal F_k (\gamma_k^\delta) - Y^\delta_k \big)\, ,
\end{equation}
for $k = 1, 2, \dots$, where we adopted the notation
$$  \mathcal F_k := \mathcal F_j, \ \ Y^\delta_k := Y^\delta_j,
    \ \ {\rm with} \ \ k = i \cdot N + j, \ \ {\rm and}\ \
    \left\{ \begin{array}{l}
    i = 0, 1, \dots \\ j = 1, \dots, N
    \end{array} \right. \, . $$

Each step of the Landweber--Kaczmarz method consists in one step of the
Landweber iteration with respect to the $j$-th component of the residual
$F(\gamma) - Y$. These steps are performed in a cyclic way for each one of
the residual components $\mathcal F_j(\gamma) - Y_j$, $j=1,\cdots,N$.

As far as the implementation of the method is concerned, it is enough to describe
the general step of the Landweber iteration. The variational
formulation of the iterative step in (\ref{eq:land-kacz}) reads
\begin{equation} \label{eq:land-cacz-var}
 \ipl \gamma_{k+1} - \gamma_k, \ h \ipr_{L^2(\Omega)} \ = \ -
   \ipl \mathcal F'_k(\gamma_k) h , \ \mathcal F_k(\gamma_k) - Y_k \ipr
   _{L^2(\Omega)} \, ,
\end{equation}
where $h \in H^1(\Omega)$ is a test function (to simplify the notation
we set $\delta = 0$, i.e. $Y^\delta_k = Y_k$ and $\gamma^\delta_k =
\gamma_k$).

In order to compute the inner product on the right hand side of
(\ref{eq:land-cacz-var}), we use the identity:
\begin{equation} \label{eq:adj-form}
  \ipl \mathcal F'(\gamma) h , z \ipr_{L^2(\Gamma_1)} \ = \
  \int_{\Omega} h \, \nabla G(\gamma) \cdot \nabla \Phi(\gamma) \, dx ,
\end{equation}
for $z \in L^2(\Gamma_1)$, where the $H^1(\Omega)$-function $\Phi(a)$ solves
$$
 \left\{ \begin{array}{rl}
     - \nabla( a(x) \nabla w) \ = \ 0 , & {\rm in} \ \Omega \\
     w \ = \ z , & {\rm on} \ \Gamma_1 \\
     w \ = \ 0 , & {\rm on} \ \partial\Omega / \Gamma_1
   \end{array} \right.
$$
and the $H^1(\Omega)$-function $G(a)$ solves
$$
 \left\{ \begin{array}{rl}
     - \nabla( a(x) \nabla w) \ = \ 0 , & {\rm in} \ \Omega \\
     a(x) w_\nu \ = \ 0 , & {\rm on} \ \partial\Omega_N \\
     w \ = \ g , & {\rm on} \ \partial\Omega_D
   \end{array} \right.
$$
Indeed, since the Fr\'echet derivative of the operator
$$ \begin{array}[t]{rcl}
     \Psi: H^2(\Omega) & \to & H^{1/2}(\partial\Omega) \\
     a & \mapsto &  a w_\nu |_{\partial\Omega}
   \end{array}
\ \ \ \ {\rm where}\ \ \ \
   \left\{ \begin{array}{rl}
     - \nabla( a(x) \nabla w) \ = \ f , & {\rm in} \ \Omega \\
     w \ = \ g , & {\rm on} \ \partial\Omega
   \end{array} \right.
 $$
in the direction $h \in H^2(\Omega)$ is given by
$$ \Psi'(a) \cdot h \ = \ (h G_\nu(a) + a \psi_\nu)\, , $$
where
$$ \left\{ \begin{array}{rcll}
     - \nabla( a(x) \nabla \psi) & = & \nabla(h(x) \nabla G(a)), & {\rm in}
     \ \Omega \\
     \psi & = & 0 , & {\rm on} \ \partial\Omega
   \end{array} \right.  $$
we have
\begin{eqnarray*}
\ipl \mathcal F'(\gamma) h, z \ipr_{L^2(\Gamma_1)} & = &
\int_{\Gamma_1} z \, \big( h (G(\gamma))_\nu + \gamma \psi_\nu  \big) \\
& \hskip-4cm = & \hskip-2cm
  \int\limits_{\Gamma_1} z \, h \, (G(\gamma))_\nu
  + \int\limits_{\Gamma_1} \Phi(\gamma) \gamma \psi_\nu
  + \int\limits_{\partial\Omega_D/\Gamma_1} \Phi(\gamma) \gamma \psi_\nu
  + \int\limits_{\partial\Omega_N} \Phi(\gamma) \gamma \psi_\nu \\
& \hskip-4cm = & \hskip-2cm
  \int_{\Gamma_1} z \, h \, (G(\gamma))_\nu
  + \int_{\Omega} \nabla( \gamma \nabla \psi) \Phi(\gamma)
  + \int_{\Omega} \gamma \, \nabla \psi \cdot \nabla \Phi(\gamma) \\
& \hskip-4cm = & \hskip-2cm
  \int_{\Gamma_1} z \, h \, (G(\gamma))_\nu
  - \int_{\Omega} \nabla(h \nabla G(\gamma)) \Phi(\gamma)
  + \int_{\partial\Omega} \psi \big( \gamma (\Phi(\gamma))_\nu \big) \\
& \hskip-4cm \ & \hskip-2cm
  - \int_{\Omega} \psi \nabla \big( \gamma \nabla \Phi(\gamma) \big) \\
& \hskip-4cm = & \hskip-2cm
  \int_{\Gamma_1} z \, h \, (G(\gamma))_\nu
  - \left[ \int_{\Gamma_1} h \, (G(\gamma))_\nu \Phi(\gamma)
         + \int_{\partial\Omega / \Gamma_1}
         h \, (G(\gamma))_\nu \Phi(\gamma) \right] \\
& \hskip-4cm \ & \hskip-2cm
  + \int_{\Omega} h \, \nabla G(\gamma) \cdot \nabla \Phi(\gamma) \\
\end{eqnarray*}
and (\ref{eq:adj-form}) follows. Therefore, the term on the right hand
side of (\ref{eq:land-cacz-var})can be evaluated by using formula
(\ref{eq:adj-form}) with $z = \mathcal F_k(\gamma_k) - Y_k$.

\begin{figure}
\centerline{ \epsfysize4cm \epsfbox{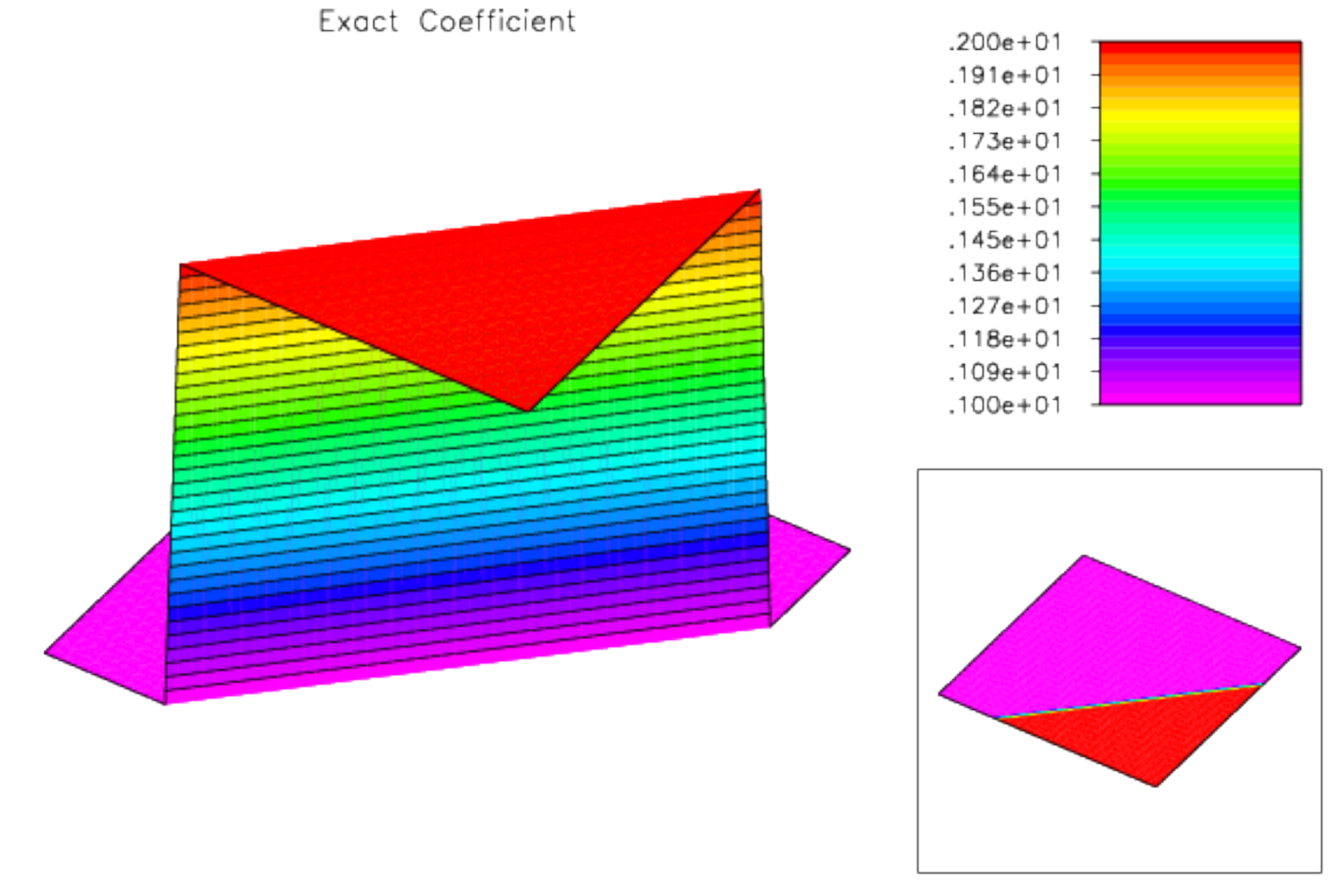} \hfill
             \epsfysize4cm \epsfbox{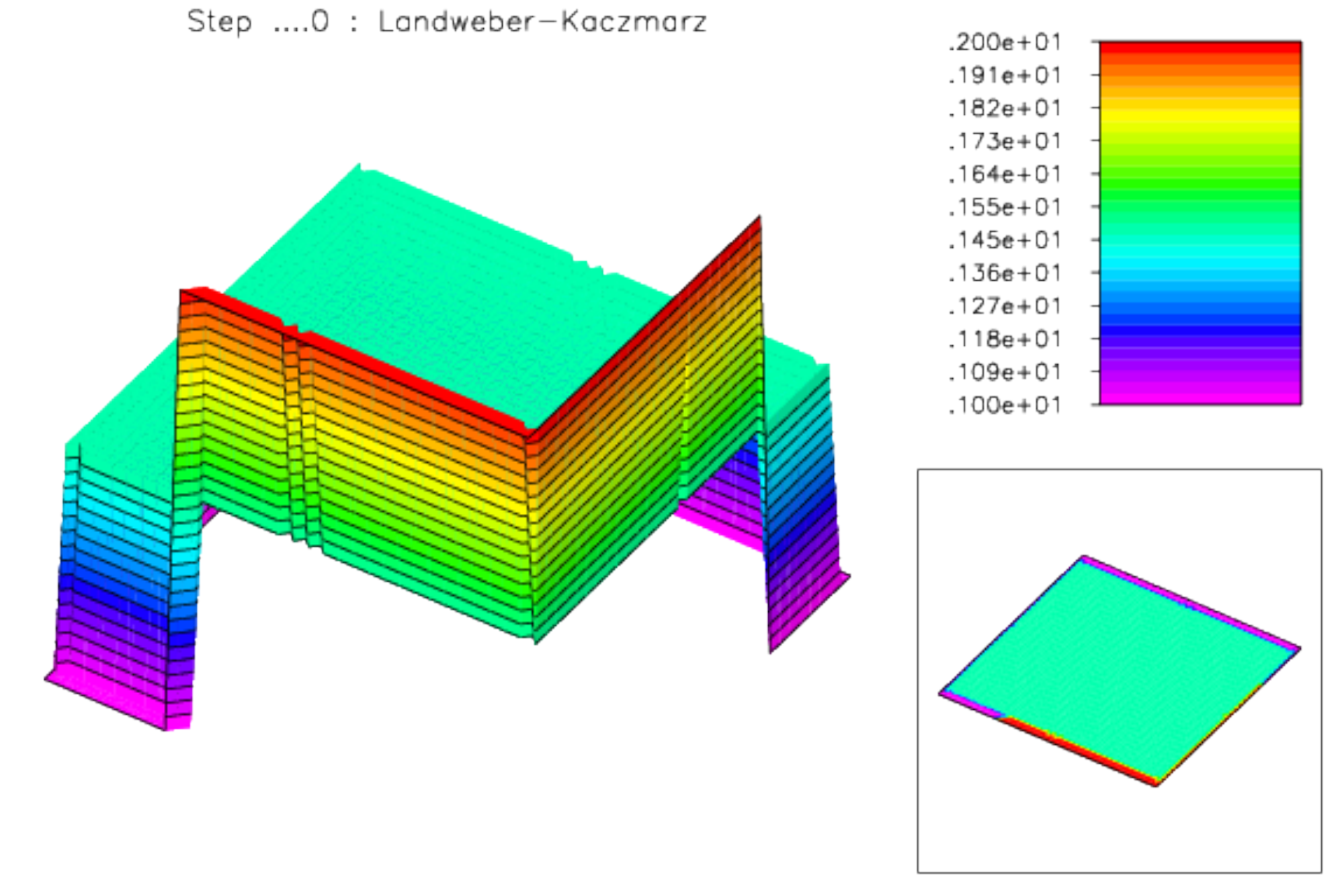}  }
\centerline{\hfil (a) \hskip5cm  (b) \hfil}
\caption{Picture (a) shows the exact coefficient $\bar\gamma(x)$ to
be reconstructed. On picture (b), the initial condition for the
Landweber-Kaczmarz iteration is shown.} \label{fig:exsol-incond}
\end{figure}

For the concrete numerical test performed in this paper, $\Omega \subset
\mathcal R^2$ is the unit square, and the boundary parts are defined as
follows
$$ \Gamma_1 \ := \  \{ (x,1) \, ;\ x \in (0,\T\frac{1}{2}) \} \, ,\ \ \
   \partial\Omega_D \ := \ \Gamma_1 \cup \{ (x,0) \, ;\ x \in (0,1) \} $$
$$ \partial\Omega_N \ := \ \{ (0,y) \, ;\ y \in (0,1) \} \cup 
   \{ (1,y) \, ;\ y \in (0,1) \} \cup
   \{ (x,1) \, ;\ x \in (\T\frac{1}{2},1) \} \, . $$
The fixed inputs $U_j \in H^{1/2}(\partial\Omega_D)$, are chosen to be
piecewise linear functions supported in $\partial\Omega_D / \Gamma_1$
$$  U_j(x) \ := \ \left\{ \begin{array}{rl}
      1 - \frac{1}{h} |x - x_j| , & |x - x_j| \le h \\
      0, & {\rm else} \end{array} \right. $$
where the points $x_j$ are equally spaced in the interval $(0,1)$. The
doping profile $C = C(x)$ to be reconstructed corresponds to the function
$\bar\gamma(x)$ shown in Figure~\ref{fig:exsol-incond}~(a). In this
figure, as well as in the forthcoming ones, $\Gamma_1$ appears in the
lower right part of the picture and $\partial\Omega_D / \Gamma_1$ appears
on the top (the origin corresponds to the upper right corner).

To generate the problem data, one has to solve the direct problem in
(\ref{eq:num-d2n}) for each input function $U_j$, $j=1,\cdots,N$. In
order to avoid the so called {\em inverse crimes}, these problems are
solved using adaptive mesh regularization and a piecewise linear finite element base
with approximately 8000 nodal points. This mesh is different from the
one used to solve the mixed elliptic boundary value problems, related
to the implementation of the Landweber--Kaczmarz method. These problems
are solved using a multigrid finite element method at uniformly refined
grids with approximately 2000 nodal points.

We still have to take into account an important issue concerning the
stability of the numerical implementation. Due to the particular geometry
of $\Omega$ (note that $\partial\Omega_D$ and $\partial\Omega_N$ meet at
angles of $\pi$ and $\pi/2$), both the solution of the direct elliptic
(mixed) problems as well as the solution of the boundary value problems
involved in the implementation of the Landweber--Kaczmarz method are not
in $H^2(\Omega)$ (see remark at the end of Subsection~\ref{ssec:ex-uniq}).

Because of this lack of regularity in the solution of the elliptic
boundary value problems, the numeric implementation of the
Landweber--Kaczmarz method has shown to be very unstable. After a
few iterative steps the sequence $\gamma_k$ became unbounded, the main
singularity appearing near the boundary (note that we assume $C$, or
equivalently $\gamma$, to be known at the boundary) close to the
contact points between $\partial\Omega_D$ and $\partial\Omega_N$.
This phenomena could be observed even if we started the iteration
with $\gamma_0(x) = \gamma^\dag(x)$, the exact solution of the inverse
problem.

In order to avoid the instability described above, we make the additional
assumption that the doping profile is known in a thin strip close to
$\partial\Omega$. Therefore, we only have to reconstruct the values
of $\gamma(x)$ at a subdomain $\tilde\Omega \subset\subset \Omega$.
With this extra assumption, the numerical implementation becomes stable
and we are able to characterize (numerically) the exact solution
$\gamma^\dag(x)$ as a fixed point of the Landweber--Kaczmarz iteration.
It is worth mentioning that this sort of assumption is very common in
the literature (see, e.g., \cite{Bo} and the references therein) and has
been used since the early investigations of the electrical impedance
tomography, in order to insure extra regularity for both numerical and
analytical approaches (see \cite{Sc1}).

In Figure~\ref{fig:exsol-incond}~(b) the initial condition for the
Landweber--Kaczmarz method is shown. Comparing the initial condition
with the exact solution, one can observe that the values of $\gamma_0(x)$
and $\gamma^\dag(x)$ coincide close to $\partial\Omega$. This is in
accordance with the assumption above. Close to the boundary
$\partial\Omega$, the values of $\gamma_k$ are not iterated, and we
actually have $\gamma_k = \gamma^\dag$ at $\Omega / \tilde\Omega$.

\begin{figure}
\centerline{ \epsfysize4cm \epsfbox{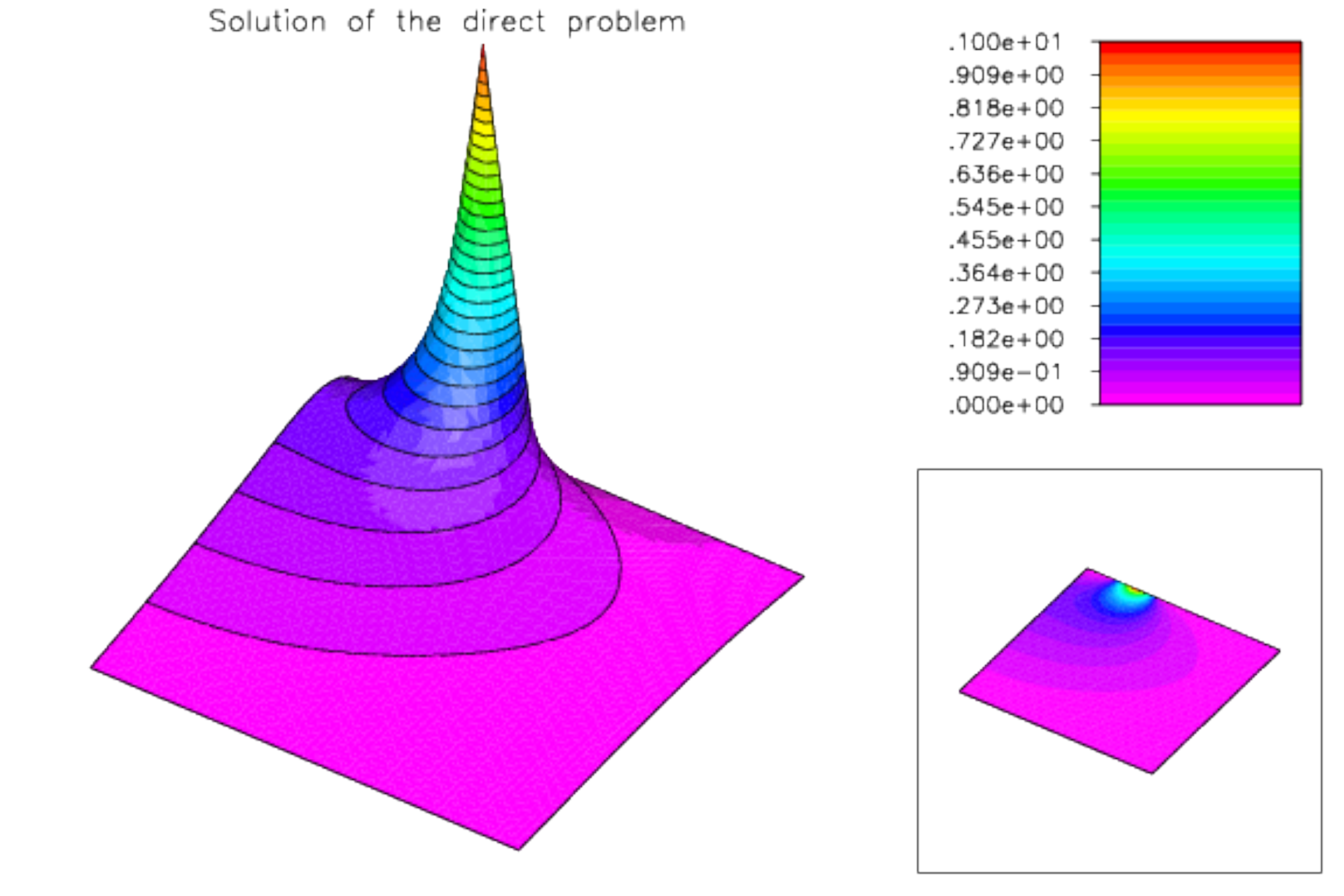} \hfill
             \epsfysize4cm \epsfbox{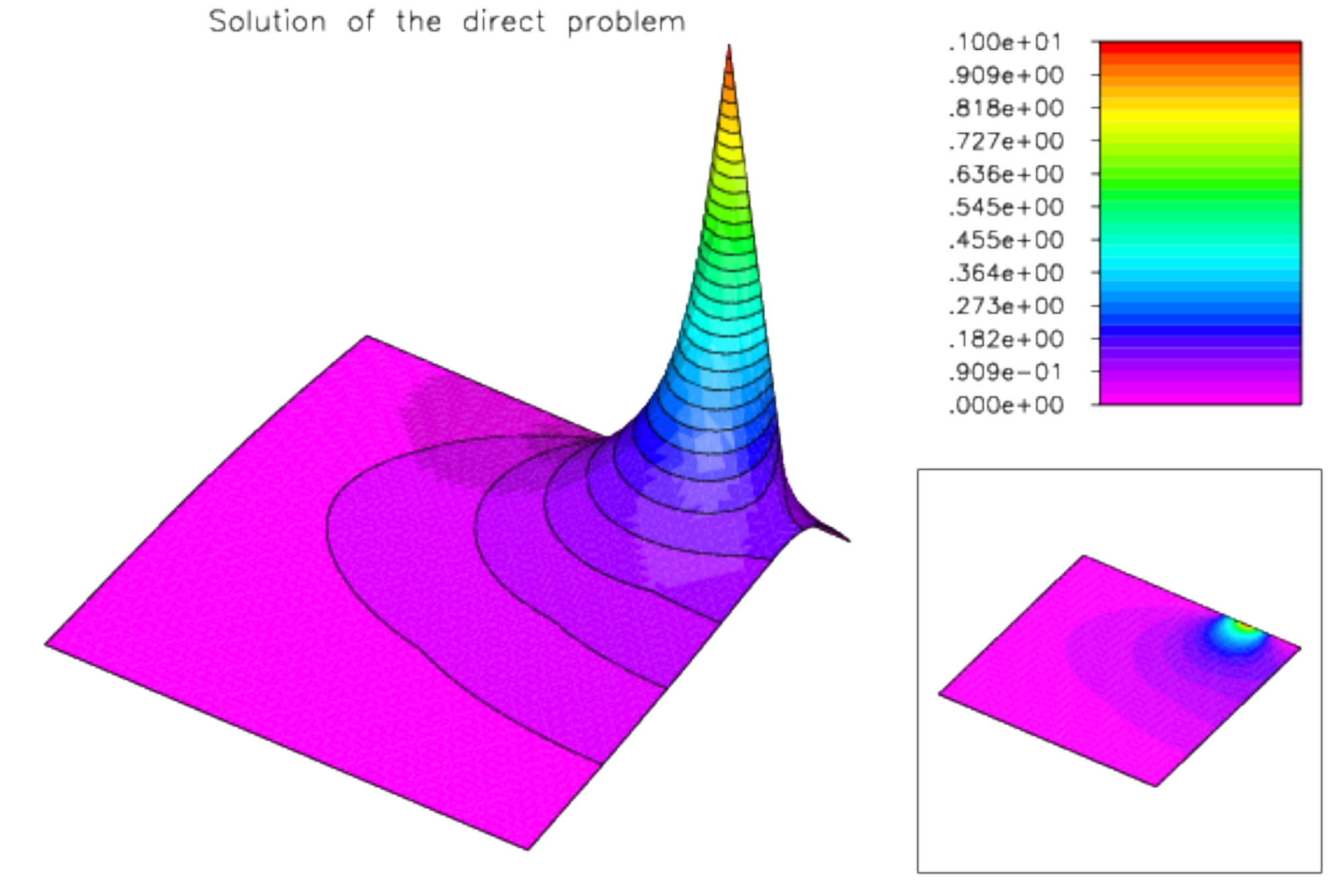}  }
\centerline{\hfil (a) \hskip5cm  (b) \hfil}
\caption{Pictures (a) and (b) show different pairs of (Dirichlet,Neumann)
data used in separate runs of the Landweber--Kaczmarz iteration.}
\label{fig:single-source}
\end{figure}

\begin{figure}
\centerline{ \epsfysize4cm \epsfbox{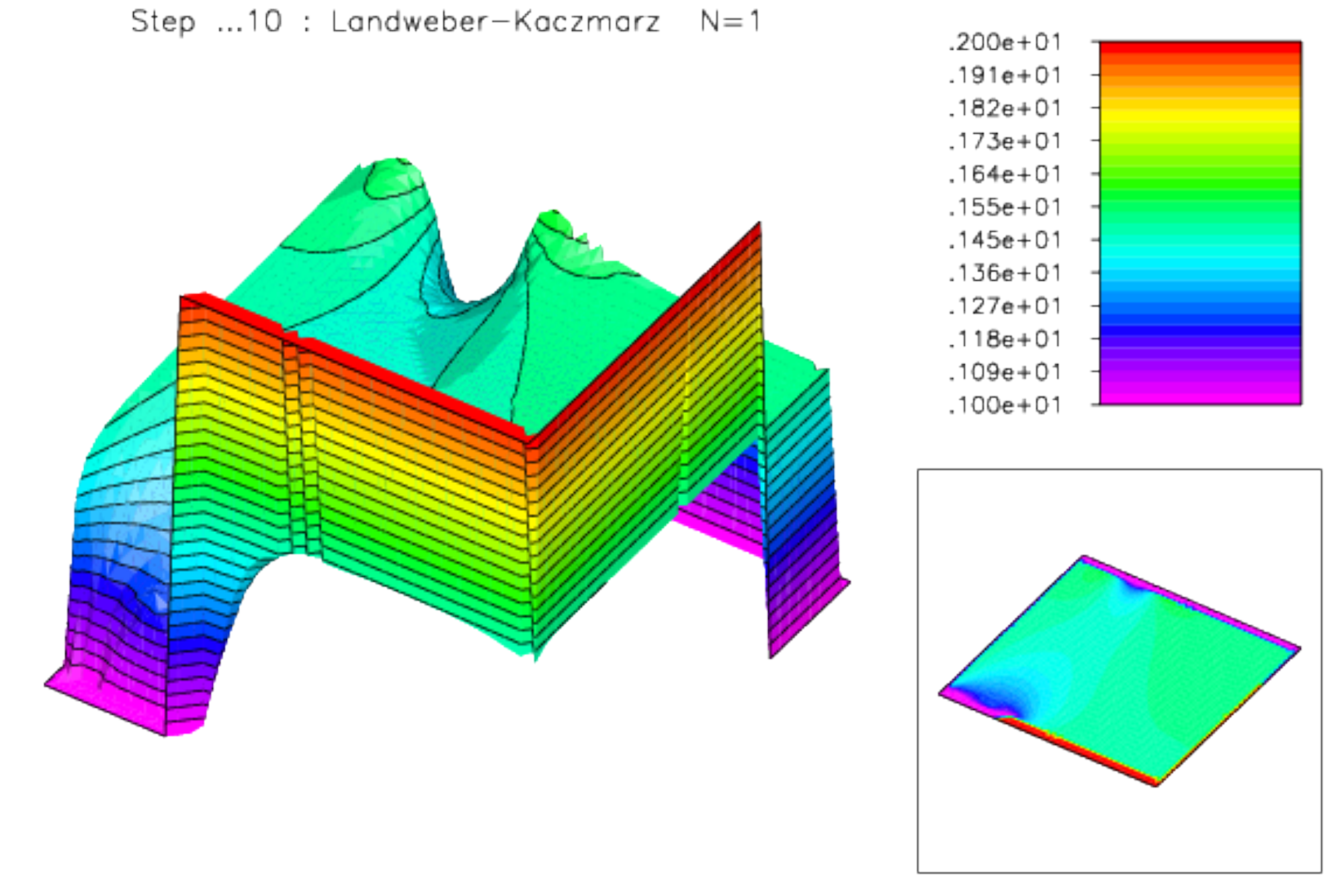} \hfill
             \epsfysize4cm \epsfbox{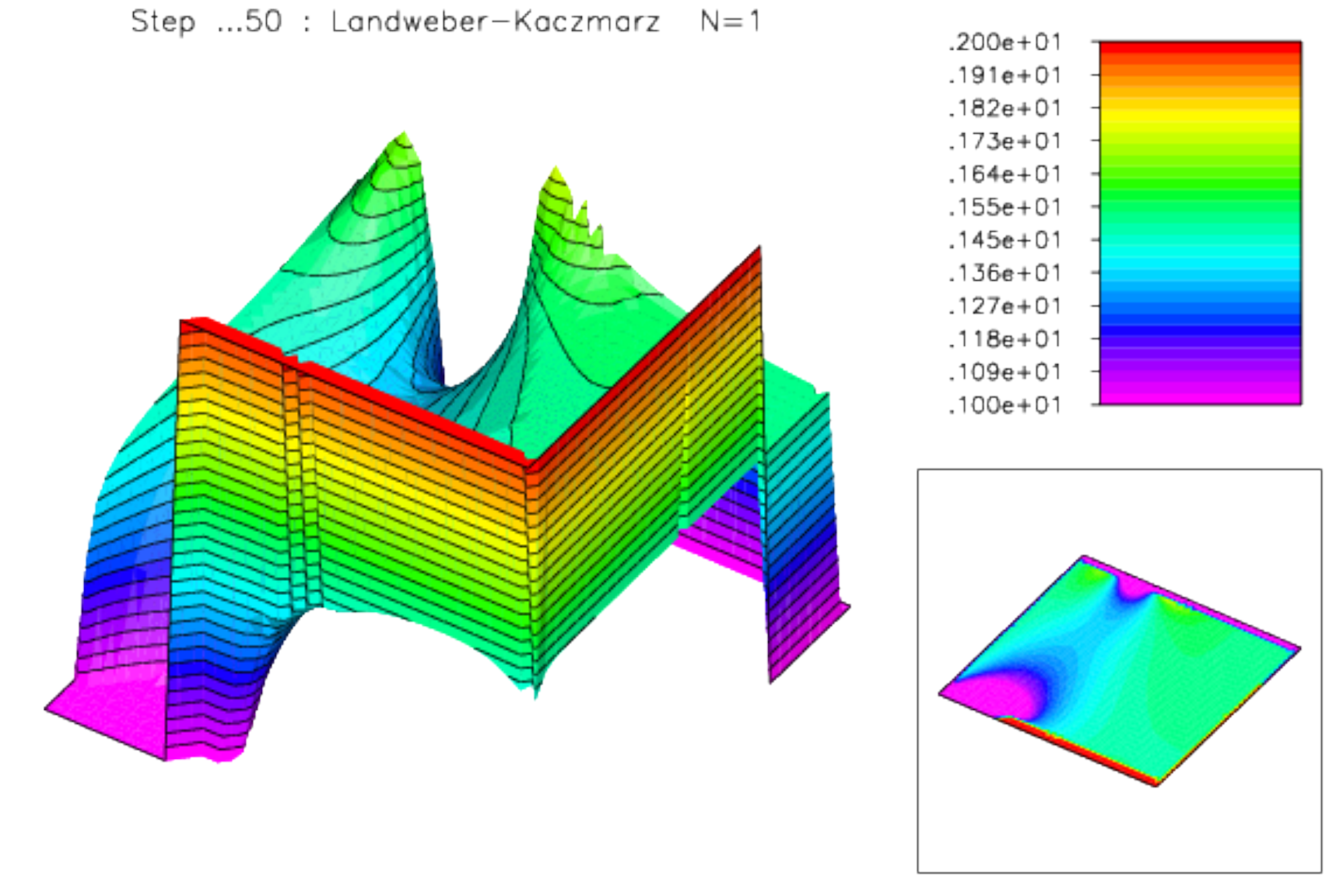}  }
\bigskip
\centerline{ \epsfysize4cm \epsfbox{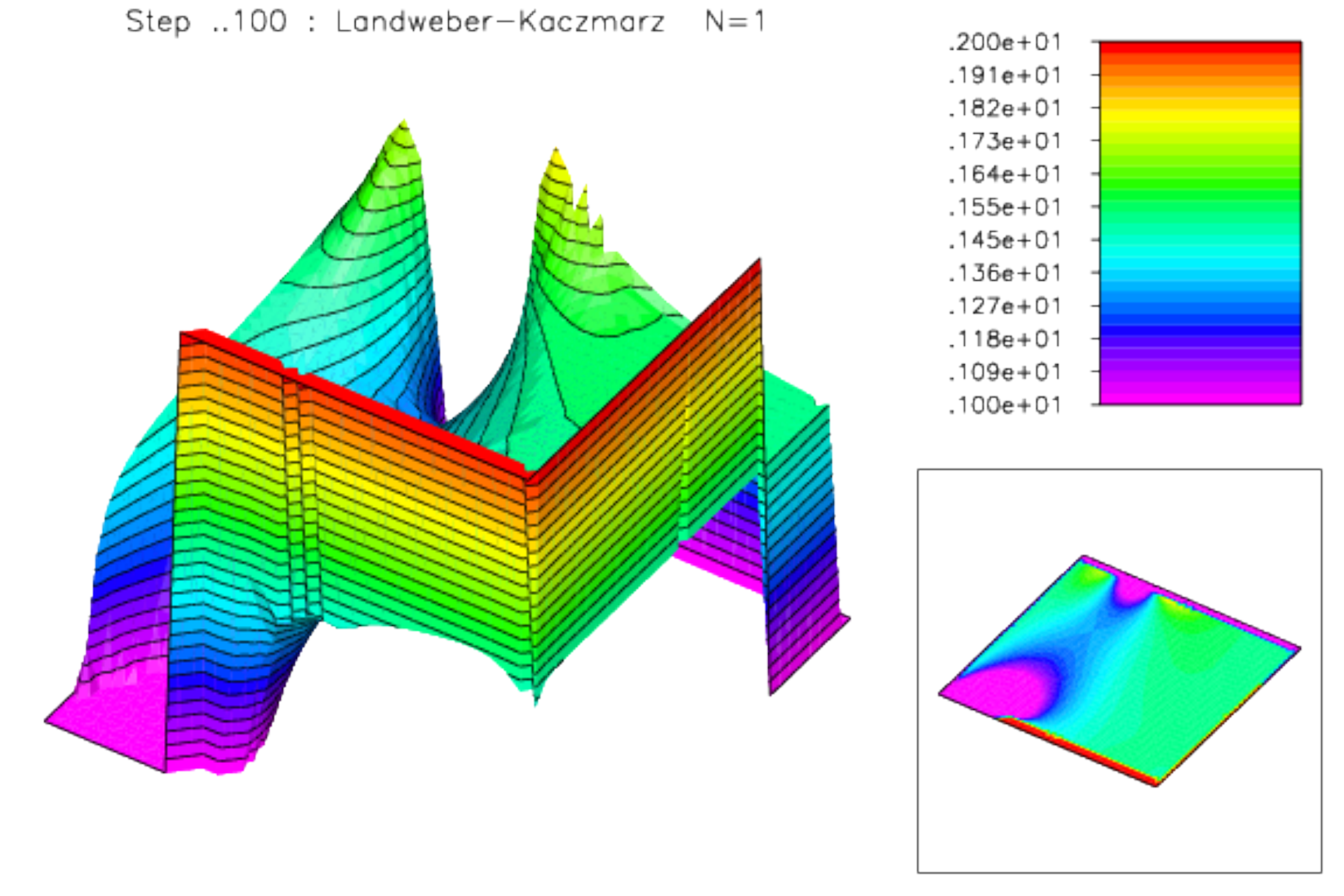} \hfill
             \epsfysize4cm \epsfbox{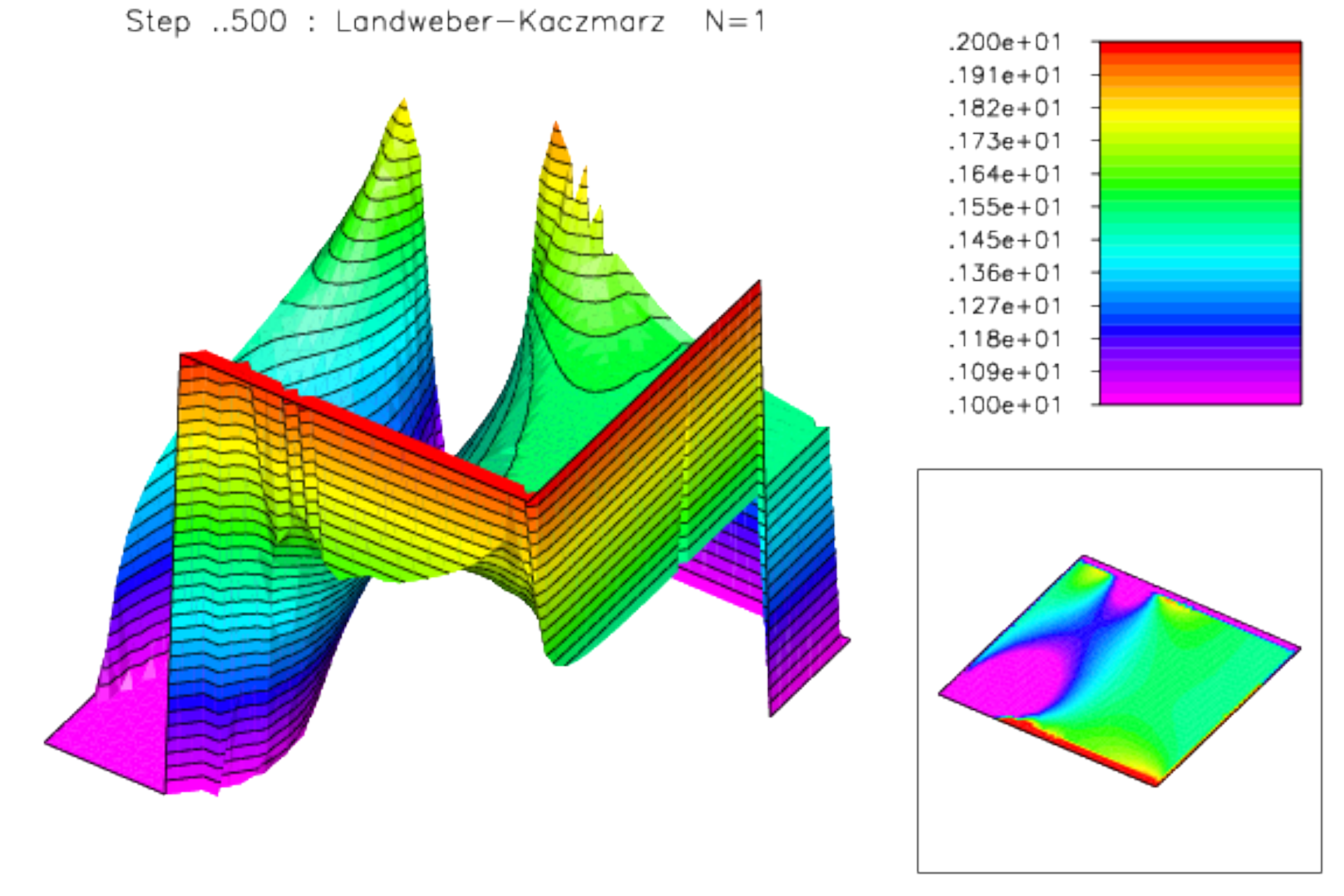}  }
\bigskip
\centerline{ \epsfysize4cm \epsfbox{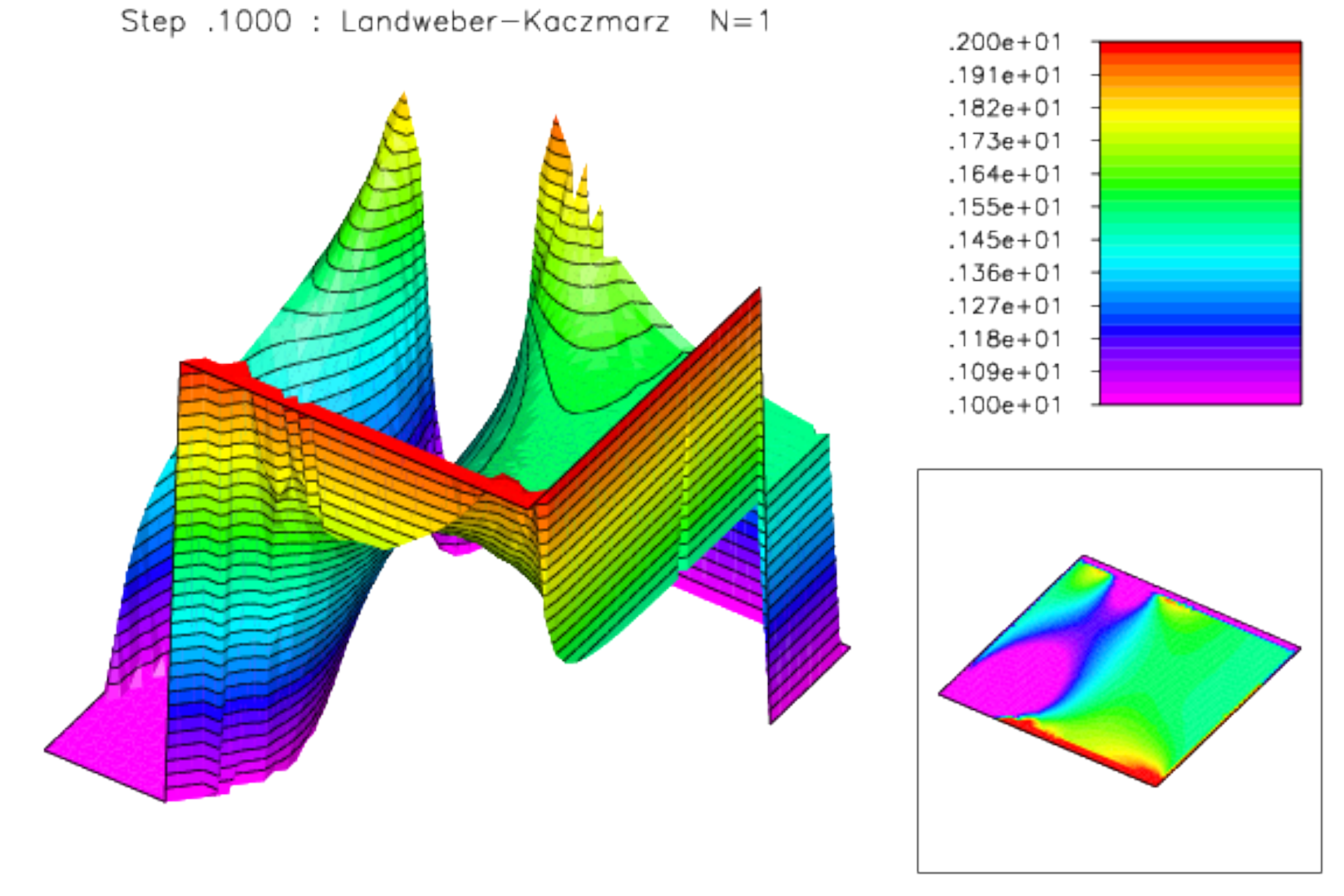} \hfill
             \epsfysize4cm \epsfbox{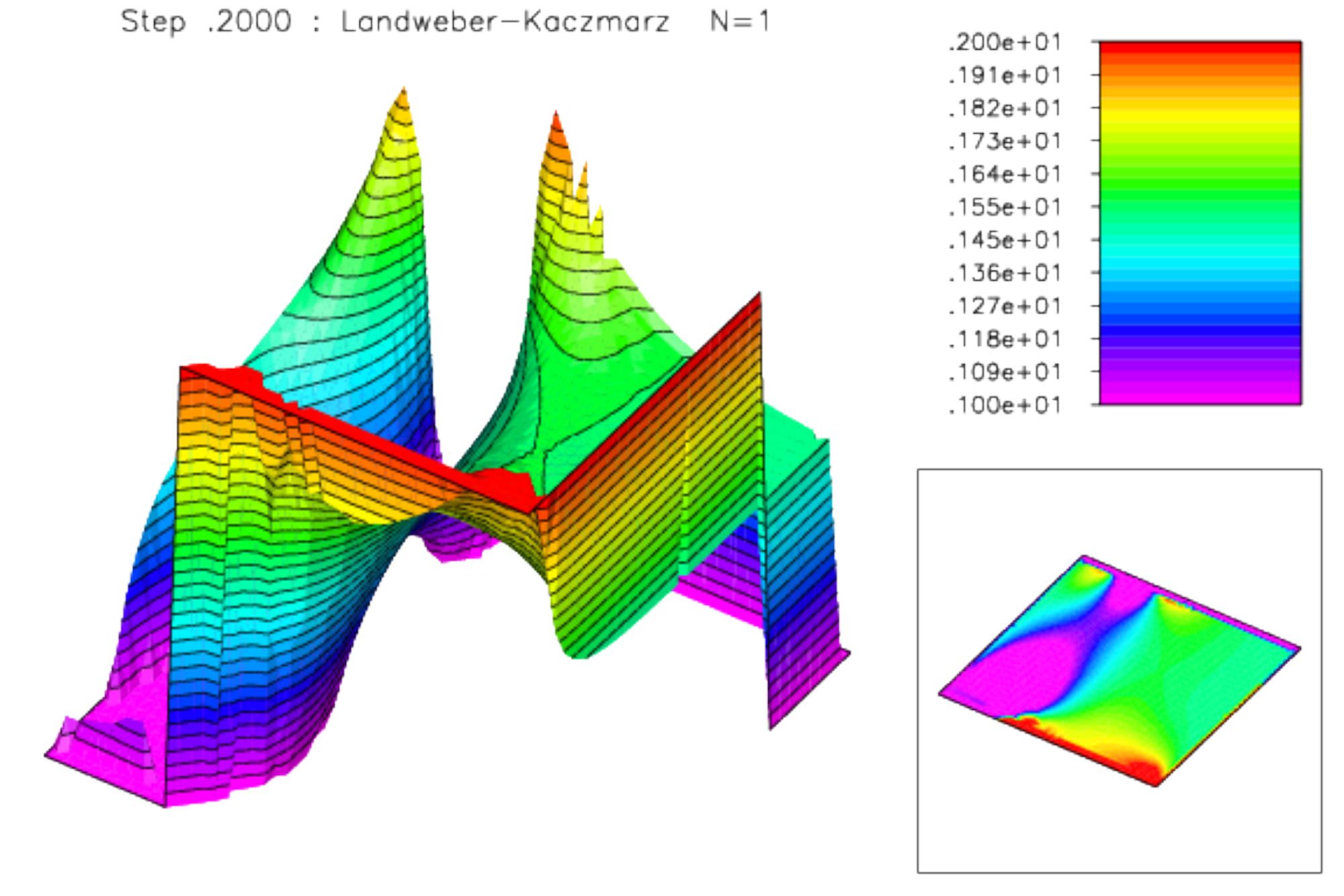}  }
\bigskip
\centerline{ \epsfysize4cm \epsfbox{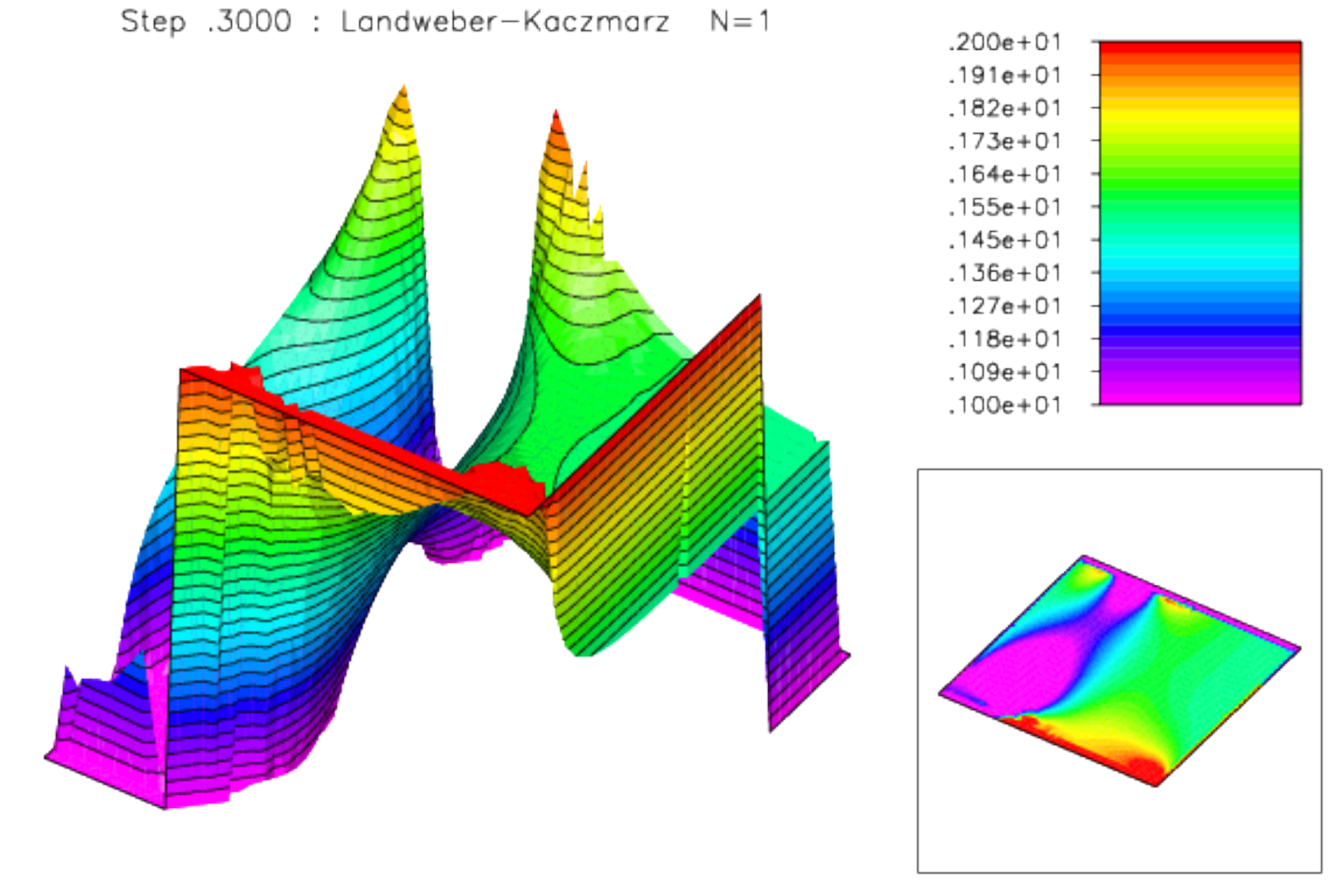} \hfill
             \epsfysize4cm \epsfbox{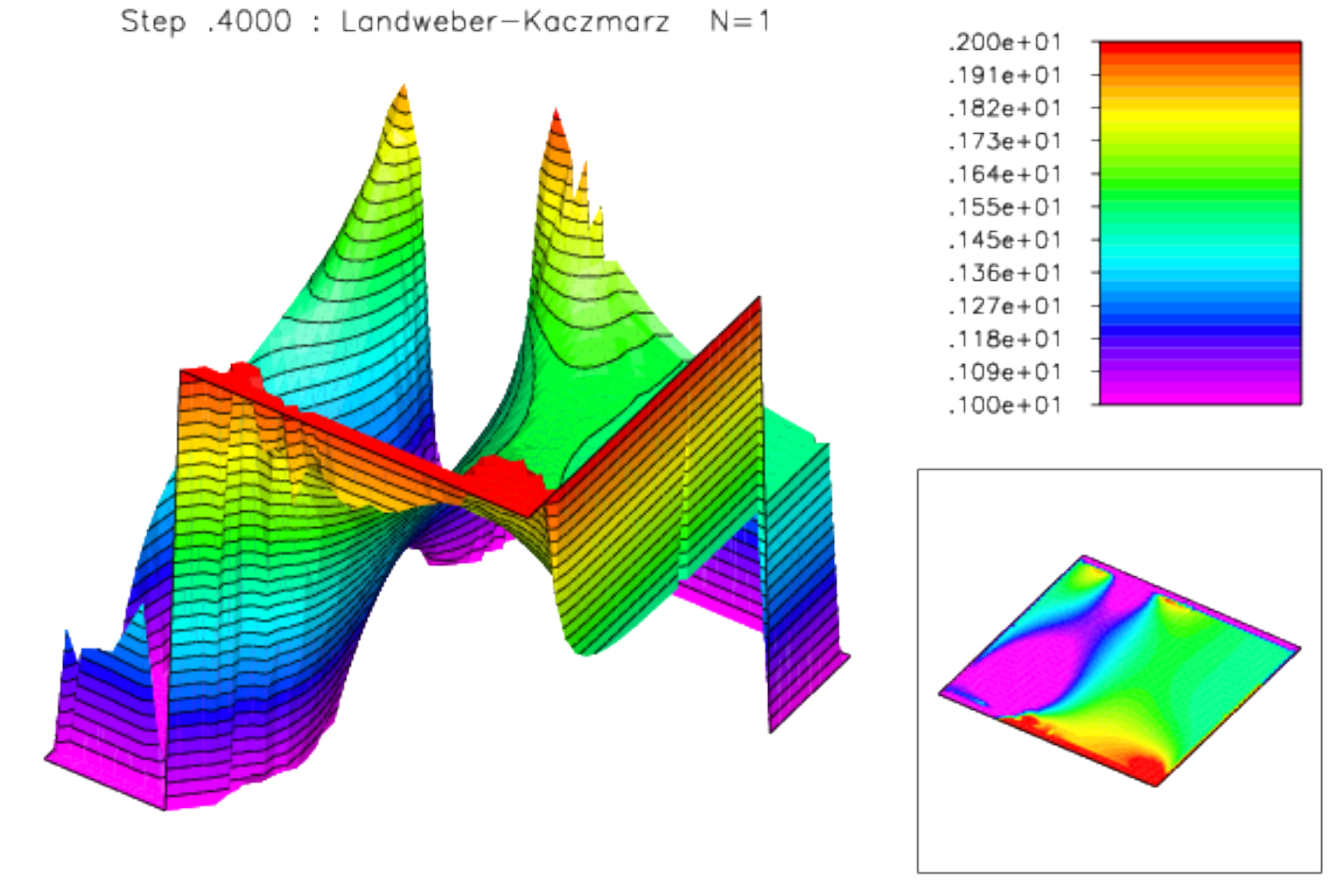}  }
\caption{Evolution of the Landweber-Kaczmarz method for one single
measurement ($N=1$). The corresponding source is shown in Figure~%
\ref{fig:single-source}~(a).} \label{fig:evol-sA}
\end{figure}

\begin{figure}
\centerline{ \epsfysize4cm \epsfbox{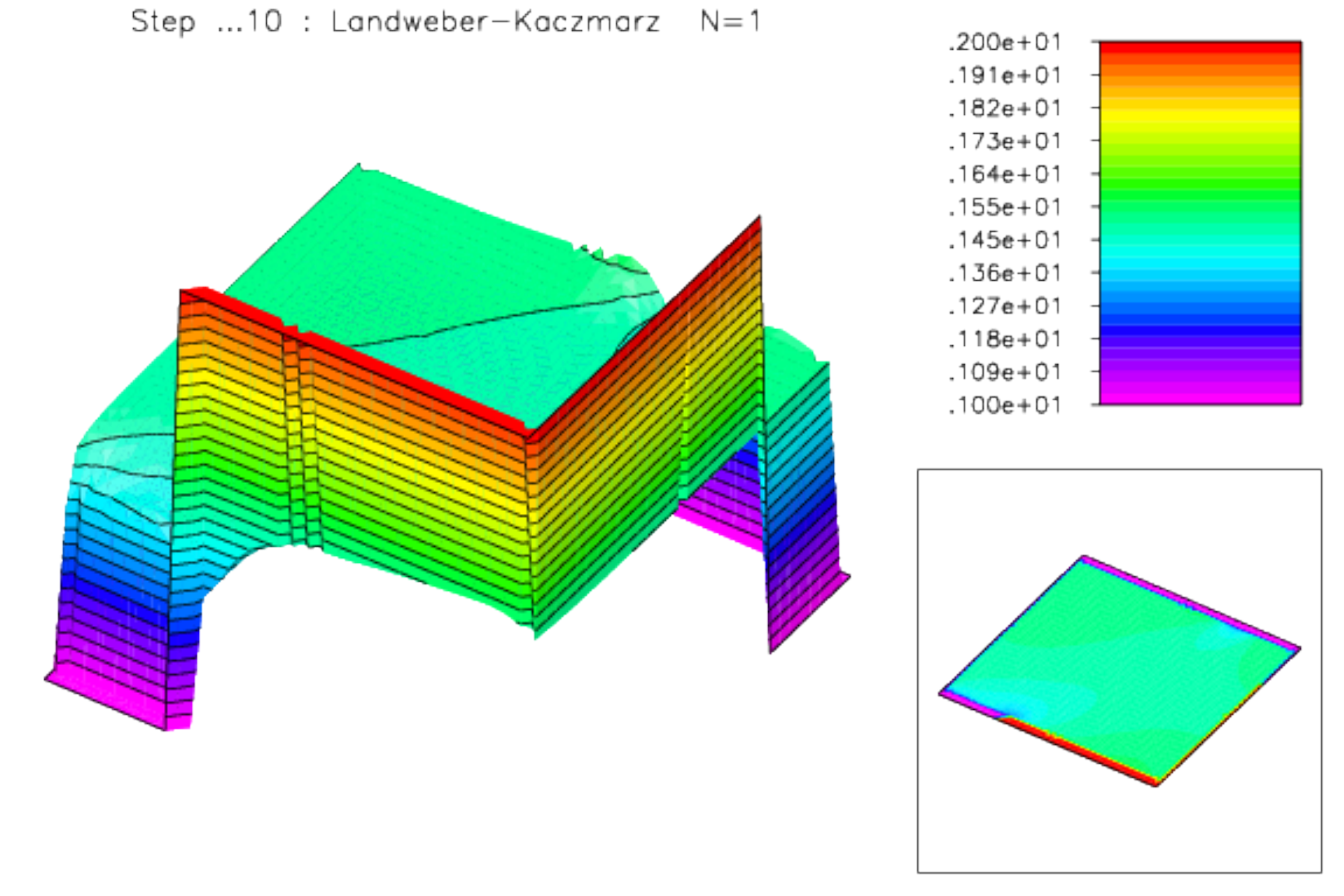} \hfill
             \epsfysize4cm \epsfbox{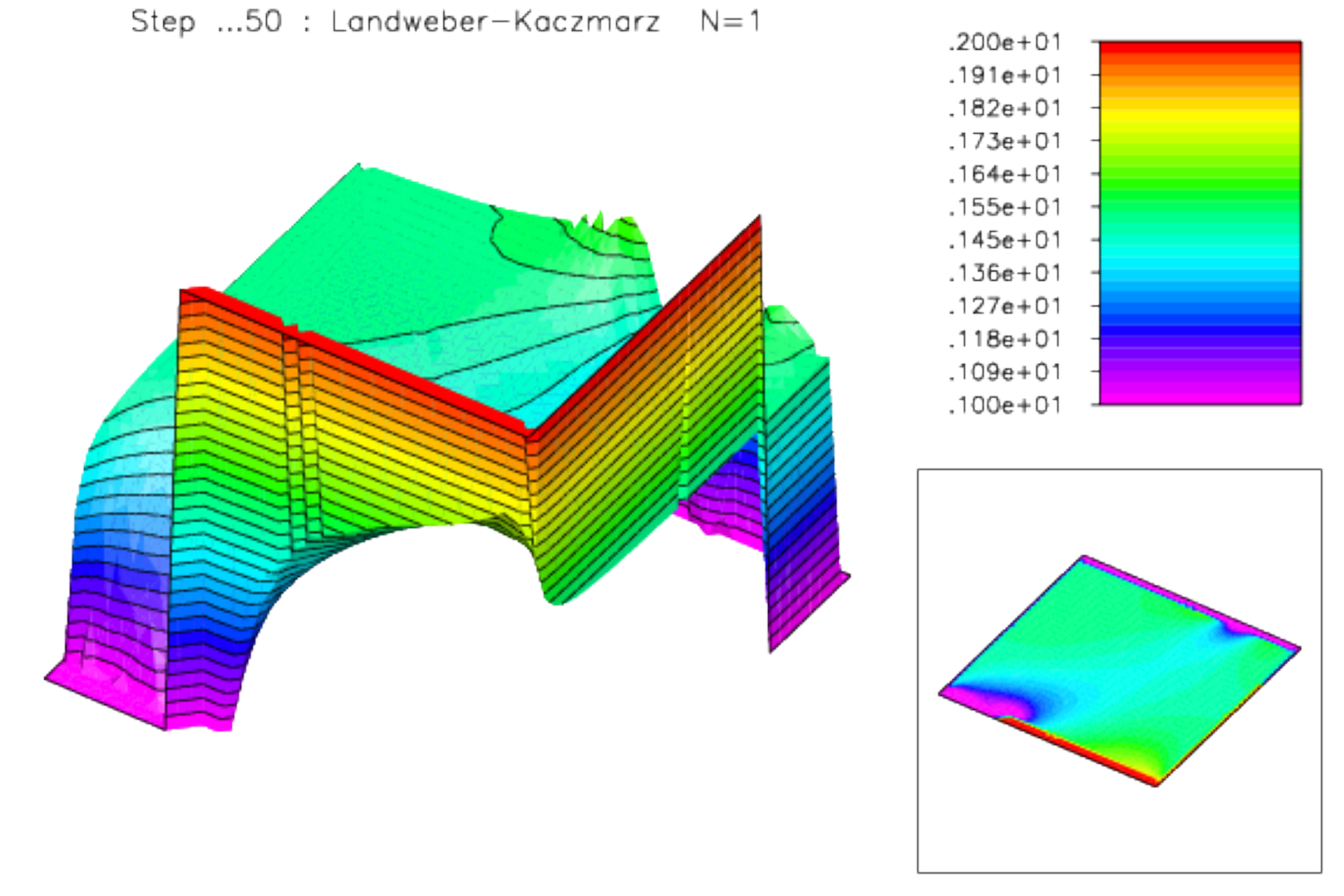}  }
\bigskip
\centerline{ \epsfysize4cm \epsfbox{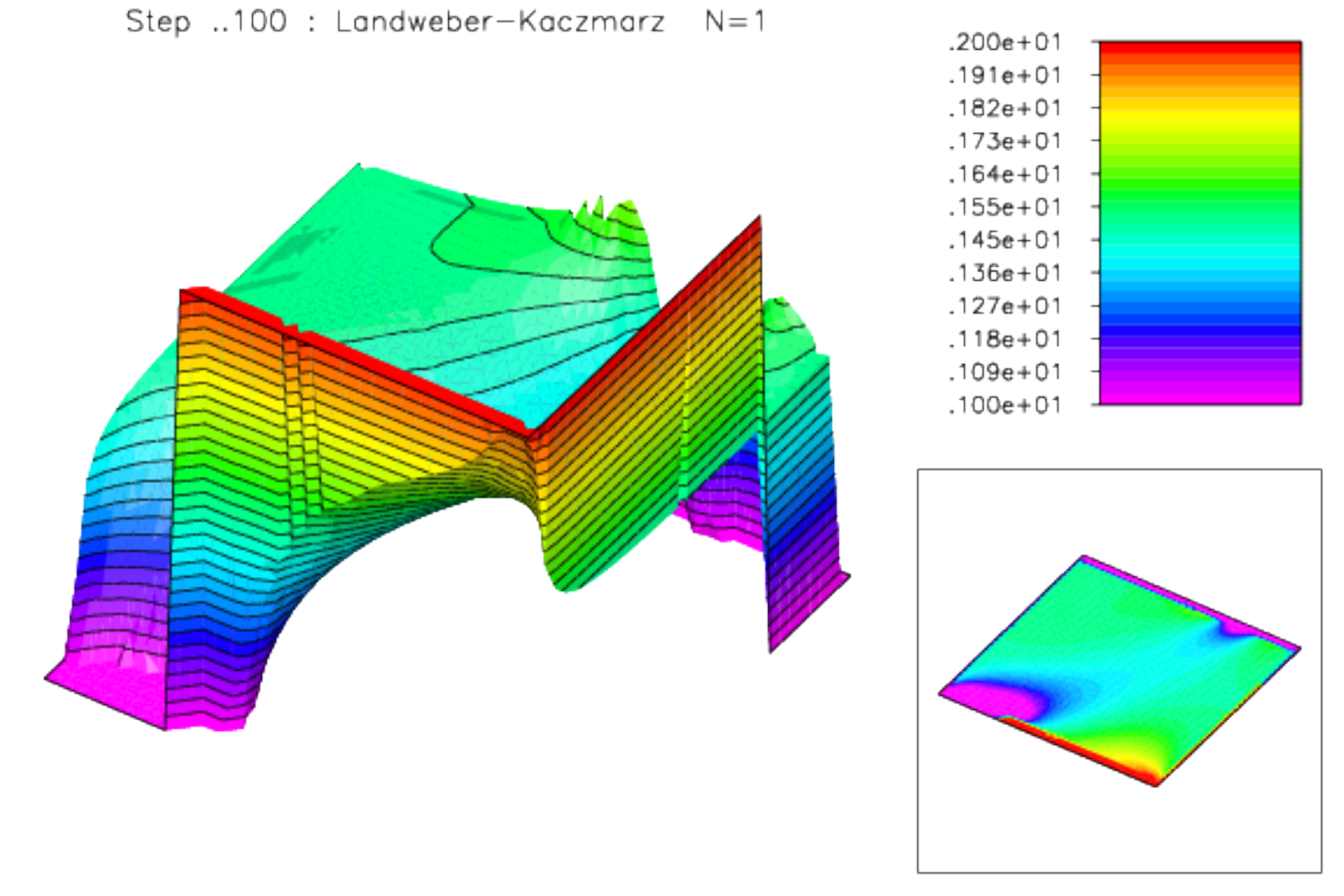} \hfill
             \epsfysize4cm \epsfbox{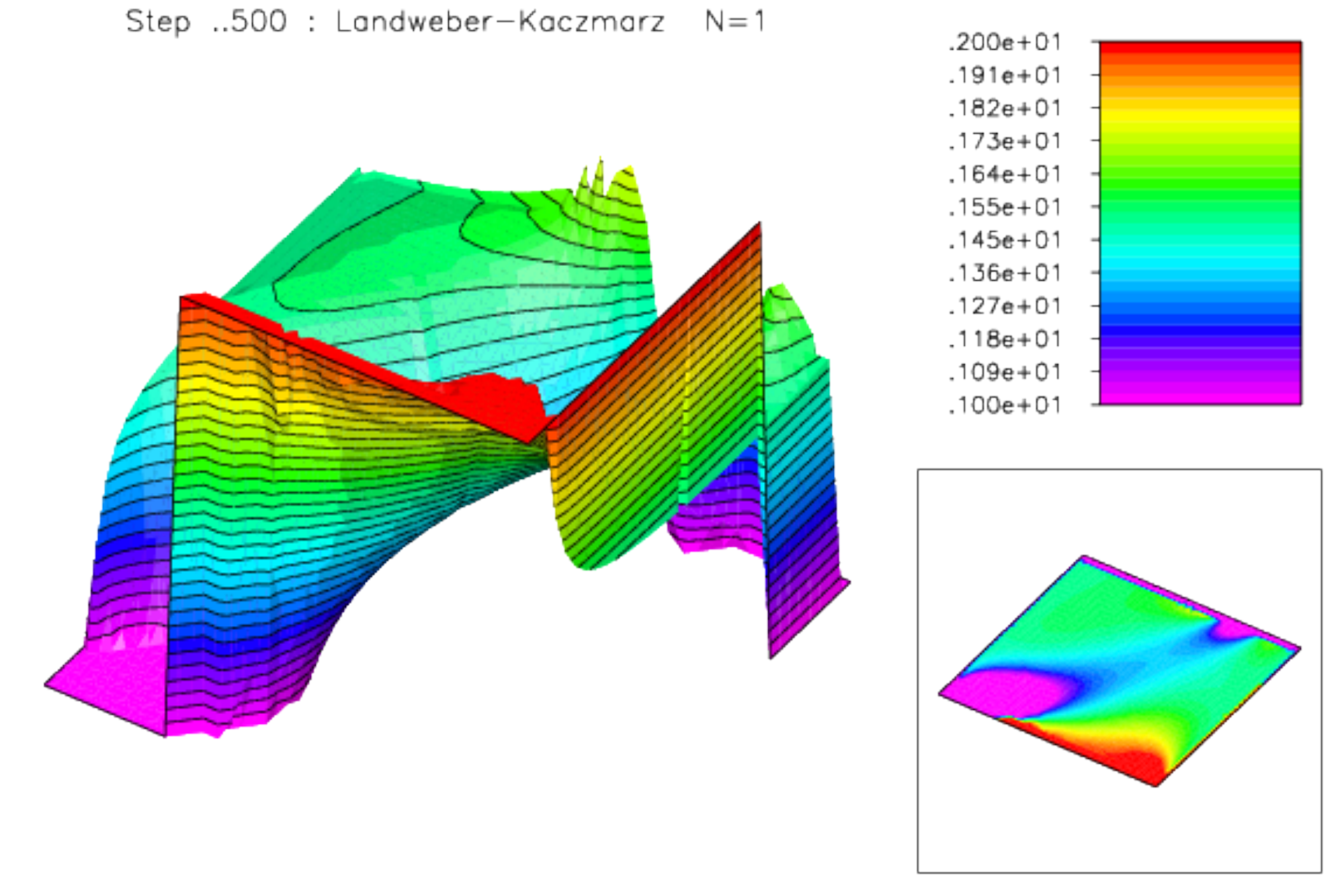}  }
\bigskip
\centerline{ \epsfysize4cm \epsfbox{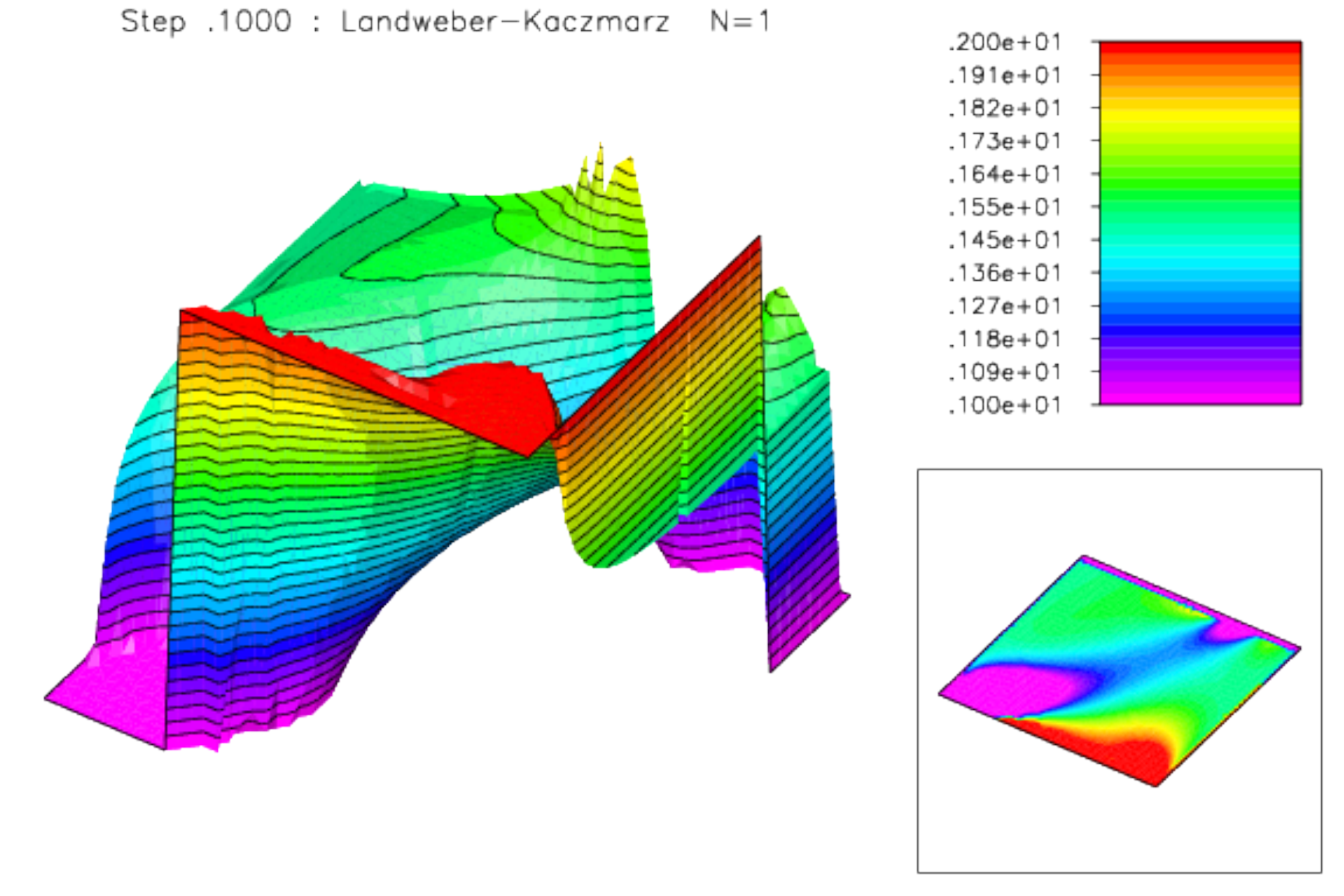} \hfill
             \epsfysize4cm \epsfbox{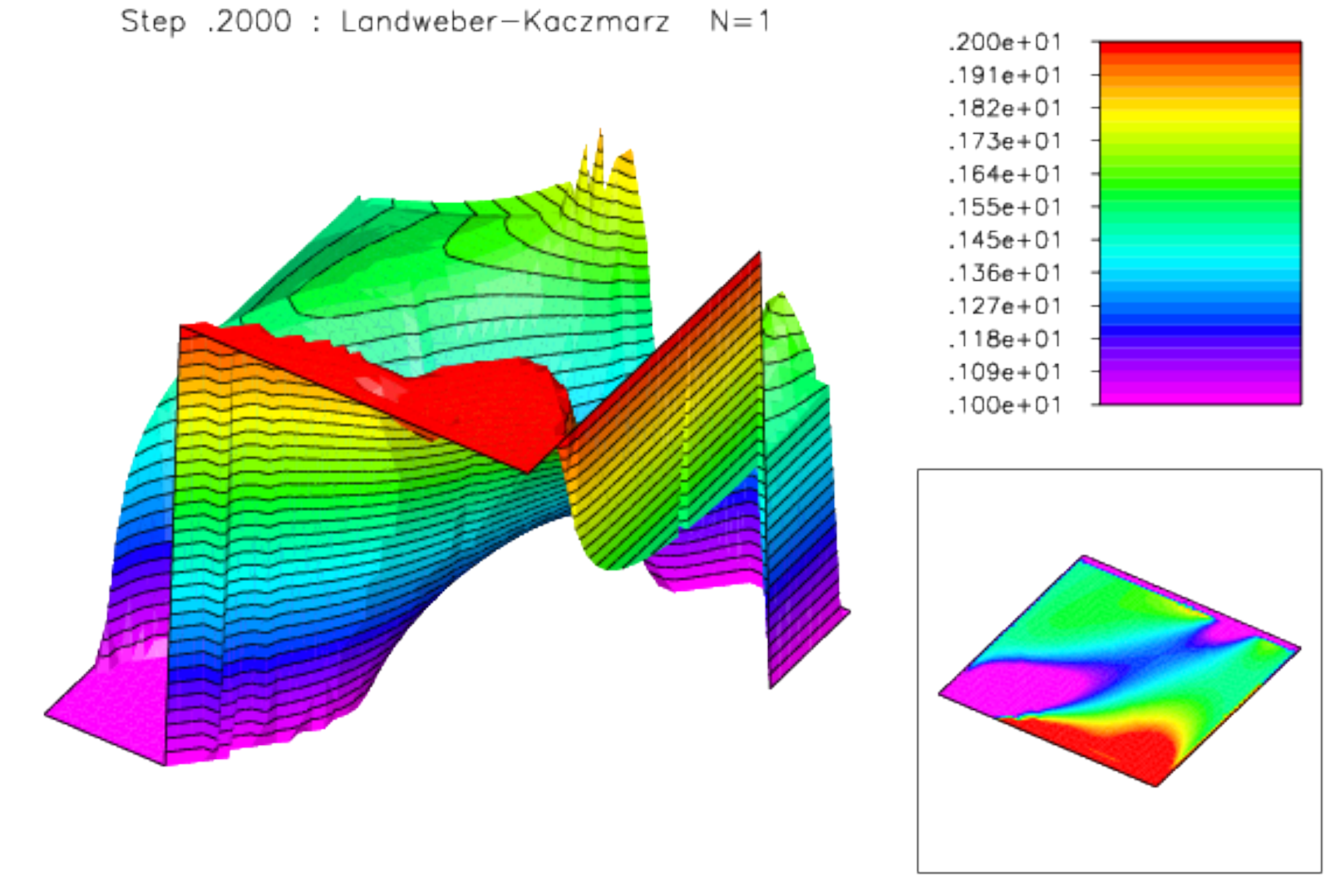}  }
\bigskip
\centerline{ \epsfysize4cm \epsfbox{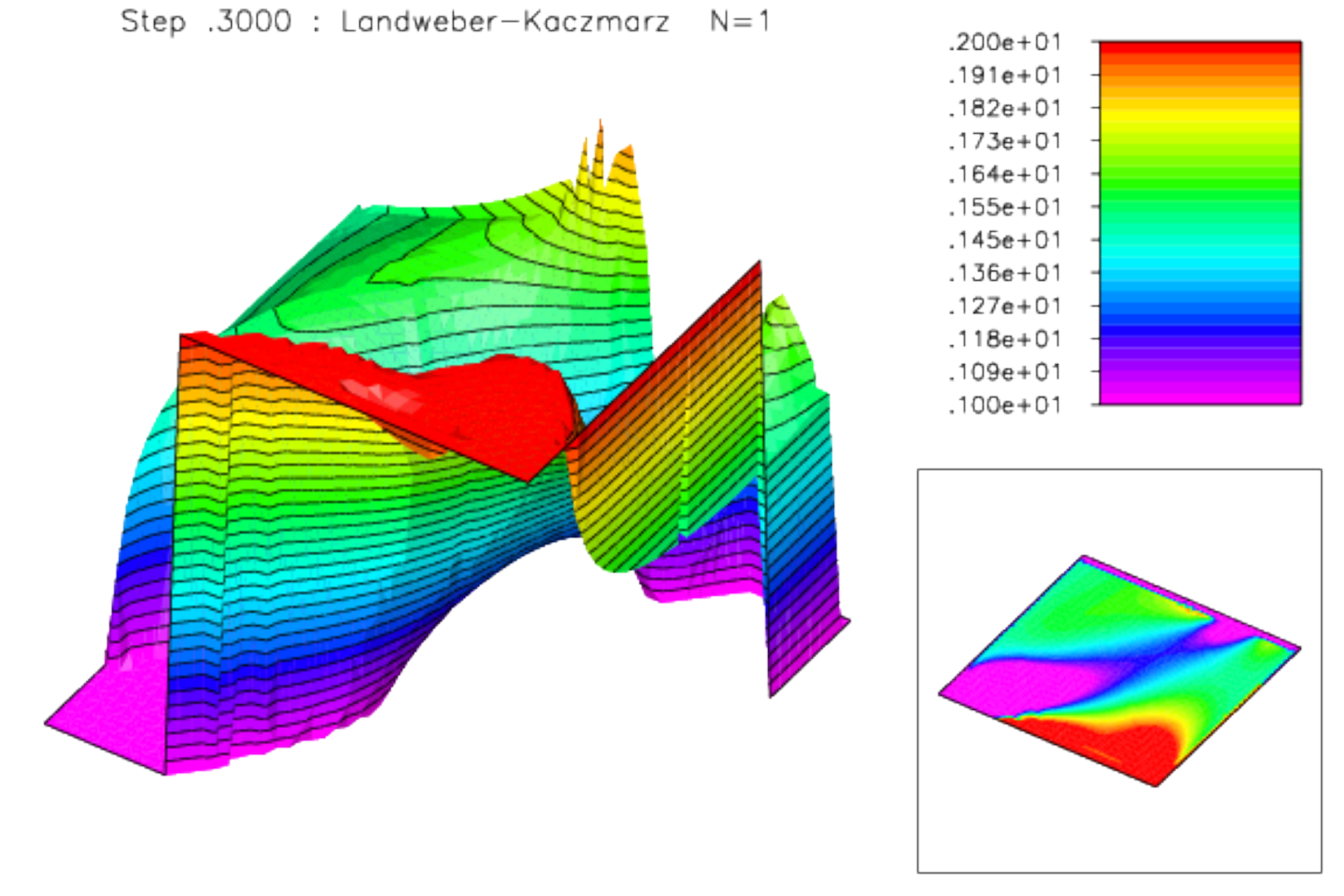} \hfill
             \epsfysize4cm \epsfbox{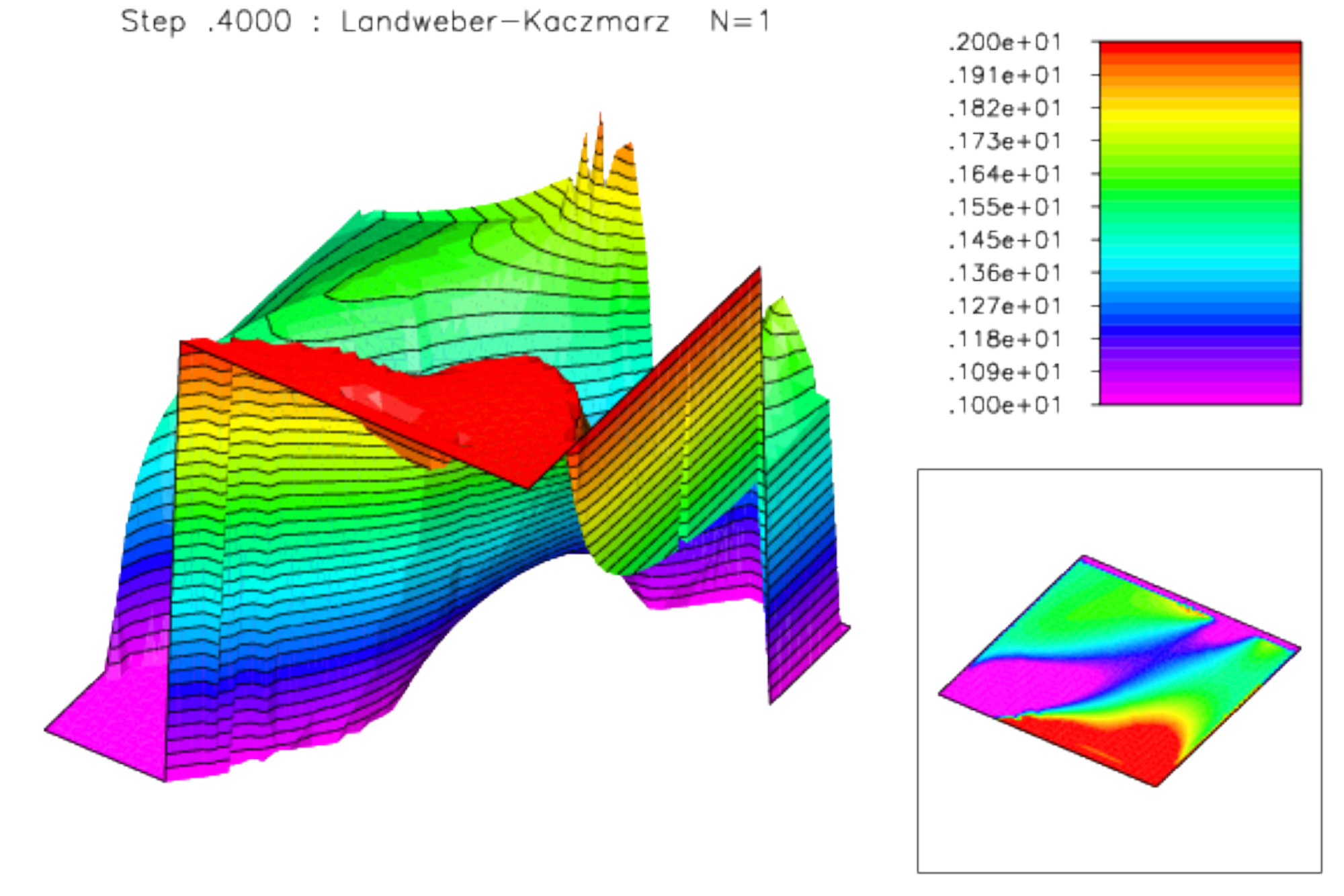}  }
\caption{Evolution of the Landweber-Kaczmarz method for one single
measurement ($N=1$). The corresponding source is shown in Figure~%
\ref{fig:single-source}~(b).} \label{fig:evol-sB}
\end{figure}

Concerning the amount of information used in the reconstruction,
we implemented (for comparison purposes) the Landweber--Kaczmarz
iteration in the case where a single pair of Dirichlet and Neumann
data was available. In this case, the Landweber--Kaczmarz method
reduces to the Landweber iterative method. This experiment is
interesting, since it shows that the quality of the reconstruction
is better at the part of the domain $\Omega$ which is closer to the
support of the applied voltage.
In Figure~\ref{fig:evol-sA} we present the evolution of the Landweber
iteration for $N=1$ and
$$  U_1(x) \ := \ \left\{ \begin{array}{rl}
      1 - 8 |x - \frac{6}{8}| , & |x-\frac{6}{8}| \le \frac{1}{8} \\
      0, & {\rm else} \end{array} \right. $$
The solution of the direct problem corresponding to this choice of $U_1$
is shown in Figure~\ref{fig:single-source}~(a).
In Figure~\ref{fig:evol-sB} we present the evolution of the Landweber
iteration for $N=1$ and
$$  U_1(x) \ := \ \left\{ \begin{array}{rl}
      1 - 8 |x - \frac{2}{8}| , & |x-\frac{2}{8}| \le \frac{1}{8} \\
      0, & {\rm else} \end{array} \right. $$
The solution of the direct problem corresponding to this choice of $U_1$
is shown in Figure~\ref{fig:single-source}~(b).

\begin{figure}
\centerline{ \epsfysize4cm \epsfbox{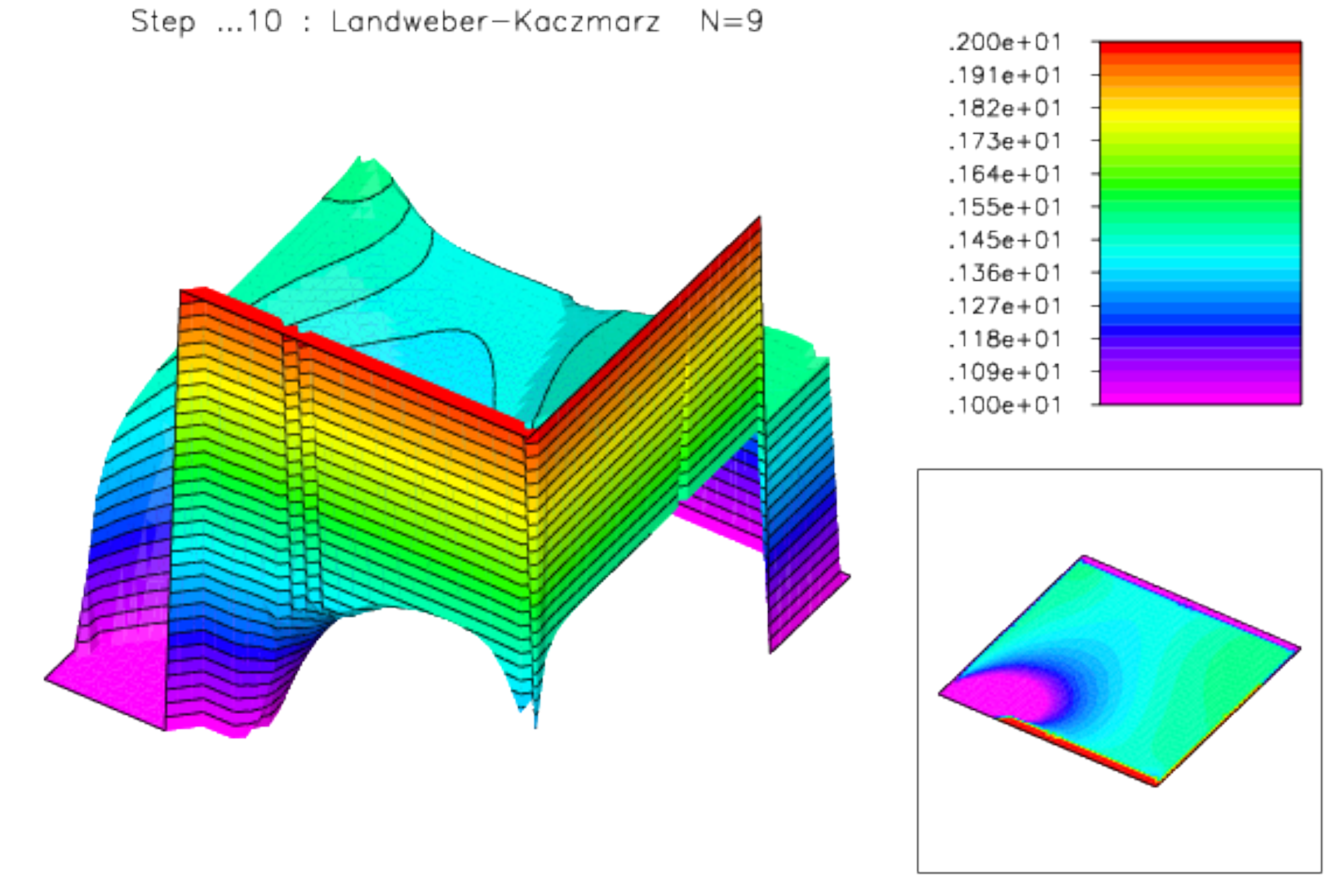} \hfill
             \epsfysize4cm \epsfbox{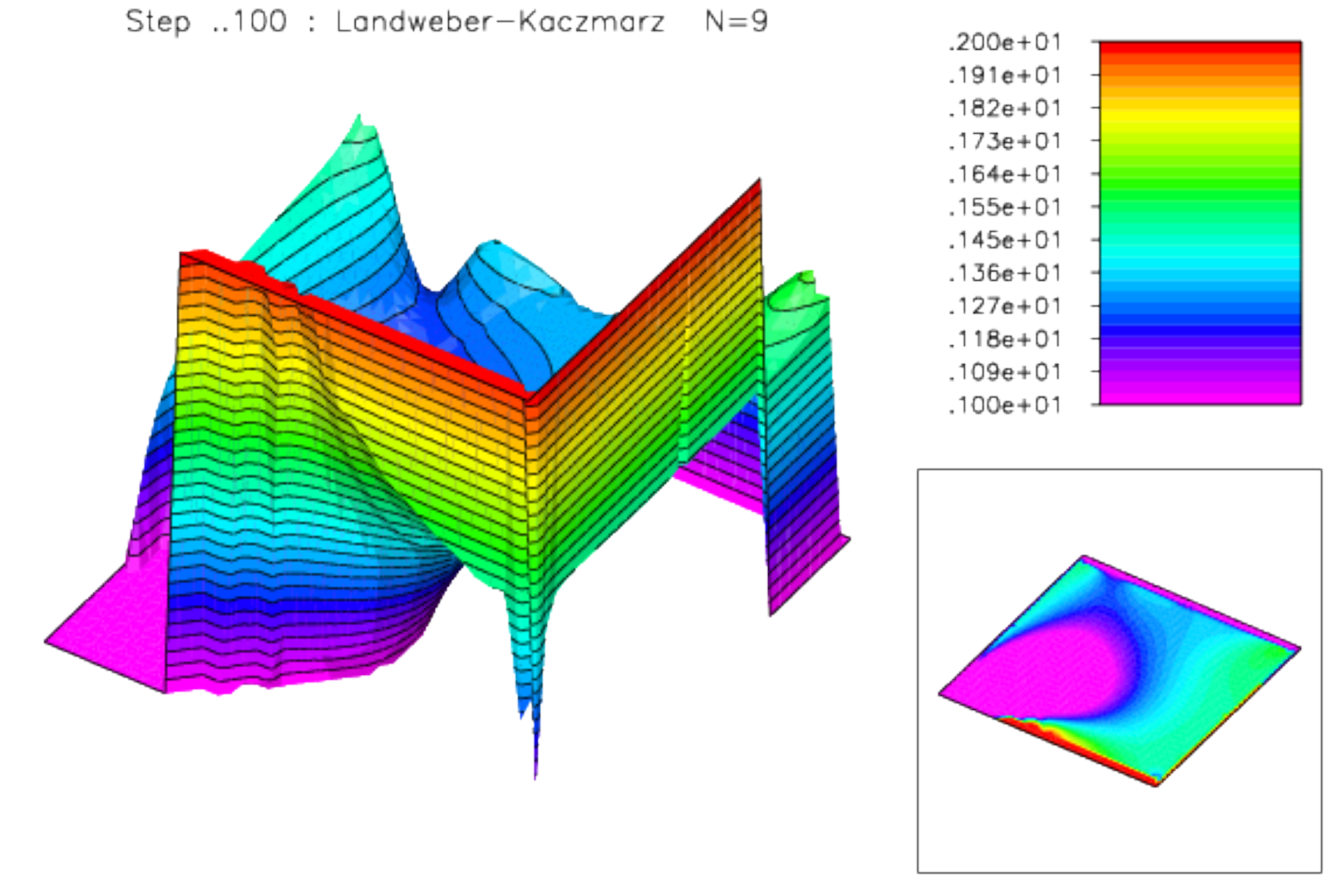}  }
\bigskip
\centerline{ \epsfysize4cm \epsfbox{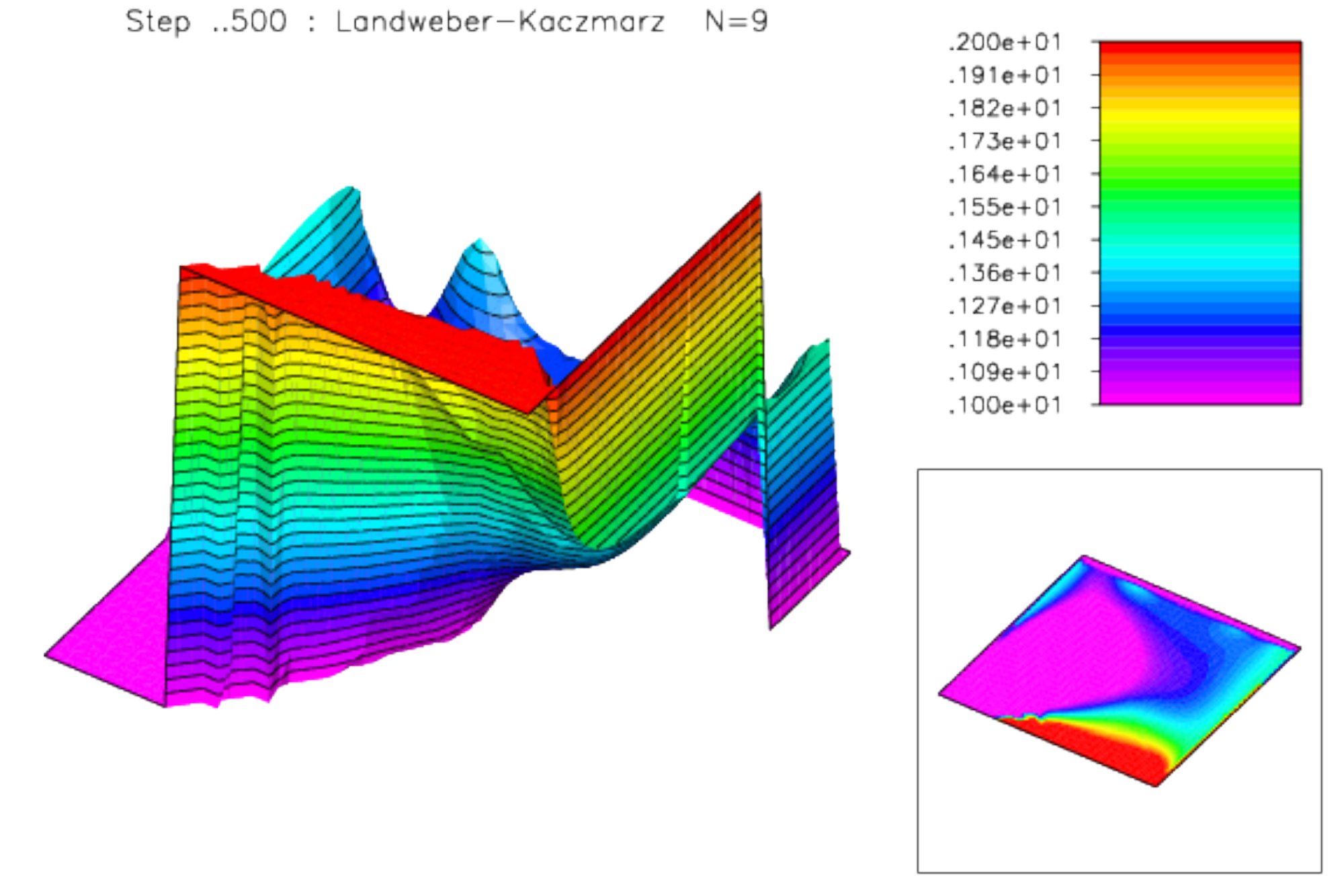} \hfill
             \epsfysize4cm \epsfbox{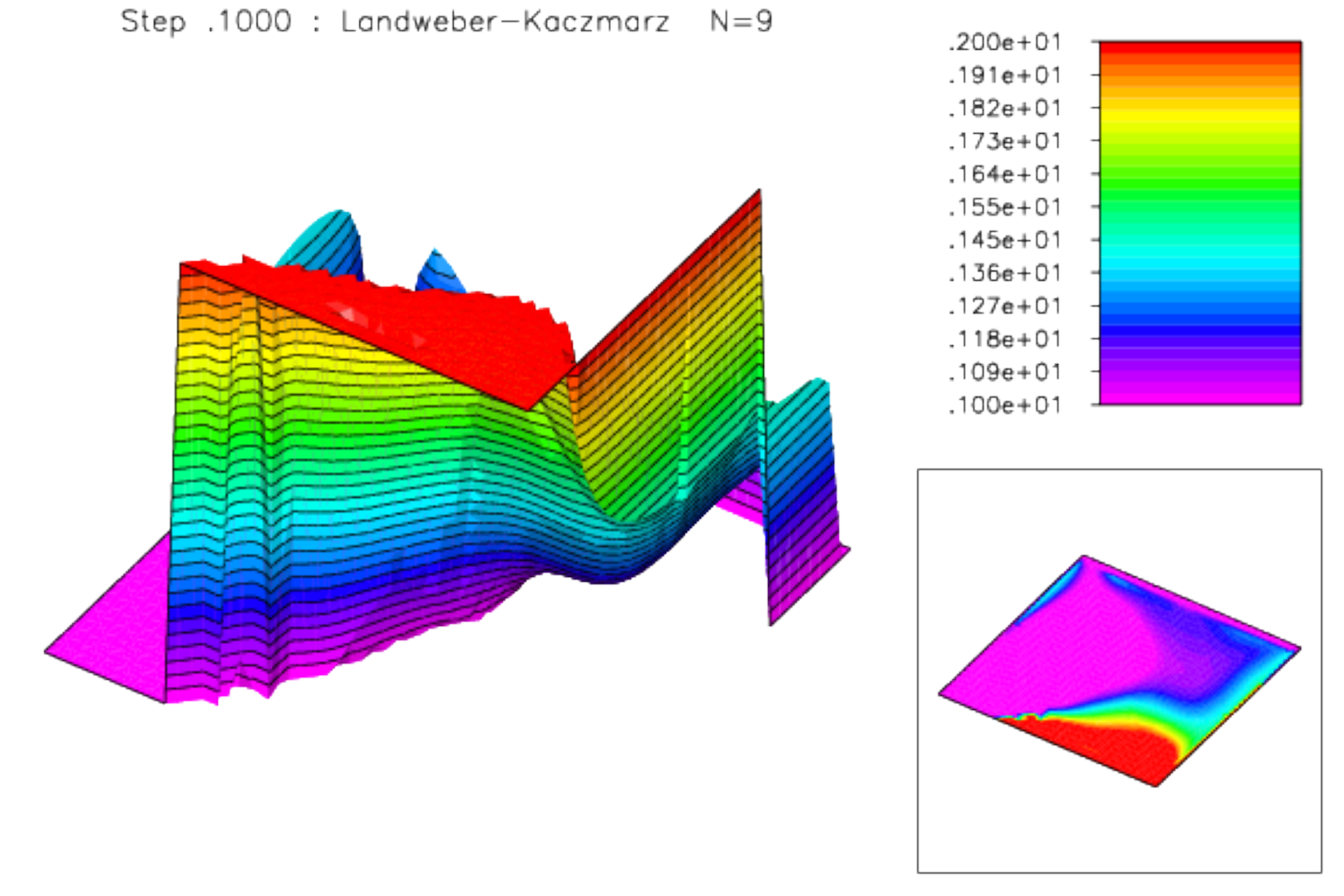}  }
\bigskip
\centerline{ \epsfysize4cm \epsfbox{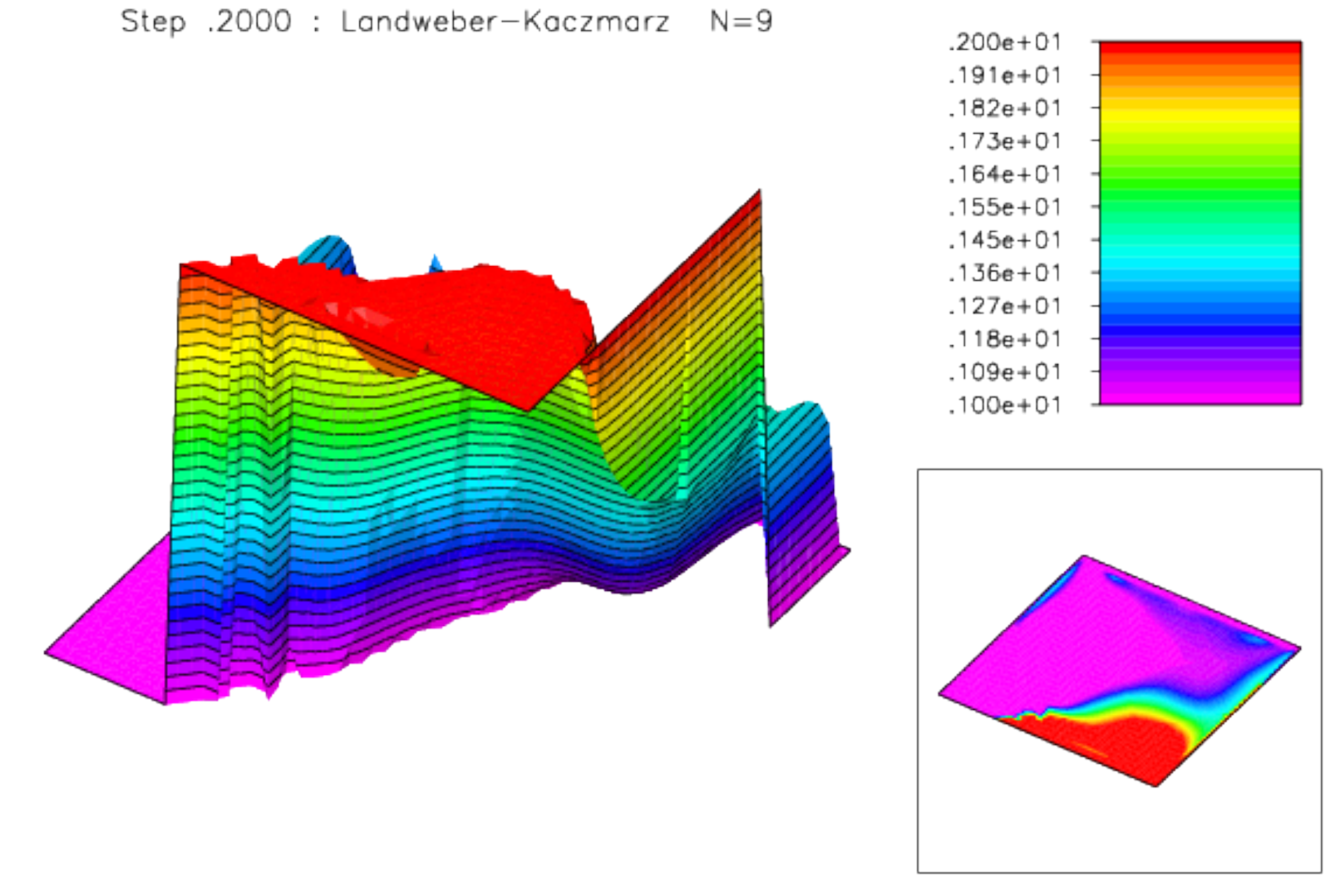} \hfill
             \epsfysize4cm \epsfbox{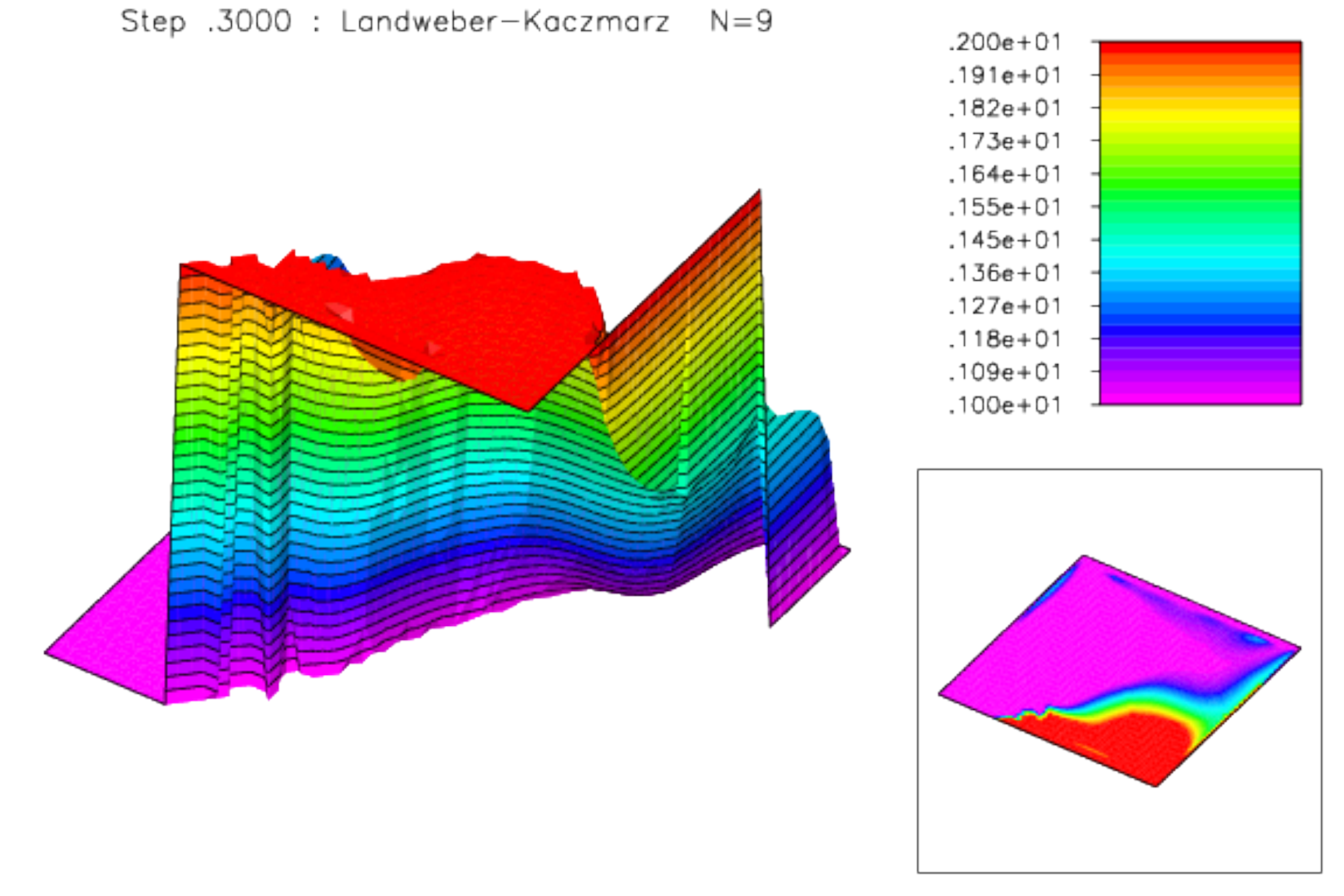}  }
\bigskip
\centerline{ \epsfysize4cm \epsfbox{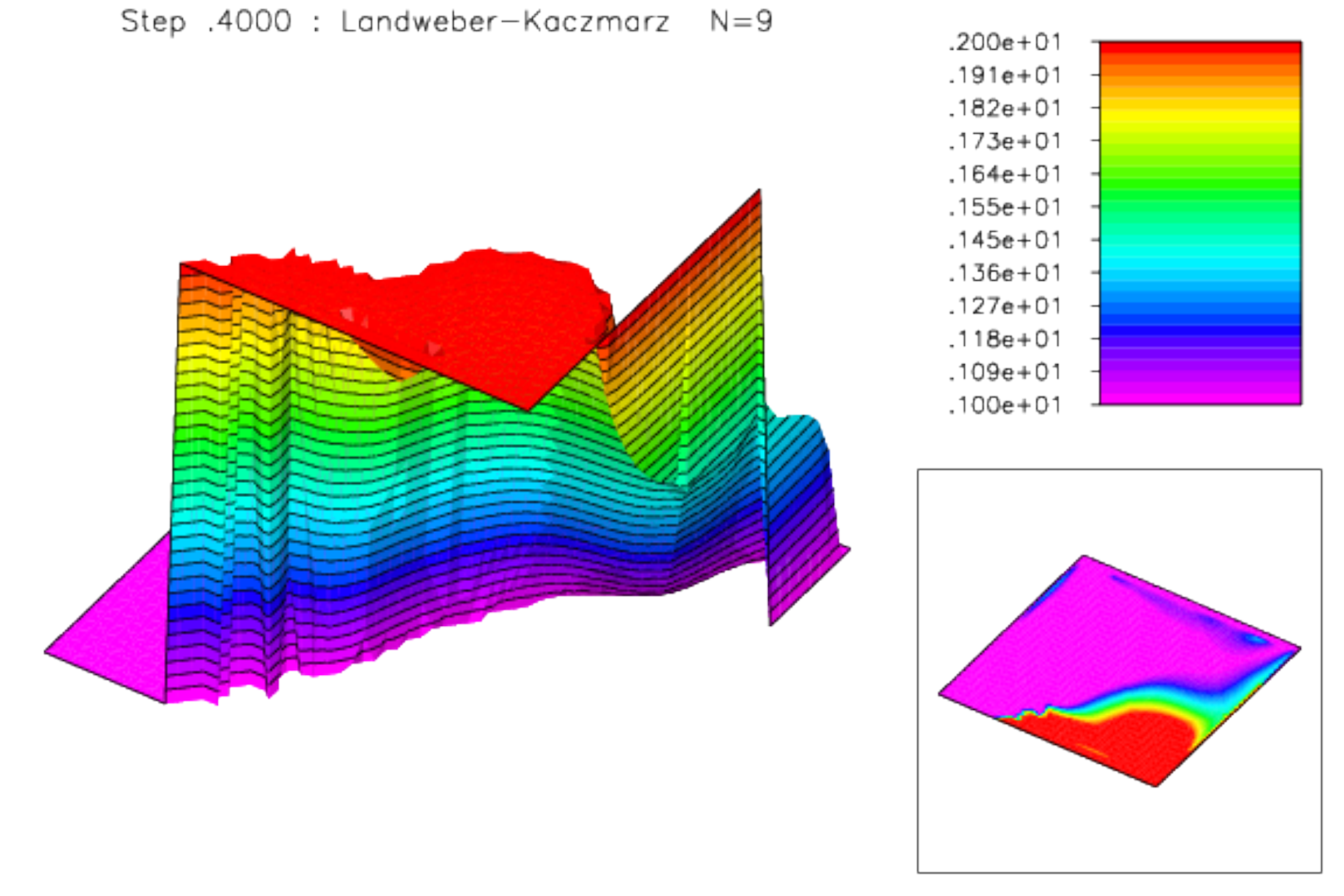} \hfill
             \epsfysize4cm \epsfbox{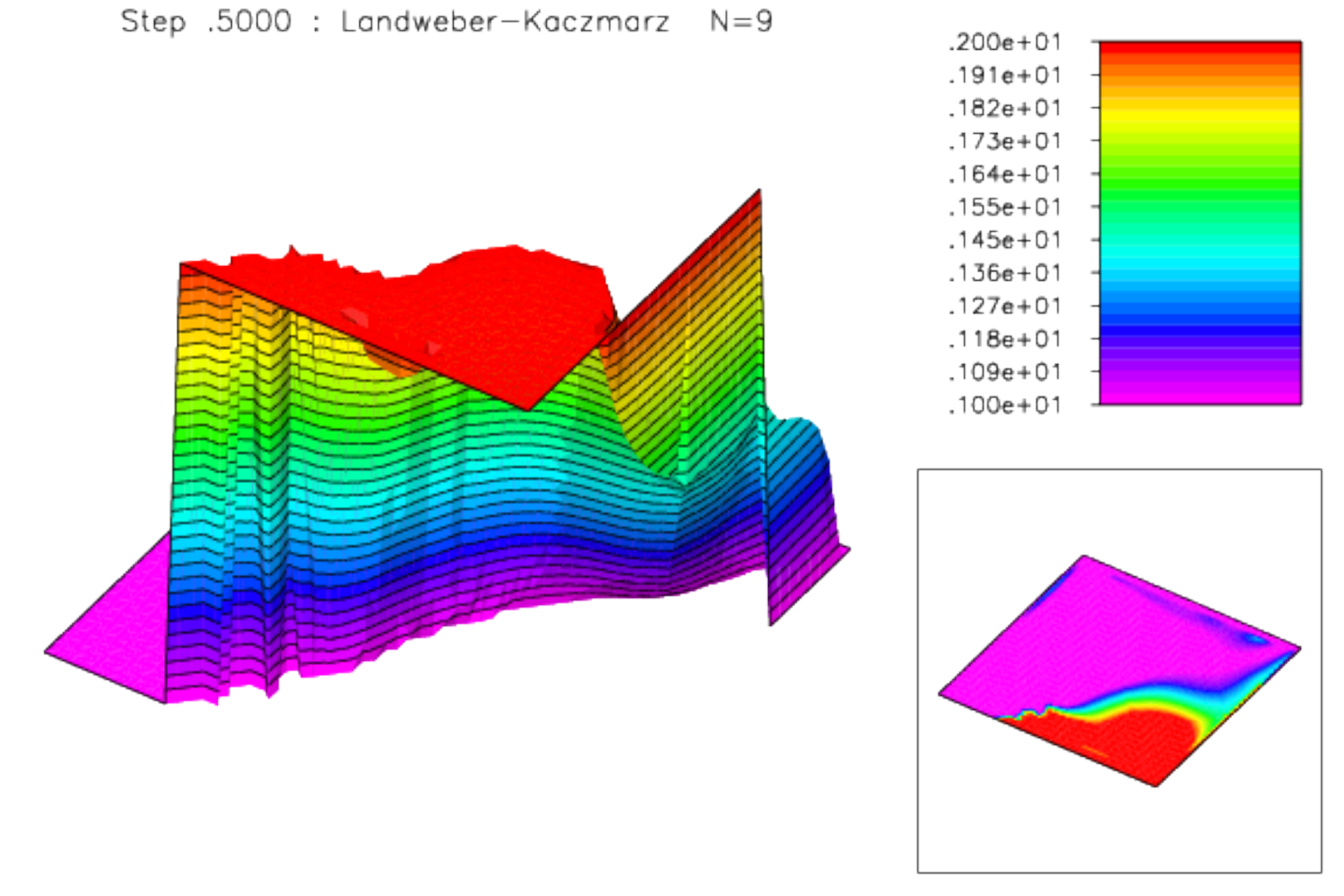}  }
\caption{Evolution of the Landweber-Kaczmarz method for $N=9$ and
exact data.} \label{fig:evol-s9}
\end{figure}

\begin{figure}
\centerline{ \epsfysize4cm \epsfbox{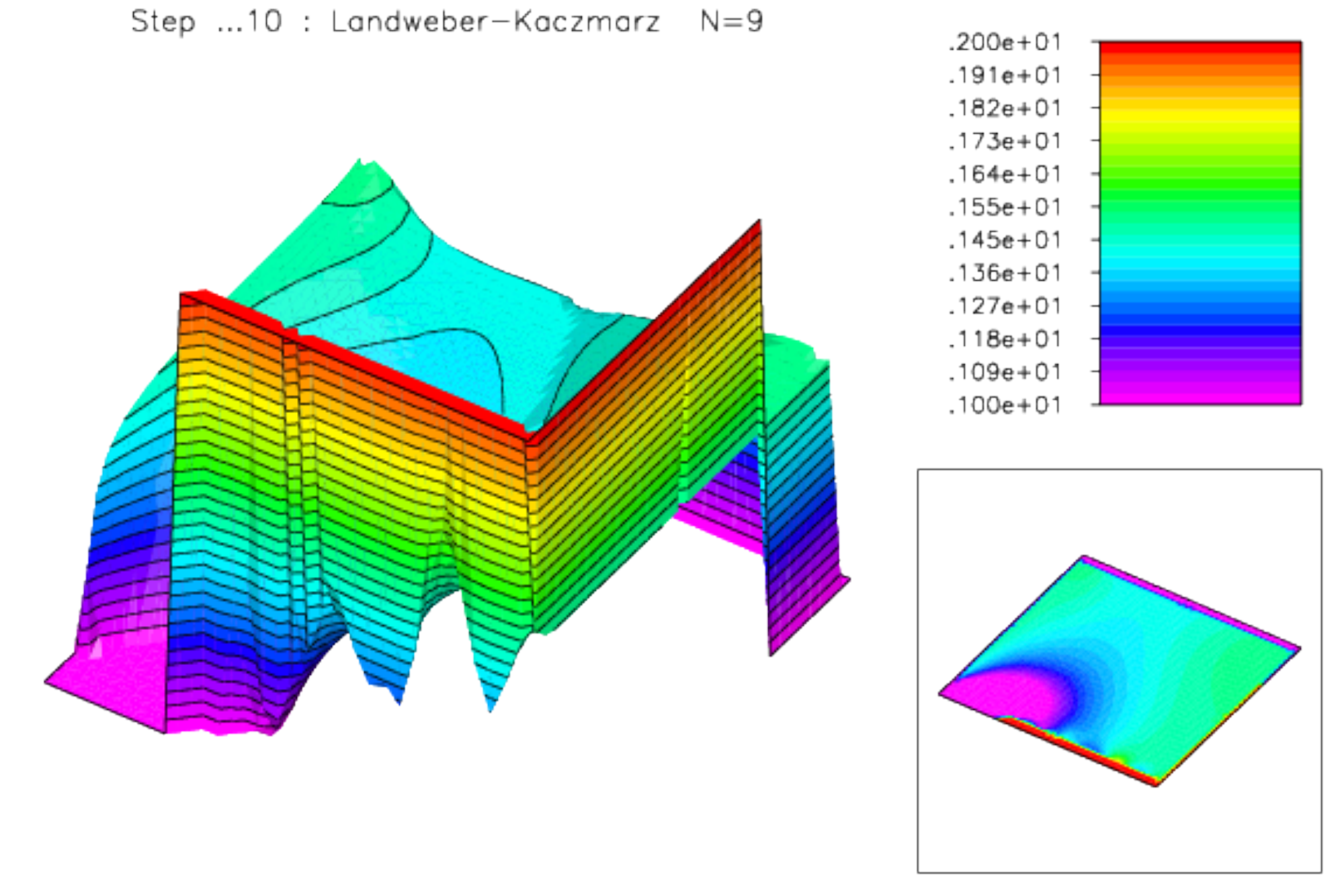} \hfill
             \epsfysize4cm \epsfbox{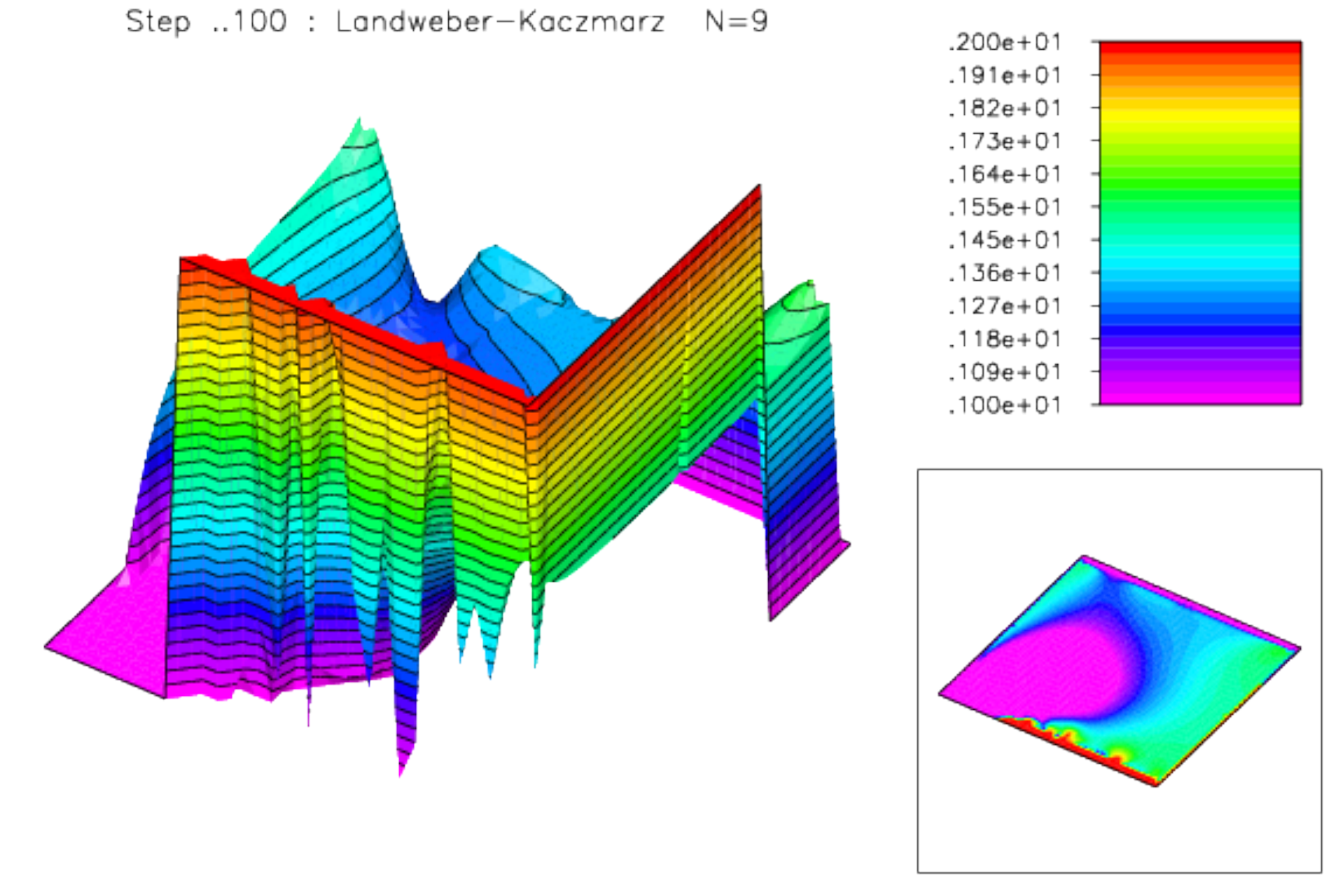}  }
\bigskip
\centerline{ \epsfysize4cm \epsfbox{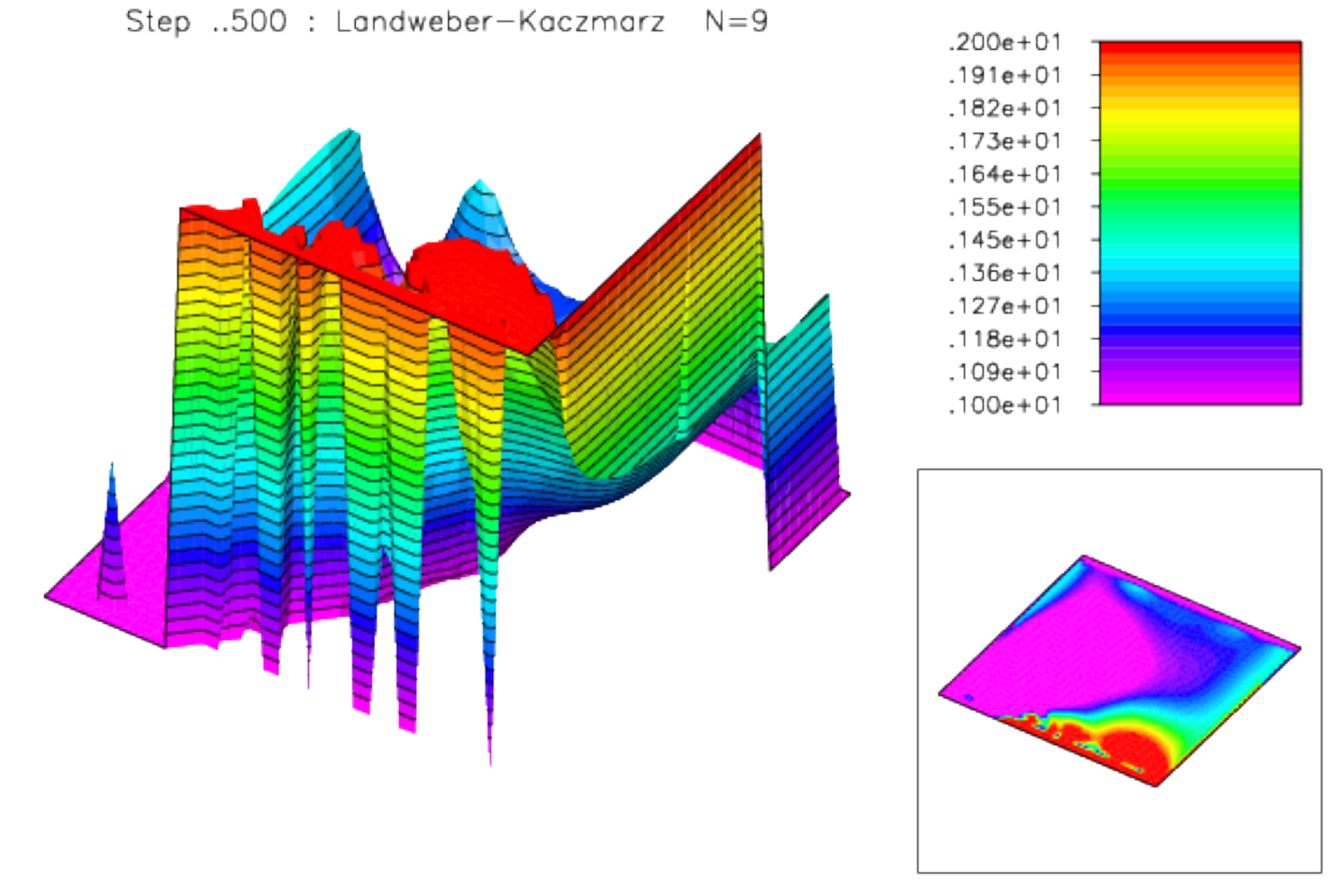} \hfill
             \epsfysize4cm \epsfbox{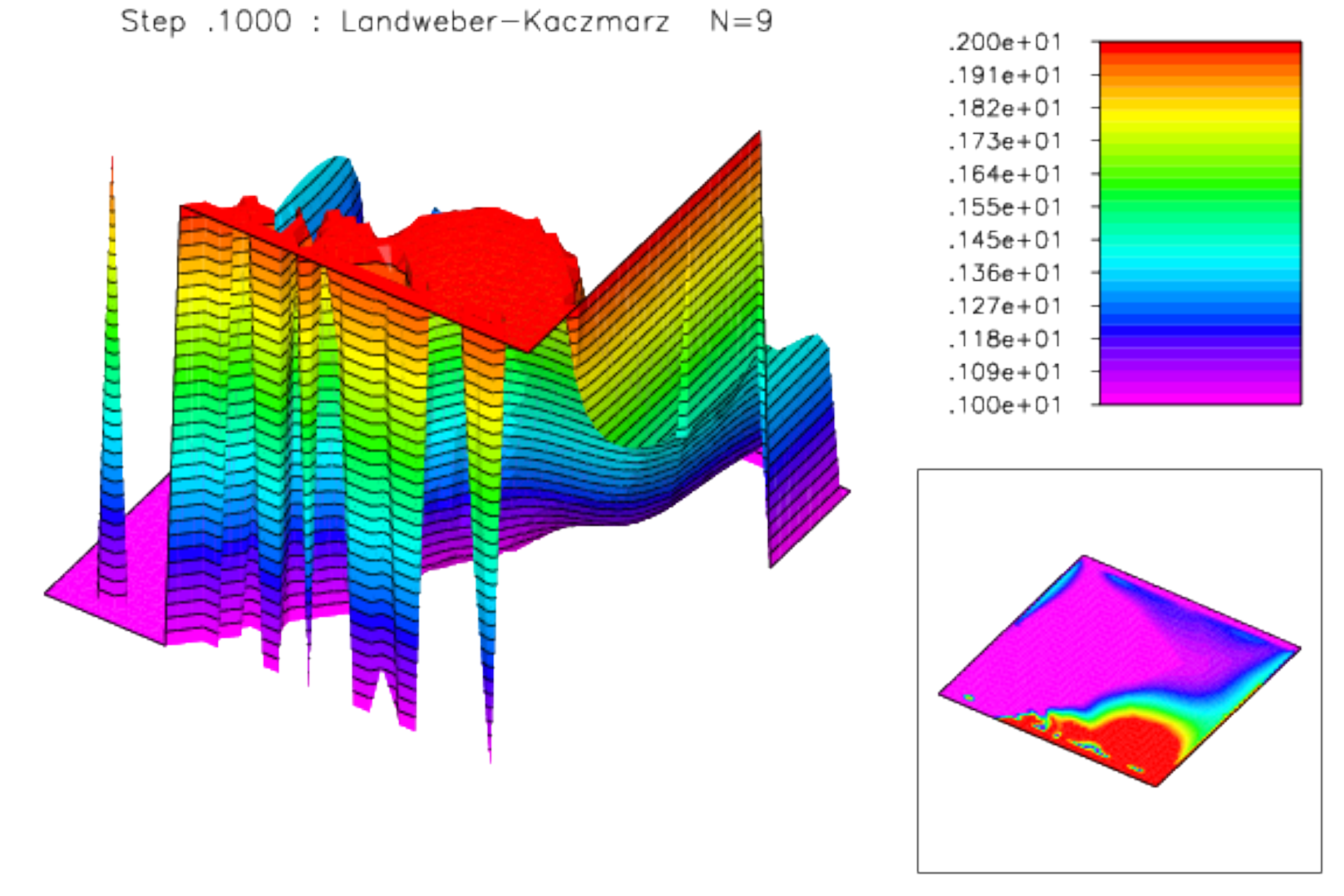}  }
\bigskip
\centerline{ \epsfysize4cm \epsfbox{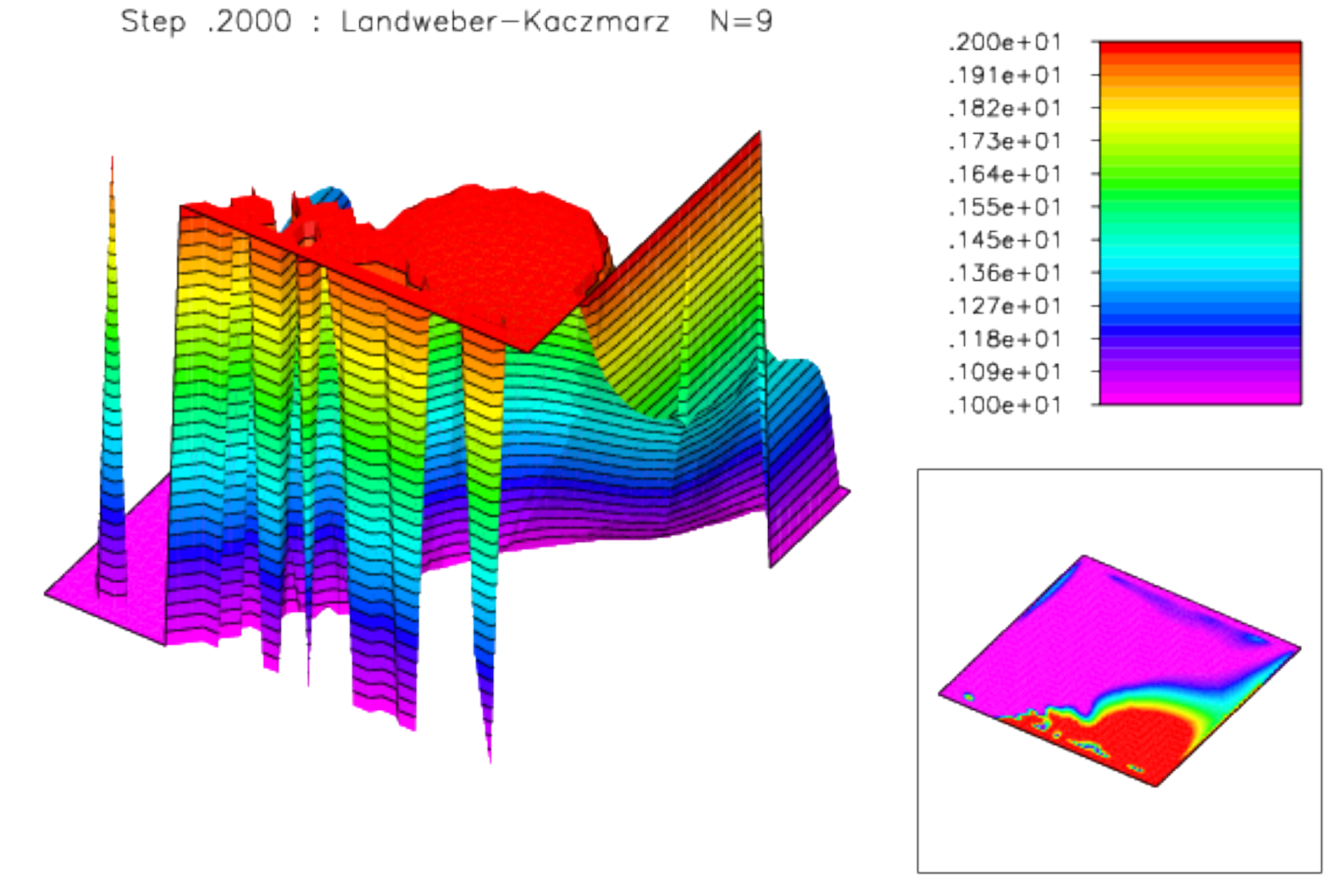} \hfill
             \epsfysize4cm \epsfbox{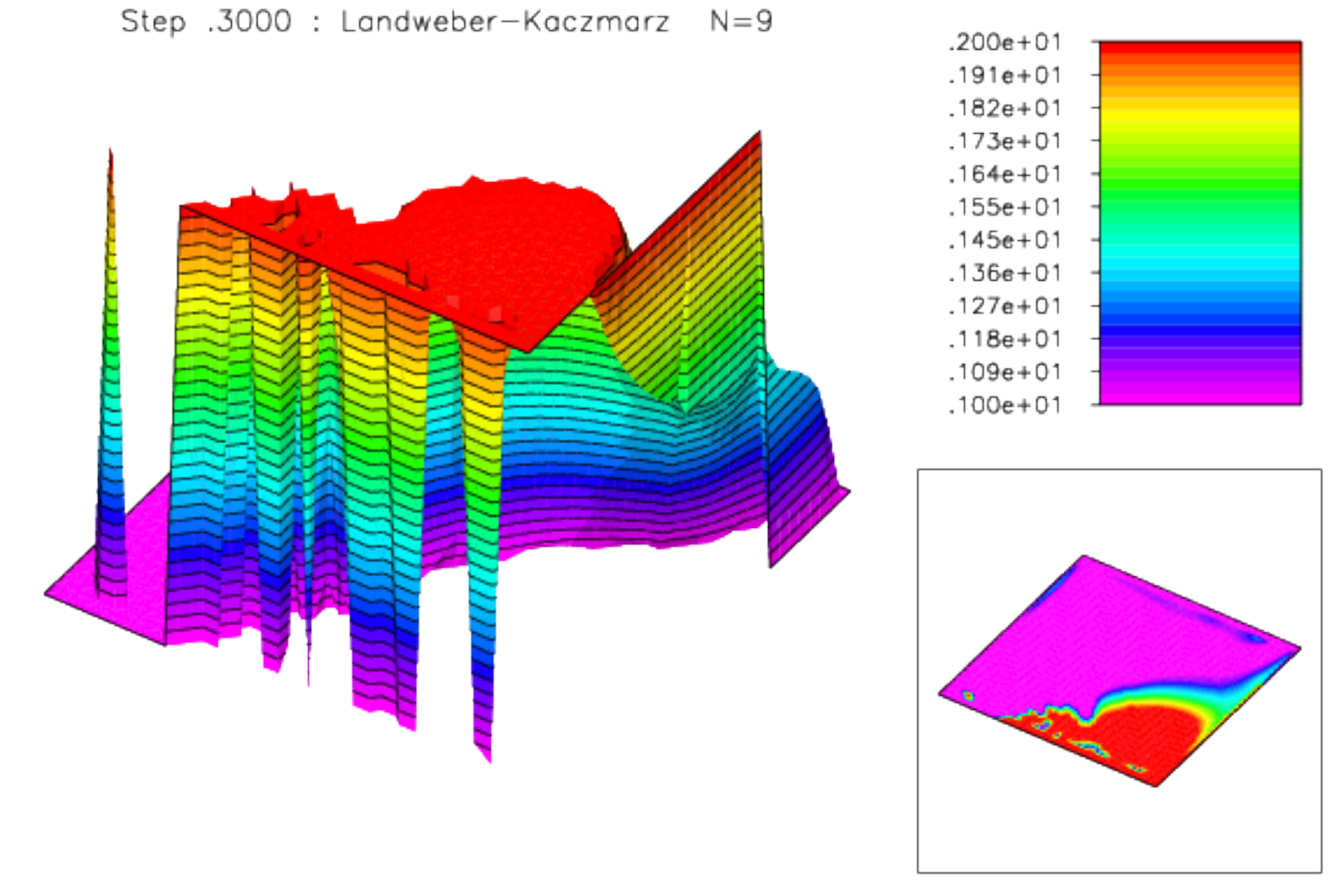}  }
\bigskip
\centerline{ \epsfysize4cm \epsfbox{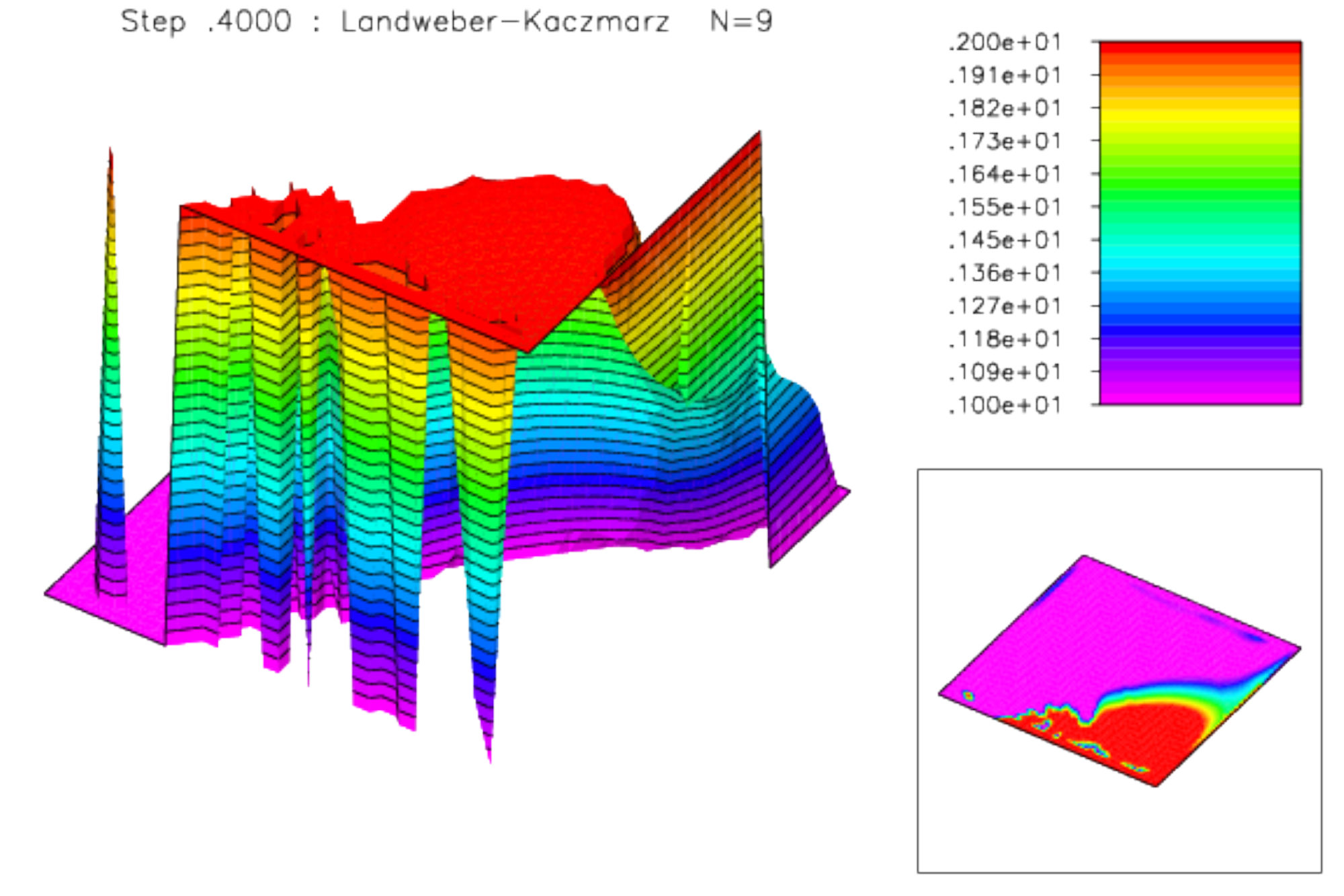} \hfill
             \epsfysize4cm \epsfbox{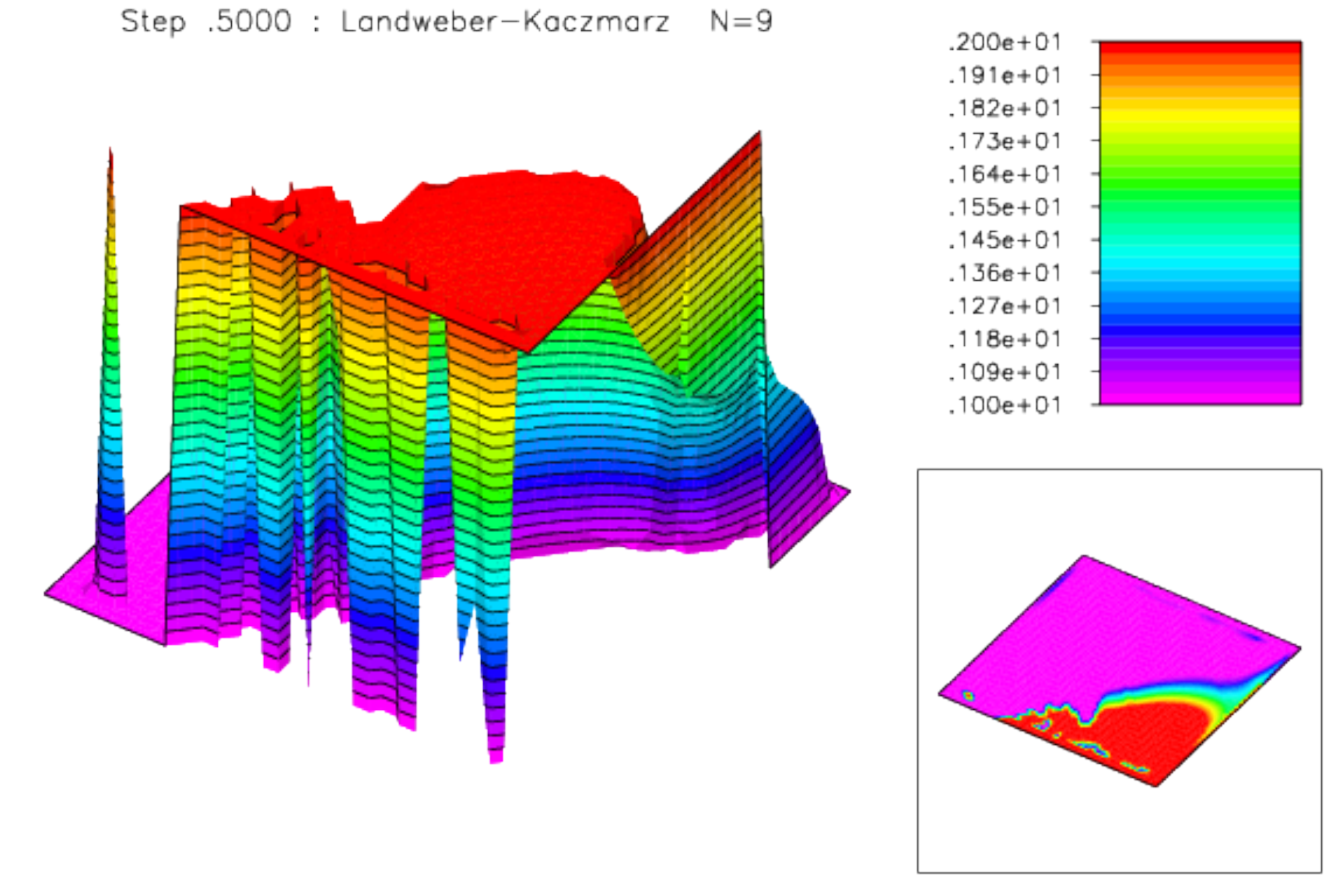}  }
\caption{Evolution of the Landweber-Kaczmarz method for $N=9$ and
noisy data. Noise level of 10\%.}
\label{fig:evol-s9-er}
\end{figure}

In Figure~\ref{fig:evol-s9} we present the reconstruction results obtained
by the Landweber--Kaczmarz iteration for $N=9$, i.e. nine pair of Dirichlet
and Neumann Data. We implemented the method with different amounts of
data (i.e. different values of $N$). For $N \ge 5$ the numerical results
were very close. The results correspond to exact data, i.e. no noise
was introduced. The numerics have shown to be sensible with respect to
noise. Even though, we were able to obtain some acceptable results
for noisy data. In Figure~\ref {fig:evol-s9-er} we present the results
obtained with a noise level of 10\% (white noise).

\section*{Acknowledgment}
M.B. and H.E. acknowledge financial support from the Austrian National
Science Foundation FWF through project SFB F 013/08.
A.L. is on leave from Department of Mathematics, Federal University of
St.\,Cata\-rina, Brazil; his work is supported by the Austrian Academy
of Sciences and CNPq, grant 305823/2003-5.
P.M. acknowledges support from the Austrian National Science Foundation
FWF through his Wittgenstein Award.

%---------------------------------------------------------------------------

\end{document}